\documentclass[12pt]{article}

\newif\ifpaper
\paperfalse % uncomment for arxiv and comment above

\usepackage{myarticle}
\usepackage{mycommands}
\usepackage{subcaption}
\usepackage[utf8]{inputenc} % allow utf-8 input
\usepackage[T1]{fontenc}    % use 8-bit T1 fonts
\usepackage{hyperref}       % hyperlinks
\usepackage{url}            % simple URL typesetting
\usepackage{booktabs}       % professional-quality tables
\usepackage{amsfonts}       % blackboard math symbols
\usepackage{nicefrac}       % compact symbols for 1/2, etc.
\usepackage{microtype}      % microtypography

\newcommand{\reg}{\varphi}

\newcommand{\X}{\R^N}                    % some space
\newcommand{\Y}{\R^M}

\newcommand{\leg}{h}                      % legendre function
\newcommand{\clLeg}{\mathscr L}           % class of Legendre functions
\newcommand{\bregmap}[1][]{D_{#1}}        % Bregman distance (as a mapping without arguments)
\newcommand{\breg}[3][]{D_{#1}(#2,#3)}    % Bregman distance
\newcommand{\dist}{\mathrm{dist}}
\newcommand{\crit}{\mathrm{crit}}
\newcommand{\Prox}[1][]{\ensuremath{P^{#1}}}            % proximal mapping
\renewcommand{\proj}{\mathrm{proj}}

\newcommand{\fun}{f}                     % objective function
\newcommand{\kp}{{k+1}}
\renewcommand{\k}{{k}}
\newcommand{\km}{{k-1}}
\newcommand{\iter}[1]{_{#1}}

\renewcommand{\rto}[1]{\overset{#1}{\to}}
\newcommand{\opt}[1]{#1^\star}

\newcommand{\lb}{\underline{L}} % lower bound constant
\newcommand{\ub}{\bar{L}} % lower bound constant
\newcommand{\setsep}{\,\vert\,}
\newcommand{\KL}{Kurdyka--{\L}ojasiewicz\xspace}

\newcommand{\dimN}{N}

\newcommand{\dimP}{P}

\newcommand{\Z}{\mathbb Z}
\newcommand{\n}{{n}}
\newcommand{\np}{{n+1}}

\renewcommand{\k}{{k}}
\newcommand{\pit}[1]{_{#1}}             % iteration index for scalars
\definecolor{miared}{rgb}{0.7,0.15,0}
\definecolor{miablue}{rgb}{0.058,0.047,.596} 
\definecolor{mygreen}{rgb}{0.058,.596,0.047}

\newcommand{\F}{\mathcal F}
\newcommand{\Fto}[1][]{\overset{#1}{\to}}

\newcommand{\PPP}{\mathcal{P}}

\usepackage[colorinlistoftodos,prependcaption]{todonotes}

\title{\vspace{-2em}
\textbf{Global Convergence of Model Function 
Based \\Bregman Proximal Minimization Algorithms }}

\author{Mahesh Chandra Mukkamala\thanks{Department of Mathematics, University of T\"{u}bingen, Germany, E-mail: \texttt{mamu@math.uni-tuebingen.de}} \quad\quad Jalal Fadili\thanks{Normandie Univ, ENSICAEN, CNRS, GREYC, France, Email: \texttt{jalal.fadili@greyc.ensicaen.fr}}   \quad\quad  Peter Ochs\thanks{Department of Mathematics,  University of T\"{u}bingen,  Germany, E-mail: \texttt{ochs@math.uni-tuebingen.de}}}
\date{}
\begin{document}

\maketitle
\vspace{-2em}
\begin{abstract}
  Lipschitz continuity of the gradient mapping of a continuously differentiable function plays a crucial role in designing various optimization algorithms. However, many functions arising in practical  applications such as low rank matrix factorization or deep neural network problems do not have a Lipschitz continuous gradient. This led to the development of a generalized notion  known as the $L$-smad property, which is based on generalized proximity measures called Bregman distances. However, the $L$-smad property cannot handle nonsmooth functions, for example, simple nonsmooth functions like $\abs{x^4-1}$ and also many practical composite problems are out of scope. We fix this issue by proposing the MAP property, which generalizes the $L$-smad property and is also valid for a large class of nonconvex nonsmooth composite problems. Based on the proposed MAP property, we propose a globally convergent algorithm called Model BPG, that unifies several existing algorithms. The convergence analysis is based on a  new Lyapunov function. We also numerically illustrate the superior performance of Model BPG on standard phase retrieval problems, robust phase retrieval problems, and Poisson linear inverse problems, when compared to a state of the art optimization method that is valid for generic nonconvex nonsmooth optimization problems.
\end{abstract}

\section{Introduction}\label{sec:intro}

  We are interested in solving the following nonconvex optimization problem:
  \begin{equation*}
      (\PPP) \qquad \inf_{\bx \in \R^N} f(\bx),
  \end{equation*}
  where $f:\R^N \to \eR$ is a proper lower semicontinuous function that is lower bounded. Special instances of the above mentioned problem include two broad classes of problems, namely, additive composite problems (Section~\ref{ssec:forward-backward}) and  composite problems (Section~\ref{ssec:general-composite-problems}).
%   \begin{align*}
%   \text{(Forward--Backward / Additive Composite Problems)} \quad & \inf_{\bx \in \R^N} f_0(\bx) + f_1(\bx) ,\\
%   \text{(Prox-Linear / Composite Problems)} \quad & \inf_{\bx \in \R^N} f_0(\bx) + g(F(\bx))\,,
%   \end{align*}
% where $f_0$ is usually nonsmooth, $f_1$ is continuously differentiable  and its gradient  satisfies a condition that generalizes the usually requires Lipschitz continuity assumption,  $g$ is a Lipschitz continuous function, and $F$ is a continuously differentiable function with its Jacobian satisfying a generalization of Lipschitz continuity condition. 
Such problems arise in numerous practical applications such as, quadratic inverse problems \cite{BSTV18}, low-rank matrix factorization problems \cite{MO2019a}, Poisson linear inverse problems \cite{BBT16}, robust denoising problems with nonconvex total variation regularization \cite{MOPS2020}, deep linear neural networks \cite{MWLCO2019}, and many more. 
  \medskip

In this paper, we design an abstract framework for globally convergent algorithms based on suitable approximations of the objective, where the convergence analysis is moreover driven by a requirement on the approximation quality. A classical special case is that of a continuously differentiable $f: \R^N \to \R$, whose gradient mapping is Lipschitz continuous over $\R^N$. For such a function, the following Descent Lemma (cf. Lemma 1.2.3 of  \cite{N1998})
\begin{equation}\label{eq:lipschitz-gradient}
  - \frac{\lb}{2} \vnorm[]{\bx-\bar{\bx}}^2 \leq f(\bx) - f(\bar{\bx}) - \scal{\nabla f(\bar{\bx})}{\bx-\bar{\bx}} \leq \frac{\ub}{2} \vnorm[]{\bx-\bar{\bx}}^2 \,,\quad  \text{ for all }  \bx,\bbx \in \R^N\,,
\end{equation}
describes the approximation quality of the objective $f$ by its linearization $f(\bar{\bx}) + \scal{\nabla f(\bar{\bx})}{\bx-\bar{\bx}}$ in terms of a quadratic error estimate with certain $\lb, \ub >0$. Such inequalities play a crucial role in designing algorithms that are used to minimize $f$. Gradient Descent is one such algorithm, which we focus here. 
% , if the surrogate function is an upper bound of the objective and certain approximation properties hold at the current iterate $\bx\iter\k$. 
We illustrate Gradient Descent in terms of sequential minimization of suitable approximations to the objective, based on the first order Taylor expansion -- the linearization of $f$ around the current iterate $\bx\iter\k \in \R^N$.  Consider the following {\em model function} at the iterate $\bx\iter\k \in \R^N$:
  \begin{equation}\label{eq:intro-model-function-GD}
  f(\bx;\bx\iter\k) := \fun(\bx\iter\k) + \scal{\nabla f(\bx\iter\k)}{\bx-\bx\iter\k} \,,
  \end{equation}
  where $\scal\cdot\cdot$ denotes the standard inner product in the Euclidean vector space $\R^N$ of dimension $N$ and $f(\cdot; \bx\iter\k)$ is the linearization of $f$ around $\bx\iter\k$. Set $\tau>0$. Now, the Gradient Descent update can be written equivalently as follows:
\begin{equation}\label{eq:gradient-descent}
  \bx\iter\kp = \argmin_{\bx\in\R^N}\, \left\{f(\bx;\bx\iter\k) + \frac{1}{2\tau}\vnorm{\bx - \bx\iter\k}^2 \right\}
    \quad\Leftrightarrow\quad
    \bx\iter\kp = \bx\iter\k - \tau \nabla f(\bx\iter\k)  \,.
\end{equation}
Its convergence analysis is essentially based on the Descent Lemma \eqref{eq:lipschitz-gradient}, which we reinterpret as a bound on the linearization error (model approximation error) of $f$. However, obviously \eqref{eq:lipschitz-gradient} imposes a quadratic error bound, which cannot be satisfied in general. For example, functions like $x^4$ or $(x^3+y^3)^2$ or $(1-xy)^2$ do not have a Lipschitz continuous gradient. The same is true in several of the before mentioned practical applications, for example, matrix factorization \cite{MO2019a} and deep linear neural networks \cite{MWLCO2019} problems. 
\medskip

This issue was recently resolved in \cite{BSTV18}, based on the initial work  in \cite{BBT16},  by introducing a generalization of the Lipschitz continuity assumption for the gradient mapping of a function, which was termed the ``$L$-smad property''. In the context of convex optimization, similar notion namely ``relative smoothness'' was proposed in \cite{LFN18}. Such a notion was also independently considered in \cite{BDX2011}, before \cite{LFN18}. However, all these approaches rely on the model function \eqref{eq:intro-model-function-GD}, which is  the linearization of the function. In this paper, we generalize to arbitrary model functions (Definition~\ref{def:model-function}) instead of the linearization of the function.
\medskip

% \paragraph{$L$-smad property \cite[Definition 2.2]{BSTV18}.} 
Now, we briefly recall the ``$L$-smad property''. The main restrictiveness of the Lipschitz continuous gradient notion arises as only quadratic model approximation errors are allowed. Even for simple functions like $x^4$ such quadratic bounds do not exist. Hence, generalized proximity measures which allow for higher order bounds are needed. To this regard, the $L$-smad property relies on generalized proximity measures known as Bregman distances. These distances are generated from so-called Legendre functions (Definition~\ref{def:legendre-function}). Consider a Legendre function $h$, then the Bregman distance between ${\bf x} \in \dom\leg$ and ${\bf y} \in \sint\dom\leg $  is given by 
\begin{equation}
  D_h({\bf x},{\bf y}) := \leg({\bf x}) - \leg({\bf y}) - \scal{{\bf x}-{\bf y}}{\nabla \leg({\bf y})}\,.
\end{equation}
A continuously differentiable function $f: \R^N \to \R$  is $L$-smad with respect to a Legendre function $h: \R^N \to \R$  over $\R^N$ with $\ub,\lb > 0$, if the following condition holds true:
\begin{equation} \label{eq:model-ineq-1}
  -\lb\breg[\leg]{\bx}{\bar{\bx}} \leq  f(\bx) - f(\bar{\bx}) - \scal{\nabla f(\bar{\bx})}{\bx-\bar{\bx}} \leq \ub\breg[\leg]{\bx}{\bar{\bx}}\,, \quad \text{ for any }\bx, \bbx \in \R^N\,.
\end{equation}
We interpret these inequalities as a generalized distance measure for the linearization error of $f$.
% Consider a function $f: \R^N \to \R$ that is continuously differentiable and a Legendre function $h: \R^N \to \R$ with $\dom\leg = \R^N$. Then, $f$ is $L$-smad with respect to a Legendre function $h$  over $\R^N$ with $\ub,\lb > 0$, if the following condition holds true
%   \begin{equation} \label{eq:model-ineq-1}
% -\lb\breg[\leg]{\bx}{\bar{\bx}} \leq  f(\bx) - f(\bar{\bx}) - \scal{\nabla f(\bar{\bx})}{\bx-\bar{\bx}} \leq \ub\breg[\leg]{\bx}{\bar{\bx}}\,, \quad \text{ for any }\bx, \bbx \in \R^N\,.
%   \end{equation}
Similar to the Gradient Descent setting, minimization of  $f(\bar{\bx}) + \scal{\nabla f(\bar{\bx})}{\bx-\bar{\bx}}  +  \frac{1}{\tau}\breg[\leg]{\bx}{\bar{\bx}}$ essentially results in the Bregman proximal gradient (BPG) algorithm's update step \cite{BSTV18}. 
\medskip

However, the $L$-smad property relies on the continuous differentiability of the function $f$, thus nonsmooth functions as simple as $\abs{x^4-1}$ or $\abs{1-(xy)^2}$ or $\log(1+ \abs{1-(xy)^2})$ cannot be captured under the $L$-smad property.  Numerous difficult nonsmooth optimization problems cannot be captured either. This motivates a more general notion than the $L$-smad property. 
\medskip

This lead us to the development of the MAP property (Definition~\ref{def:map-property}), where MAP abbreviates Model Approximation Property.  Consider a function $f:\R^N \to \R$ that is proper lower semicontinuous, and a Legendre function $h : \R^N \to \R$ with $\dom \leg  = \R^N$. We abbreviate ``lower semicontinuous''  as  ``lsc''. For certain $\bbx \in \R^N$, we consider generic model function $f(\bx;\bbx)$ that is proper lsc and approximates the function around the model center $\bbx$, while preserving the local first order information (Definition~\ref{def:model-function}). The MAP property is satisfied with the constants $\bar L >0$ and $\underline{L} \in \R$ if for any $\bbx\in \R^N$ the following holds: 
\begin{equation} \label{eq:model-ineq}
-\lb\breg[\leg]{\bx}{\bar{\bx}} \leq  \fun(\bx)- \fun(\bx;\bbx) \leq \ub\breg[\leg]{\bx}{\bar{\bx}} \,, \quad \forall \bx\,\in\,\R^N\,.
\end{equation}
Note that we do not require the continuous differentiability of the function $f$. Our MAP property is inspired from \cite{DDJ2018}, however, their work considers only the upper bound, and also they rely on decomposition of function into two components.
\medskip

We illustrate the MAP property with a simple example. Consider a composite problem  $f(x) = g(F(x)) := \abs{x^4-1}$, where $F(x) := x^4-1$ is a continuously differentiable function over $\R$, and $g(x) := \abs{x}$ is a Lipschitz continuous function  over $\R$.   Note that neither the Lipschitz continuity of the gradient nor the $L$-smad property is valid for this problem. However, the MAP property is valid here. At certain  ${\bar x} \in \R$, we consider the model function that is given by $f(x; \bar{x}) := g(F(\bar{x}) + \nabla F(\bar{x})(x-\bar{x}))$, where $\nabla F(\bar{x})$ is the Jacobian of $F$ at ${\bar x}$. Then, with $\ub = \lb = 4$, the MAP property is satisfied:
\begin{equation} \label{eq:model-ineq-2a}
  -\lb\breg[\leg]{x}{\bar{x}} \leq  g(F(x)) - g(F(\bar{x}) + \nabla F(\bar{x})(x-\bar{x})) \leq \ub\breg[\leg]{x}{\bar{x}} \,, \text{ for all }x,{\bar x} \in \R\,,
\end{equation}
where $h(x) = 0.25x^4$ and the generated Bregman distance is $\breg[\leg]{x}{\bar{x}} = 0.25x^4 - 0.25{\bar x}^4 - {\bar x}^3(x - {\bar x})$.  We  provide further details in Example~\ref{ex:running-example} and in Example~\ref{ex:running-example-1}. 
\medskip

We considered the above given composite problem for illustration purposes, and we emphasize that our framework is applicable for large classes of nonconvex problems (see Section~\ref{sec:examples}).  Similar to the BPG setting, minimization of $f(\bx ; \bbx)  +  \frac{1}{\tau}\breg[\leg]{\bx}{\bar{\bx}}$ essentially results in Model BPG algorithm's update step. The precise definition of the model function is provided in Definition~\ref{def:model-function}, the MAP property in full generality is provided in Definition~\ref{def:map-property}, and the Model BPG algorithm is provided in Algorithm~\ref{alg:acc-BregMin-bt}.
\medskip

We now discuss our main contributions and the related work. 

\subsection{Contributions} Our main contributions are the following.
\begin{itemize}[leftmargin=*]
  \item We introduce the MAP property, which generalizes the Lipschitz continuity assumption of the gradient mapping and the $L$-smad property \cite{BSTV18, BBT16}. Earlier proposed notions were restricted to additive composite problems. The MAP property is essentially an extended Descent Lemma that is valid for generic composite problems (see Section~\ref{sec:examples}), based on Bregman distances. Our theory is applicable to generic nonconvex nonsmooth objectives, and is not restricted to composite objectives. MAP like property was also partially considered in \cite{DDJ2018}, however with focus on stochastic optimization. The MAP property relies on the notion of \textit{model function}, that serves as a function approximation, and preserves the local first order information of the function.  Our work extends the foundations laid by \cite{drusvyatskiy2019nonsmooth,  DDJ2018} that consider generic model functions (potentially nonconvex), and \cite{OFB19} which considers convex model functions. 
  \item Based on the MAP property, Model based Bregman Proximal Gradient (Model BPG) algorithm (Algorithm~\ref{alg:acc-BregMin-bt}) is proposed. Several existing algorithms such as Proximal Gradient Method \cite{CP11a}, Bregman Proximal Gradient Method \cite{BSTV18} (or Mirror Descent \cite{BT03}), Prox-Linear algorithm \cite{DP2019}, and many other algorithms can be seen as a special case. Moreover, novel  algorithms arise depending on the definition of the model function. We emphasize that Model BPG is practical, simple to implement and also does not require special knowledge about the problem such as the so-called information zone  \cite{BST2018}. Close variants of Model BPG already exist in the literature, such as  line search based Bregman proximal gradient method \cite{OFB19}, and mirror descent variant \cite{DDJ2018}, however, the convergence of the full sequence of iterates was not known.
  \item The standard global convergence analysis, in the sense that the full sequence of iterates converges to a single point, relies on descent properties of function values evaluated at the iterates of an algorithm. However, using function values can be restrictive, and alternatives are sought \cite{P2016}. To fix this issue, we introduce a new  Lyapunov function, through which we prove the global convergence of the full sequence of iterates generated by Model BPG. We eventually show that the sequence generated by Model BPG converges to a critical point of the objective function, which is potentially nonconvex and nonsmooth. Notably, the usage of a Lyapunov function is popular for inertial algorithms \cite{OCBP14, MOPS2020} and through our work we aim to popularize Lyapunov functions also for noninertial algorithms. Usage of Lyapunov functions is also popular in the context of dynamical systems \cite{HC2011}.
  \item The global convergence analysis of Bregman proximal gradient (BPG) \cite{BSTV18} relies on the full domain of the Bregman distance. However, there are many Bregman distances for which the domain is restricted. We show in this paper, that under certain assumptions that are typically satisfied in practice, the global convergence of the full sequence of iterates generated by Model BPG using generic Bregman distances can indeed be obtained (Theorem~\ref{thm:full-conv}, \ref{thm:global-conv-obj-func}). In general, this requires the limit points of the sequence to lie in the interior of domain of the employed Legendre function. While this is certainly a restriction, nevertheless, the considered setting is highly nontrivial and novel in the general context of nonconvex nonsmooth optimization. Moreover, it allows us to avoid the common restriction of requiring (global) strong convexity of the Legendre function, which is a severe drawback that rules out many interesting applications in related approaches (e.g., see Section~\ref{ssec:poisson-linear-inverse-problems}). 
  \item We provide a comprehensive numerical section showing the superior performance of Model BPG compared to a state of the art optimization algorithm, namely, Inexact Bregman Proximal Minimization Line Search (IBPM-LS) \cite{ODBP13}, on standard phase retrieval problems, robust phase retrieval problems and Poisson linear inverse problems.
\end{itemize}

\subsection{Related work} 
Our work is fundamentally based on three pillars, namely, Bregman distances, model functions, and \KL (KL) inequality. Bregman distances are certain generalized proximity measures, which generalize Euclidean distances. Model functions serve as function approximations which preserve local first order information about the function. The KL inequality is a certain regularity property of the function, which is crucial for global convergence analysis, and is typically satisfied by objectives that arise in practice. We provide below the related work based on these three topics. 
\paragraph{Bregman distances.} 
Recently, there has been huge surge of work on Bregman distances \cite{M18,JN11,B15,benning2017choose,corona2019enhancing,geiping2018composite,pang2018decomposition,eskandani2018hybrid}.  This is due to the flexibility one gains in modelling the proximity measures. The seminal Mirror Descent algorithm \cite{BT03} incorporates Bregman distances in the update step. At-times the special structure of the Bregman distance can result in closed form update steps, simple case being Gradient Descent with Euclidean Distance. Also, for instance in the minimization problem obtained for deblurring an image under Poisson noise, one can obtain a closed form expression for an optimization subproblem using a Bregman distance generated by Burg's entropy \cite{OFB19}. Bregman distances for structured matrix factorization problems were considered in \cite{MO2019a, TV2020, DBA2019, LZTW2019, HG2020} and certain extensions to deep linear neural networks were considered in \cite{MWLCO2019}. Bregman distances allow for many optimization algorithms, which were previously thought to be completely different to co-exist in a single algorithm, thus making the analysis simpler. The crucial observation that Bregman distances can indeed be used to generalize the notion of Lipschitz continuous gradient was considered in \cite{BBT16}. However, their setting was restricted to convex problems. This was later mitigated in \cite{BSTV18}, via the $L$-smad property for nonconvex problems. Recently, the related notions such as relative smoothness \cite{LFN18}, and relative continuity \cite{L17} were proposed based on Bregman distances. Before \cite{BBT16} and \cite{LFN18}, the work in \cite{BDX2011} also considered a generalization of the Lipschitz continuous gradient notion. More related references also include \cite{Nguyen17, LOC2019}. As mentioned in the introduction, the $L$-smad property can also be restrictive, and thus we propose the MAP property to generalize the $L$-smad property even further. Closely related work is \cite{DDJ2018}, however, their focus was on developing stochastic algorithms.

% Most algorithms can be interpreted as sequential minimization of surrogate functions. While majorization--minimization methods employ an upper bound as a surrogate function, the key is a certain approximation quality of the model function to the objective.

\paragraph{Model functions.} The MAP property relies on the concept of the model function, which is essentially a function approximation that preserves the local first order information. In smooth optimization, it is common to use the Taylor approximation of a certain order as model function. In nonsmooth optimization, we can only speak of ``Taylor-like'' models \cite{NPA08,Noll13,drusvyatskiy2019nonsmooth,OFB19}, which is a (nonunique) approximation that satisfies certain error bound or a growth function \cite{drusvyatskiy2019nonsmooth,OFB19}. The class of model functions used in \cite{NPA08,Noll13} only satisfy a lower bound, and bundle methods are developed, which is a different class of algorithms that we do not discuss here. The growth functions in \cite{drusvyatskiy2019nonsmooth,OFB19} that measure the approximation quality of the model function, which is also used in this paper, can be interpreted as a generalized first-order oracle. It has been shown in \cite{OFB19} that the \textit{concept of model functions unifies several algorithms} for smooth and nonsmooth optimization, for example, Gradient Descent, Proximal Gradient Descent, Levenberg Marquardt's method, ProxDescent, certain variable metric versions of these algorithms and some related majorization--minimization based algorithms. More recently, model functions were considered in the context of the Conditional Gradient method in \cite{OM2019}. A particularly interesting class of model functions is the one for which the approximation quality measure is formed by Bregman distances \cite{BBT16, BSTV18,OFB19}, which is our main focus in this paper.

\paragraph{\KL inequality.} Based on the MAP property, we propose Model BPG algorithm. In order to prove the global convergence of the full sequence of iterates generated by Model BPG algorithm, the \KL inequality \cite{Kurd98,Loj63,Loj93,BDL06,BDLS07} is key. This inequality is satisfied by most functions that appear in practical applications, in particular, semi-algebraic functions \cite{BCR98}, globally subanalytic functions \cite{BDL06b}, or more generally, functions that are definable in an o-minimal structure \cite{BDLS07,Dries98}. Usually, the essential conditions required for global convergence analysis can be collected in an abstract manner, and are clearly summarized and studied in \cite{ABS13,BST14}. Basically, the conditions that need to be verified are called ``sufficient descent condition'', ``relative error condition'', and ``continuity condition''. The sequence satisfying such conditions is at-times called \textit{gradient-like descent sequence} \cite{BSTV18}, which we detail in Section~\ref{sec:gradient-like-descent} in the appendix. In order to prove the global convergence of the full sequence of iterates generated by Model BPG, it suffices to prove that it is a gradient-like descent sequence. Sequences arising using several algorithms such as Bregman Proximal Gradient (BPG) or Proximal Gradient method are gradient-like descent sequences. In the context of additive composite problems, global convergence analysis of BPG was provided in \cite{BSTV18}. However, their setting is restrictive as the employed Legendre function is assumed to be strongly convex with full domain and the model framework is not considered.  In this paper, we do not have such restrictions, thus our framework is highly general and is applicable to broad classes of nonconvex nonsmooth problems (see Section~\ref{sec:examples}).

% \subsection{Preliminaries and Notations} \label{sec:prelim}
\subsection{Preliminaries and notations.} We work in a Euclidean vector space $\X$ of dimension $N\in\N$ equipped with the standard \emphdef{inner product} $\scal\cdot\cdot$ and induced \emphdef{norm} $\vnorm\cdot$. For a set $C\subset\X$, we define $\vnorm[-]{C}:=\inf_{{\bf s}\in C}\, \vnorm{{\bf s}}$. We skip basic definitions here, instead we provide them in Section~\ref{sec:additional-preliminaries} in the appendix and all notations are primarily taken from \cite{Rock98}.
\medskip

Legendre functions defined below generate the Bregman distances, which are generalized proximity measures compared to the Euclidean distance.

\begin{DEF}[{Legendre function \cite[Def. 1]{BBT16}}]\label{def:legendre-function}
  Let $h: \R^N  \to \eR$ be a proper lsc 
  convex function. It is called: (i) essentially smooth, if $h$ is differentiable on $\sint\dom\leg$, with moreover $\vnorm[]{\nabla h({\bf x}\iter\k)} \to \infty$ for every sequence $\seq[\k\in\N]{{\bf x}_\k} \in \sint\dom\leg$ converging to a boundary point of $\dom\leg$ as $k\to\infty$; (ii) of Legendre type if $\leg$ is essentially smooth and strictly convex on $\sint\dom\leg$.
\end{DEF}

Some properties of Legendre function include the following: 
\[
  \dom \partial h = \sint\dom\leg, \text{ and }\, \partial h(\bx) = \{\nabla h(\bx)\},\, \forall \bx \in \sint\dom\leg.
\]
Legendre function is also referred as kernel generating distance \cite{BSTV18}, or a reference function \cite{LFN18}. Generic reference functions used in  \cite{LFN18} are more general compared to Legendre functions, as they do not require essential smoothness. 
\medskip

The Bregman distance associated with any Legendre function $h$ is defined by
\begin{equation}
  D_h({\bf x},{\bf y}) = \leg({\bf x}) - \leg({\bf y}) - \scal{{\bf x}-{\bf y}}{\nabla \leg({\bf y})}, \quad \forall\, {\bf x} \in \dom\leg,\, {\bf y} \in \sint\dom\leg \,.
\end{equation}
In contrast to the Euclidean distance, the Bregman distance is lacking symmetry. 
\paragraph{Examples.} Prominent examples of Bregman distances can be found in \cite[Example 1, 2]{BBT16}. We provide some examples below. For any vector $\bx \in \R^N$, the $i^{\,\text{th}}$ coordinate is denoted by $\bx_i$.
\begin{itemize}
  \item Bregman distance generated from $h({\bf x}) = \frac12\vnorm[]{{\bf x}}^2$ is equivalent to the Euclidean distance. 
  \item Let $\bx, \bbx \in \R_{++}^N$, for  $h({\bf x}) = -\sum_{i=1}^N\log({\bf x}_i)$ (Burg's entropy), the generated Bregman distance is  
  \[
    D_h({\bf x},{\bf {\bar x}}) = \sum_{i=1}^N\left(\frac{{\bf x}_i}{{\bf {\bar x}}_i} - \log\left(\frac{{\bf x}_i}{{\bf {\bar x}}_i}\right) - 1\right)\,.
  \]
  Such distances are helpful in Poisson linear inverse problems \cite{BBT16,OFB19}. 
  \item   Let $\bx \in \R_{+}^N$, $\bbx \in \R_{++}^N$,  for  $h({\bf x}) = \sum_{i=1}^N{\bf x}_i\log({\bf x}_i)$ (Boltzmann--Shannon entropy), with $0\log(0):= 0$, the Bregman distance is given by 
  \[
    D_h({\bf x},{\bf {\bar x}}) = \sum_{i=1}^N{\bf x}_i (\log({\bf x}_i)-\log({\bf {\bar x}}_i)) - ({\bf x}_i-{\bf {\bar x}}_i)\,.
  \] 
  Such distances are helpful to handle simplex constraints \cite{BT03}. 
  \item Phase retrieval problems \cite{BSTV18} use the Bregman distance based on the Legendre function $h : \R^N \to \R$ that is given by
  \[
    h({\bf x})= \frac{1}{4}\vnorm[2]{{\bf x}}^4 + \frac{1}{2}\vnorm[2]{{\bf x}}^2\,.
  \]
  \item  Matrix factorization problems \cite{MO2019a, TV2020} use the Bregman distance based on the Legendre function $h : \R^{N_1} \times \R^{N_2} \to \R$ that is given by
  \[
  h(\bx, \by)  = c_1\left(\frac{\vnorm[]{\bx}^2 + \vnorm[]{\by}^2}{2}\right)^2 + c_2\left(\frac{\vnorm[]{\bx}^2 + \vnorm[]{\by}^2}{2}\right)\,,
  \]
  with certain $c_1,c_2>0$ and $N_1,N_2 \in \N$. Based on the work in \cite{MO2019a}, related Bregman distances for deep linear neural networks were also explored later in \cite{MWLCO2019}.
\end{itemize}

\section{Problem setting and Model BPG algorithm}\label{sec:problem}

We solve possibly nonsmooth and nonconvex optimization problems of the form
\begin{equation} \label{eq:problem}
(\PPP) \quad\quad    \inf_{\bx\in \X}\, \fun(\bx)\,,
\end{equation}
that satisfy the following assumption, which we impose henceforth. 
\begin{ASS}\label{ass:problem} The objective function $\map{\fun}{\X}{\eR}$ is proper, lower semi-continuous (possibly nonconvex nonsmooth) and a coercive function, i.e., as $\norm{\bx} \to \infty$ we have $f(\bx) \to \infty$.  
% \begin{itemize}
% \item is proper, lower semi-continuous (possibly nonconvex nonsmooth),  
% \item and satisfies $\Argmin_{\bx\in \X} \fun(\bx) \neq\emptyset$.
% \end{itemize}
\end{ASS}
Due to \cite[Theorem 1.9]{Rock98}, the function $f$ satisfying Assumption~\ref{ass:problem} is bounded from below, and $\Argmin_{\bx\in \X} \fun(\bx)$ is nonempty and compact. We denote the following:
\[
  v(\PPP) := \inf_{\bx \in \R^N} f(\bx) > -\infty\,.
\]
We denote the set of critical points with respect to the limiting subdifferential (see Appendix~\ref{sec:additional-preliminaries}) as 
\[
  \crit f := \left\{ \bx \in \R^N : \;  {\bf 0} \in \partial f(\bx)  \right\}\,.
\]

We require the following technical definitions.
\begin{DEF}[{Growth function \cite{drusvyatskiy2019nonsmooth,OFB19}}] \label{def:growth-function} 
  A differentiable univariate function $\map{\omega}{\R_+}{\R_+}$ is called \emphdef{growth function} if it satisfies $\omega(0)=\omega_+^\prime(0) = 0$,  where $\omega^\prime_+$ denotes the one sided (right) derivative of $\omega$. If, in addition, $\omega_+^\prime(t) >0$ for $t>0$ and\xspace equalities $\lim_{t\dto 0} \omega_+^\prime(t) = \lim_{t\dto 0} \omega(t)/\omega_+^\prime(t) = 0$ hold, we say that $\omega$ is a \emphdef{proper growth function}.
\end{DEF}
Example of a proper growth function is  $\omega(t) = \frac{\eta}{r}t^r$ for  $\eta,r >0$. Lipschitz continuity and H\"{o}lder continuity can be interpreted with growth functions or, more generally, with uniform continuity \cite{OFB19}. We use the notion of a growth function to quantify the difference between a  model function (defined below) and the objective function.

\begin{DEF}[Model Function]\label{def:model-function}
  Let $f$ be a proper lower semi-continuous (lsc) function. A function $\map{f(\cdot,\bbx)}{\R^N}{\eR}$ with $\dom f(\cdot,\bbx) = \dom\fun$ is called \emphdef{model function} for $f$ around the \emphdef{model center} $\bbx \in \dom f$, if there exists a growth function $\omega_{\bbx}$ such that the following is satisfied:
  \begin{equation}\label{eq:model-function-main}
    \abs{f(\bx) - f(\bx;\bbx)} \leq \omega_{\bbx}(\vnorm[]{\bx - \bbx})\,,\quad \forall\, \bx\in \dom\fun.
  \end{equation}
\end{DEF}
Model function is essentially a first-order approximation to a function $f$ (see Lemma~\ref{lem:first-order-info}), which explains the naming as "Taylor-like model" by \cite{drusvyatskiy2019nonsmooth}. The qualitative approximation property is represented by the growth function. We refer to \eqref{eq:model-function-main} as a bound on the model error, and the symbol $\omega_{\bbx}$ denotes the dependency of the growth function on the model center $\bbx$.
\medskip

Few remarks are in order, which we provide below:
\begin{itemize}
  \item  Informally, the model function approximates the function well near the model center. Convex model functions are explored in \cite{OFB19, OM2019}, however in our setting, the model functions can be nonconvex. 
  \item Nonconvex model functions were considered in \cite{drusvyatskiy2019nonsmooth}, however only subsequential convergence was shown. Their work is focussed on the termination criterion of the algorithms, however, they do not present an implementable algorithm. 
\end{itemize}

If the growth function constants are independent of $\bbx$, this results in a uniform approximation. However, typically the growth function depends on the model center, as we illustrate below.
\begin{EX}[Running Example]\label{ex:running-example}
  Let $f(\bx) = \abs{g(\bx)}$ with $g(\bx) = \vnorm[]{\bx}^4 - 1$. With $\bbx \in \R^N$ as the model center, we consider the following model function:
  \[
    f(\bx;\bbx) := \abs{g(\bbx) + \scal{\nabla g(\bbx)}{\bx - \bbx}}\,.
  \]
  As per the proof provided in Section~\ref{sec:proof-running-example} in the appendix, the model error  is given by
  \begin{align*}
    \abs{f(\bx) - f(\bx;\bbx)} &\leq 24\vnorm[]{\bbx}^2\vnorm[]{\bx- \bbx}^2 + 8\vnorm[]{\bx- \bbx}^4\,,
  \end{align*}
  where the growth function is $\omega_{\bbx}(t) = 24\vnorm[]{\bbx}^2 t^2 + 8t^4$.
\end{EX}

The above example illustrates that a constant in the growth function $\omega_{\bbx}(t)$ is dependent on the model center. It is often of interest to obtain a uniform approximation for the model error $\abs{f(\bx) - f(\bx;\bbx)}$, where the growth function is not dependent on the model center. In general, obtaining such a uniform approximation is not trivial, and may even be impossible. Moreover, typically finding an appropriate growth function is not trivial.
\medskip

For this purpose, it is preferable to have a global bound on the model error, for which such a bound can be easily verified, the dependency on the model center is more structured, and the constants arising do not have any dependency on the model center. In the context of additive composite problems, previous works such as \cite{BBT16, LFN18, BSTV18} relied on Bregman distances to upper bound the model error and verified the model error property with a simple convexity test based on second order information (c.f. \cite[Proposition 1]{BBT16}).  Based on this idea, we propose the following MAP property, which is valid for a huge class of generic nonconvex problems and also generalizes the previous works. We emphasize that the MAP property is valid for a large class of nonsmooth functions. MAP like property that is valid for composite problems was also explored in \cite{DDJ2018}. We provide the precise connections to previous works and examples in Section~\ref{sec:examples}.

\begin{DEF}[MAP: Model Approximation Property]\label{def:map-property}
  Let $h$ be a Legendre function that is continuously differentiable over $\sint\dom\leg$. A proper lsc function $f$ with $\dom f \subset \scl\dom h$ and $\dom\fun \cap\, \sint\dom\leg \neq \emptyset$, and model function $f(\cdot, \bbx)$ for $f$ around $\bbx\in \dom f \cap \sint\dom h$ satisfies the \emphdef{Model Approximation Property (MAP) at $\bbx$}, with the constants $\bar L >0$,\,$\underline{L} \in \R$, if for any $\bbx\in \dom f\cap\sint\dom h$ the following holds: 
  \begin{equation} \label{eq:model-ineq}
-\lb\breg[\leg]{\bx}{\bar{\bx}} \leq  \fun(\bx)- \fun(\bx;\bbx) \leq \ub\breg[\leg]{\bx}{\bar{\bx}} \,, \quad \forall \bx\,\in\,\dom\fun \cap \dom\leg\,.
  \end{equation}
\end{DEF}
\begin{REM} We provide the following remarks.
  \begin{itemize}
    \item The design of a model function is independent of an algorithm. However, algorithms can be governed by the model function, for example, Model BPG in Algorithm~\ref{alg:acc-BregMin-bt}. The property of a model function is rather an analogue to differentiability or a (uniform) first-order approximation. Note that for ${\bar{\bx}} \in \sint\dom h$, the Bregman distance $D_h(\bx,{\bar{\bx}})$ is bounded  by $o(\vnorm{\bx-{\bar{\bx}}})$, which is a growth function. Therefore, the MAP property requires additional algorithm specific properties of the model function. In particular, we require the constants $\bar L$ and $\underline{L}$ to be independent of $\bbx$, which provides a global consistency between the model function approximations.
    \item The condition $\dom\fun \subset \scl\dom\leg$ is a minor regularity condition. For example, if $\dom\fun = [0,\infty)$ and $\dom\leg = (0,\infty)$ (e.g., for  $h$ in Burg's entropy), such a function $h$ can still be used in MAP property. However, the $L$-smad property \cite{BSTV18} would require $\bx,\bbx$ in \eqref{eq:model-ineq} to lie in $\sint\dom\leg$ (see also Section~\ref{ssec:forward-backward}). 
    \item Note that the choice of ${\underline L}$ is unrestricted in MAP property. For nonconvex $f$, $\lb$ is typically a positive real number. For convex $f$ typically the condition ${\underline L} \geq 0$ holds true. However, note that the values of $\lb,\ub$ are governed by the model function. In the context of convex additive composite problems, ${\underline L} < 0$ can hold true for relatively strongly convex functions \cite{LFN18}. 
  \end{itemize}
\end{REM}
\begin{EX}[Running Example -- Contd]\label{ex:running-example-1}
  We continue Example~\ref{ex:running-example} to illustrate the MAP property. Let $h(\bx) = \frac{1}{4}\vnorm[]{{\bf x}}^4$, we clearly have 
  \[
    g(\bx) - g(\bbx) - \scal{\nabla g(\bbx)}{\bx - \bbx} \leq 4D_h(\bx,\bbx)\,,\quad \forall\, \bx \in \R^N\,,
  \]
  which in turn results in the following upper bound for the model error
  \begin{align*}
    \abs{f(\bx) - f(\bx;\bbx)} \leq \abs{g(\bx) - g(\bbx) - \scal{\nabla g(\bbx)}{\bx - \bbx}} \leq  4D_h(\bx,\bbx)\,.
  \end{align*}
  The upper bound is obtained in terms of a Bregman distance. Clearly, the constants arising do not have any dependency on the model center. 
\end{EX}
We now present Model BPG that we analyze for the setting of Assumption~\ref{ass:model}.

\mybox{%
\begin{ALG}[{Model BPG:} Model based Bregman Proximal Gradient]\label{alg:acc-BregMin-bt} \
\begin{itemize}
  \item \emph{Initialization:} Select $\bx\iter0 = \bx\iter1 \in \dom\fun \cap \sint\dom\leg$. Choose ${\underline \tau}, {\bar \tau}$ such that $0<{\underline \tau}<{\bar \tau}< (1/{\ub})$. 
  % and  ${\sigma} > -\alpha {\bar \tau}$ holds. %[OLD]
  \item \emph{For each $\k\geq 1$:} Choose $\tau\iter\k \in  [{\underline \tau},{\bar \tau}]$ and  compute
  \begin{align}
      \label{eq:alg-BregMin-bt:update}
      \bx\iter\kp \in &\ \Argmin_{\bx\in\X}\, \left\{ \fun(\bx;\bx\iter\k) + \frac{1}{\tau\iter\k} \breg[\leg]{\bx}{\bx\iter\k} \right\}\,.
  \end{align}
\end{itemize}\vspace{-0.8em}
\end{ALG}
}

\begin{ASS}\label{ass:model}
  Let  $h$ be a Legendre function that is $\mathcal{C}^2$ over $\sint\dom\leg$. Moreover, the conditions $\dom\fun \cap \sint\dom\leg \neq \emptyset$ and $\crit\fun \cap \sint\dom\leg \neq \emptyset$ hold true. 
  \begin{enumerate}
    % \item The Legendre function $h$ is $\sigma$-strongly convex for certain $\sigma >0$. %[OLD]
    \item The exist $\bar L >0$,\,$\underline{L} \in \R$  such that for any $\bbx \in \dom\fun\, \cap\, \sint\dom\leg$, the function $f$ with $\dom f \subset \scl\dom h$, and model function $f(\cdot, \bbx)$ for $f$ around the model center $\bbx$ satisfies the \emphdef{MAP property at $\bbx$} with the constants $\bar L,\underline{L}$.
    % \item There exists $\alpha \in \R$ such that for any $\bbx \in \dom\fun\, \cap\, \sint\dom\leg$, the model function $\map{\fun({\cdot;\bbx})}{\R^N}{\eR}$ is $\alpha$-semi-convex, i.e., $f(\cdot;\bbx) - \alpha\frac{\vnorm[]{\cdot}^2}{2}$ is convex. 
    \item For any $\bbx \in \dom\fun \cap \sint\dom\leg$, the following qualification condition holds true:
    \begin{equation}\label{eq:qualification-condition}
      \partial^{\infty}_{\bx}f(\bx;\bbx) \cap (-N_{\dom h}(\bx)) = \{{\bf 0}\}\,,\quad \forall\, \bx \in \dom\fun \cap \dom\leg\,.
    \end{equation}
    % \begin{equation}\label{eq:qualification-condition}
    %   N_{\dom f}(\bx) \cap (-N_{\dom h}(\bx)) = \{{\bf 0}\}\,,\quad \forall\, \bx \in \dom\fun \cap \dom\leg\,.
    % \end{equation}
    % and the following supercoercivity condition holds true:
    % \begin{equation}
    %   \forall \lambda > 0\,,\quad \lim_{\norm{\bx} \to \infty} \frac{f(\bx) + \lambda h(\bx)}{\norm{\bx}} = \infty\,.
    % \end{equation}
    \item For all $\bx, \by \in \dom\fun$, the condition
    \[
      ({\bf 0}, \bv) \in \partial^{\infty}f(\bx;\by) \quad\text{implies}\quad \bv = {\bf 0}\,,\quad\text{and}\quad (\bv, {\bf 0}) \in \partial^{\infty}f(\bx;\by) \quad\text{implies}\quad\bv = {\bf 0}
    \]
    hold true. Moreover, $f(\bx; \by)$ is regular \cite[Definition 7.25]{Rock98} at any $(\bx, \by) \in \dom\fun \times \dom\fun$. 
    \item The function $f(\bx;\bbx)$ is a proper, lsc function and is continuous over $(\bx,\bbx) \in \dom\fun \times \dom\fun$.
  \end{enumerate} 
\end{ASS}

By $\partial_\bx f(\bx;\bbx)$ we mean the limiting subdifferential of the model function $\bx \mapsto f(\bx;\bbx)$ with $\bbx$ fixed and $\partial f(\bx;\by)$ denotes the limiting subdifferential w.r.t $(\bx,\by)$; dito for the horizon subdifferential. 
\medskip

% \paragraph{Discussion on Assumption~\ref{ass:model}.} Often, typical Legendre functions can be equipped with an additional quadratic term, which makes them strongly convex and in such a case $\sigma$ will be a positive real number. Notably, certain Legendre functions are strongly convex with respect to a different norm, for e.g., such as those distances based on Boltzmann–Shannon entropy. In this paper we focus on  Legendre functions that are strongly convex with respect to $\ell^2$-norm. However, in Section~\ref{ssec:poisson-linear-inverse-problems} we provide a brief extension of our theory with a Legendre function that is not strongly convex with respect to $\ell^2$-norm. % [OLD]

\paragraph{Discussion on Assumption~\ref{ass:model}.}  The condition $\dom f \subset \scl\dom h$ is not a restriction as one can always add an indicator function to $f$ such that the iterates never leave $\scl\dom\leg$. The qualification condition in \eqref{eq:qualification-condition} is required for the applicability of the subdifferential summation rule (see \cite[Corollary 10.9]{Rock98}).   Assumption~\ref{ass:model}(iii) and \cite[Corollary 10.11]{Rock98} ensures that for all $\bx, \by \in \dom\fun$,  the following holds true:
\begin{equation}\label{eq:main-seperable-remark}
  \partial f(\bx;\by) =  \partial_\bx f(\bx;\by) \times \partial_\by f(\bx;\by)\,,\,
  \partial^{\infty} f(\bx;\by) =  \partial^{\infty}_\bx f(\bx;\by) \times \partial^{\infty}_\by f(\bx;\by)\,.\tag{Assumption~\ref{ass:model}(iii)'}
\end{equation}
We emphasize that Assumption~\ref{ass:model}(iii) is only required for the implication \eqref{eq:main-seperable-remark}. Certain classes of functions mentioned in Section~\ref{sec:examples} satisfy \eqref{eq:main-seperable-remark} directly, instead of Assumption~\ref{ass:model}(iii). Assumption~\ref{ass:model}(iv) is typically satisfied in practice and plays a key role in Lemma~\ref{lem:critical-fixed-2}. Based on Assumption~\ref{ass:model}(iii), for any fixed $\bbx \in \dom\fun$, the model function $\fun({\bx;\bbx})$ is  regular at any $\bx \in \dom f$. Using this fact, we deduce that the model function preserves the first order information of the function, in the sense that for $\bx \in \dom\fun$ the condition $\partial_{\by} \fun({\by};\bx)|_{{\by} = \bx} = \widehat{\partial} f(\bx)$  holds true, which we prove in  Lemma~\ref{lem:first-order-info} in the appendix. Many popular algorithms such as Gradient Descent, Proximal Gradient Method, Bregman Proximal Gradient Method, Prox-Linear method are special cases of Model BPG depending on the choice of the model function and the choice of Bregman distance, thus making it a unified algorithm (also c.f. \cite{OFB19}).  Examples of model functions are provided in Section~\ref{sec:examples}, for which we verify all the assumptions. Other related model functions can also be found in \cite[Section 5]{OFB19}. 
\medskip

Let $\tau>0$, $\bbx \in \dom\fun \cap \sint\dom\leg$,  the update mapping from \eqref{eq:alg-BregMin-bt:update} of Model BPG  is defined by 
\begin{equation}\label{eq:update-mapping}
    T_{\tau}(\bbx) := \Argmin_{\bx\in\X}\, \fun(\bx;\bbx) + \frac{1}{\tau} \breg[\leg]{\bx}{\bbx}\,.
\end{equation}
Denote $\eps_k := \left(\frac{1}{\tau\iter\k}-\ub\right)>0$ and clearly $ {\underline \eps}\leq \eps\iter\k \leq {\bar \eps}$, where ${\bar \eps} := \frac{1}{{\underline{\tau}}} - \ub$ and ${\underline \eps} := \frac{1}{{\bar{\tau}}} - \ub$. 
\medskip

Well-posedness of the update step \eqref{eq:alg-BregMin-bt:update}  is given by the following result.
\begin{LEM}\label{lem:well-posedness-1}
  Let Assumption~\ref{ass:problem}, \ref{ass:model} hold true and let $\bbx \in \dom\fun \cap \sint\dom\leg$. Then, for all $0<\tau<\frac{1}{\ub}$ the set $T_{\tau}(\bbx)$ is a nonempty compact subset of $\dom\fun \cap \sint\dom\leg$.
\end{LEM}
\begin{proof}
  Firstly, note that as a consequence of MAP property due to Assumption~\ref{ass:model} and nonnegativity of Bregman distances, the following condition is satisfied
  \begin{equation}\label{eq:temp-map-property}
    f(\bx) \leq f(\bx;\bbx) + \frac{1}{\tau} D_h(\bx, \bbx)\,, \quad \forall \,\bx \in \dom f \cap \dom \leg\,.
  \end{equation}
  If the set $\dom f \cap \dom \leg$ is bounded, the objective $f(\cdot;\bbx) + \frac{1}{\tau} D_h(\cdot, \bbx)$ is coercive. Otherwise, the coercivity of $f$ implies that the objective $f(\cdot;\bbx) + \frac{1}{\tau} D_h(\cdot, \bbx)$ is coercive, due to \eqref{eq:temp-map-property}. Then, the result follows from a simple application of \cite[Lemma 3.6]{LOC2019}  and \cite[Theorem 1.9]{Rock98}.
\end{proof}
% \begin{LEM}\label{lem:well-posedness-1}
%   Let Assumption~\ref{ass:problem}, \ref{ass:model} hold true. Let $\bbx \in \dom\fun \cap \sint\dom\leg$, and also let ${\sigma} > -\alpha {\bar \tau}$ hold true. Choose $\tau$ such that $0< \tau \leq {\bar \tau} < (1/\ub)$ holds true for certain ${\bar \tau}>0$. Then, $T_{\tau}(\bbx)$ is a singleton set, and $T_{\tau}(\bbx) \subset \dom\fun \cap \sint\dom\leg$. 
% \end{LEM}
% The proof is the same as the proof of  \cite[Lemma 3.6]{LOC2019}, hence we skip it.
% We provide the proof of Lemma~\ref{lem:well-posedness-1} in Section~\ref{sec:well-posedness} in the appendix.
\medskip

We would like to highlight that Model BPG results in monotonically nonincreasing function values, which we prove below.
\begin{PROP}[Sufficient Descent Property in Function values]\label{prop:function-descent-property}Let Assumptions~\ref{ass:problem}, \ref{ass:model} hold. Also, let $\seq[\k\in\N]{\bx\iter\k}$ be a sequence generated by Model BPG, then the following holds for $k\geq 1$
  \begin{equation}\label{prop:BPG-function-1:1}
    f(\bx\iter\kp)\leq f(\bx\iter\k) - \eps\iter\k D_h(\bx\iter\kp,\bx\iter\k)\,.
  \end{equation}
  \end{PROP} 
  We provide the proof of Proposition~\ref{prop:function-descent-property} in Section~\ref{sec:function-descent} in the appendix.
\medskip

\begin{REM}\label{rem:boundedness}
  Under Assumptions~\ref{ass:problem}, \ref{ass:model}, the coercivity of $f$ along with Proposition~\ref{prop:function-descent-property} implies that the iterates of Model BPG lie in the compact convex set $\{\bx : f(\bx) \leq f(\bx_0)\}$, thus bounded. 
\end{REM}

\section{Global convergence analysis of Model BPG algorithm}\label{sec:convergence-analysis}
The convergence analysis of most algorithms in nonconvex optimization is based on a descent property. Usually, the objective value is shown to decrease, for example, as in Proposition~\ref{prop:function-descent-property} and in the analysis of additive composite problems \cite[Lemma 4.1]{BSTV18}. However, function values proved to be restrictive, primarily because the same techniques as additive composite problems do not work anymore for general composite problems, and alternatives like \cite{P2016} are sought after.   

\subsection{New Lyapunov function}
Here, we discuss one of our main contribution. We propose a Lyapunov function as our measure of progress. The Lyapunov function $F^\leg_{\ub}$ is given by 
\begin{equation}\label{eq:lyapunov-function}
  \map{F^\leg_{\ub}}{\X\times\X}{\eR}\,,\quad 
  (\bx,\bbx) \mapsto 
    \fun(\bx;\bbx) +  \ub  \breg[\leg]{\bx}{\bbx}\,,
\end{equation}
and  $\dom F^\leg_{\ub} = \dom\fun^2 \times \dom D_h\,.$  The set of critical points of the above given Lyapunov function is given by 
\begin{equation}
  \crit F^\leg_{\ub} := \left\{ \left(\bx, \bbx\right) \in \R^N \times \R^N: \;  ({\bf 0}, {\bf 0}) \in \partial F^\leg_{\ub}(\bx, \bbx)  \right\}\,.
\end{equation}
Usage of Lyapunov functions is  a popular strategy in the analysis of inertial methods \cite{OCBP14,MOPS2020}. Even though our algorithm is non-inertial in nature, we show that the above defined Lyapunov function is suitable for the global convergence analysis. Certain previous works such as \cite{MLF2019} considered a Lyapunov function based analysis for (non-inertial) Forward–Douglas–Rachford splitting method. Also, Lyapunov function based analysis is popular in the context of dynamical systems \cite{HC2011}. 
\medskip

The motivation for using the Lyapunov function $F^\leg_{\ub}$ instead of the function $f$ is the following. In each iteration of Model BPG, we optimize the model function with a proximity measure, and the analysis with our proposed Lyapunov function reflects this explicitly, unlike the function value. The proposed Lyapunov function is related to the Bregman-Moreau envelope \cite{LOC2019} of the model function $f(\cdot; \bbx)$ where $\bbx \in \dom f \cap \sint\dom\leg$. Under certain special case of the model function (Section~\ref{ssec:forward-backward}), such a Bregman-Moreau envelope is related to the  Bregman forward-backward envelope \cite{ATP2019}. In the context where the Bregman distance is set to the Euclidean distance, the related works which consider value function based analysis is provided \cite{BP2016,P2016,STP2017}.
\medskip

We now look at some properties of $F^\leg_{\ub}$.
\begin{PROP}\label{prop:lya-fun-1}
  The Lyapunov function defined in \eqref{eq:lyapunov-function} satisfies the following properties:
  \begin{itemize}
  \item[$\rm{(i)}$] For all $\bx\in \dom \fun \cap \dom\leg$ and $\by  \in \dom\fun \cap \sint\dom\leg$, we have
  $
    f(\bx) \leq F^\leg_{\ub}(\bx,\by)\,.
  $
  \item[$\rm{(ii)}$]  For all $\bx \in \dom\fun \cap \sint\dom\leg$, we have
  $
    F^\leg_{\ub}(\bx,\bx) = \fun(\bx)\,.
  $
  \item[$\rm{(iii)}$] Moreover, we have
  \begin{equation}
  \inf_{(\bx,\by)\,\in\, \R^N\times\R^N} F^\leg_{\ub}(\bx,\by) \geq v(\PPP) > -\infty\,.
  \end{equation}
  \end{itemize}
  \end{PROP}
\begin{proof}
\begin{itemize}
  \item[$\rm{(i)}$] This follows from  \textit{MAP} property and the definition of $F^\leg_{\ub}$ . 
  \item[$\rm{(ii)}$] Substituting $\by=\bx$ in \eqref{eq:lyapunov-function} gives the result.
\item[$\rm{(iii)}$] By \textit{MAP} property, for all $(\bx,\by) \in \dom F^\leg_{\ub}$ we have the following:
\[
  v(\PPP)\leq \fun(\bx) \leq \fun(\bx;\by) + \ub D_h(\bx,\by)\,.
\]
Furthermore, we obtain the following:
\[
  \inf_{\bx \in \dom\fun\,\cap\,\dom\leg}\fun(\bx) \leq \inf_{(\bx,\by) \in \dom F^\leg_{\ub}}\left(\fun(\bx;\by) + \ub D_h(\bx,\by)\right)\,.
\]
The statement follows using $\inf_{\bx \in \R^N}\fun(\bx) = v(\PPP) > -\infty$ due to Assumption~\ref{ass:problem}\,.
\end{itemize} 
\end{proof}

Equipped with the Lyapunov function $F^\leg_{\ub}$, we focus now on the global convergence result of Model BPG. Our global convergence analysis is broadly divided into the following five parts.
\begin{itemize}
  \item \textbf{Sufficient descent property.} In Section~\ref{ssec:suff-desc}, we show that the sequence generated by Model BPG results in monotonically nonincreasing Lyapunov function values. 
  \item \textbf{Relative error condition.} In Section~\ref{ssec:relative-error}, based on certain additional assumptions, we show that the infimal norm of the subdifferential of the Lyapunov function can be upper bounded by an entity that depends on the difference of successive iterates, and that entity tends towards zero asymptotically, implying stationarity in the limit. 
  \item \textbf{Subsequential convergence.} In Section~\ref{ssec:subsequential-convergence}, we explore the behavior of limit points obtained from the sequence generated by Model BPG.  We prove $F_{\ub}^\leg$-attentive convergence along converging subsequences. Moreover, we prove that the set of $F_{\ub}^\leg$-attentive limit points is compact, connected and $F_{\ub}^\leg$ is constant on this set. When all limit points of the sequence generated by Model BPG lie in $\sint\dom\leg$, we show that all the limit points are critical points of the Lyapunov function.
  \item \textbf{Global convergence to stationarity point of the  Lyapunov function.} Under the condition that the Lyapunov function satisfies \KL property, we show in Section~\ref{ssec:global-convergence-lyapunov-func} that the full sequence generated by Model BPG converges to a point $\bx$ such that $(\bx, \bx)$ is the critical point of the Lyapunov function. However, the relation of $\bx$ to the function $f$ is not imminent here.
  \item \textbf{Global convergence to stationarity point of the function.}  In Section~\ref{ssec:global-convergence-function}, we prove that the update mapping is continuous and also show that fixed points of the update mapping are critical points of $f$. We exploit these properties to deduce that the full sequence of iterates generated by Model BPG converges to a critical point of $f$. 
\end{itemize}

\subsection{Sufficient descent property}\label{ssec:suff-desc}
We have already proved the sufficient descent property in terms of function values in Proposition~\ref{prop:function-descent-property}.  Here, we prove the sufficient descent property of the Lyapunov function.
\begin{PROP}[Sufficient descent property]\label{prop:descent-property}Let Assumptions~\ref{ass:problem}, \ref{ass:model} hold. Also, let $\seq[\k\in\N]{\bx\iter\k}$ be a sequence generated by Model BPG, then the following holds for $k\geq 1$
\begin{equation}\label{prop:BPG-Lyapunov-1:1}
F_{\ub}^{\leg}(\bx\iter\kp,\bx\iter\k) \leq F_{\ub}^{\leg}(\bx\iter\k,\bx\iter\km) - \eps\iter\k D_h(\bx\iter\kp,\bx\iter\k)\,.
\end{equation}
\end{PROP} 
\begin{proof}
By global optimality of $\bx\iter\kp$ as in \eqref{eq:alg-BregMin-bt:update}, we have
\begin{align*}
\fun(\bx\iter\kp;\bx\iter\k)+  \frac{1}{\tau\iter\k}D_h(\bx\iter\kp,\bx\iter\k) \leq \fun(\bx\iter\k;\bx\iter\k) = \fun(\bx\iter\k)\,.
\end{align*}
We have the following inequality from the MAP property
\begin{align*}
\fun(\bx\iter\k;\bx\iter\k) = \fun(\bx\iter\k) \leq  \fun(\bx\iter\k;\bx\iter\km)  + \ub D_h(\bx\iter\k,\bx\iter\km) \,.
\end{align*}
Thus, the result follows from the definition of $F_{\ub}^{\leg}$ in \eqref{eq:lyapunov-function}.
% \begin{align*}
%  \fun(\bx\iter\kp;\bx\iter\k) +  \ub D_h(\bx\iter\kp,\bx\iter\k) \leq    \fun(\bx\iter\k;\bx\iter\km)  + \ub D_h(\bx\iter\k,\bx\iter\km) - \eps\iter\k D_h(\bx\iter\kp,\bx\iter\k)\,.
% \end{align*}
\end{proof}
\begin{PROP}\label{p:SuffDesc}
  Let Assumptions~\ref{ass:problem}, \ref{ass:model} hold and let $\seq[\k\in\N]{\bx\iter\k}$ be a sequence generated by Model BPG. The following assertions hold: 
    \begin{itemize}
      \item[$\rm{(i)}$] The sequence $\left\{ F_{\ub}^{\leg}\left(\bx\iter\kp , \bx\iter\k\right) \right\}_{k \in \N}$ is nonincreasing and converges to a finite value.
      \item[$\rm{(ii)}$] $\sum_{k = 1}^{\infty} D_h(\bx\iter\kp,\bx\iter\k)  < \infty$, and hence the sequence $\left\{ D_h(\bx\iter\kp,\bx\iter\k)  \right\}_{k \in \N}$ converges to zero.
      \item[$\rm{(iii)}$] For any $n\in \N$, the condition 
      \[
        \min_{1 \leq k \leq n} D_h(\bx\iter\kp,\bx\iter\k)  \leq \frac{F_{\ub}^{\leg}\left(\bx_{1} , \bx_{0}\right)-v(\PPP)}{{\underline \eps} n}
      \]
      holds true.
    \end{itemize}       
\end{PROP}
\begin{proof}
  \begin{itemize}
    \item[$\rm{(i)}$]  Nonincreasing property follows trivially from Proposition~\ref{prop:descent-property} and as $\eps\iter\k >0$. We know from Proposition~\ref{prop:lya-fun-1}(iii) that the Lyapunov function is lower bounded, which implies convergence of $\left\{ F_{\ub}^{\leg}\left(\bx\iter\kp , \bx\iter\k\right) \right\}_{k \in \N}$ to a finite value. 
    \item[$\rm{(ii)}$] Let $n$ be a positive integer. Summing \eqref{prop:BPG-Lyapunov-1:1} from $k = 1$ to $n$ and using ${\underline \eps} \leq {\eps\iter\k}$ we get
      \begin{equation} \label{P:SuffDesc0:2}
        \sum_{k = 1}^{n} D_h(\bx\iter\kp,\bx\iter\k)  \leq \frac{1}{{\underline \eps}}\left(F_{\ub}^{\leg}\left(\bx_{1} , \bx_{0}\right) - F_{\ub}^{\leg}\left(\bx_{n + 1} , \bx_{n}\right)\right) \leq \frac {1}{{\underline \eps}}\left(F_{\ub}^{\leg}\left(\bx_{1} , \bx_{0}\right) -v(\PPP)\right), 
      \end{equation}
       since $F_{\ub}^{\leg}\left(\bx_{n + 1} , \bx_{n}\right) \geq v(\PPP)$. Taking the limit as $n \rightarrow \infty$, we obtain the first  assertion, from which  we immediately deduce that $\left\{D_h(\bx\iter\kp,\bx\iter\k) \right\}_{k \in \N}$ converges to zero.
    \item[$\rm{(iii)}$] From \eqref{P:SuffDesc0:2} we also obtain,
      \begin{equation*}
        n\min_{1\leq k \leq n} \left(D_h(\bx\iter\kp,\bx\iter\k) \right) \leq \sum_{k = 1}^{n} \left(D_h(\bx\iter\kp,\bx\iter\k) \right) \leq \frac{1}{{\underline \eps}}\left(F_{\ub}^{\leg}\left(\bx_{1} , \bx_{0}\right) -v(\PPP)\right),
      \end{equation*}
      which after division by $n$ yields the result.
  \end{itemize}
  \vspace{-0.2in}
\end{proof}

\subsection{Relative error condition}\label{ssec:relative-error}
For the purposes of analysis, we require the following assumption.
\begin{ASS}\label{ass:relative-error} We have the following conditions:
\begin{enumerate}
  \item  Consider any bounded set $B\subset \dom\fun$. There exists $c>0$ such that for any $\bx, \by \in B$ we have 
    \[
      \inf_{\bv \in \partial_{\by}\fun(\bx;\by)}\vnorm[]{\bv} \leq c\vnorm[]{\bx - \by} \,.
    \]
  \item The function $\leg$ has bounded second derivative on any bounded subset $B\subset \sint\dom\leg$. 
  \item  For bounded $\seq[\k\in\N]{{\bf u}\iter\k}$, $\seq[\k\in\N]{\bv\iter\k}$ in $\sint\dom\leg$, the following holds as $\k\to\infty$:
  \[
    \breg[\leg]{{\bf u}\iter\k}{\bv\iter\k} \to 0  \quad \iff \quad
    \vnorm[]{{\bf u}\iter\k - \bv\iter\k} \to 0 \,.
  \]
\end{enumerate}
\end{ASS}

Through Example~\ref{ex:bounded-second-order}, we illustrate Assumption~\ref{ass:relative-error}(i), which governs the variation of the model function w.r.t. model center.   Assumption~\ref{ass:relative-error}(ii) is a standard condition required for the analysis of Bregman proximal methods \cite{BSTV18, OFB19, MOPS2020}.  Assumption~\ref{ass:relative-error}(iii) essentially states that the asymptotic behavior of vanishing Bregman distance is equivalent to that of vanishing Euclidean distance (cf. \cite[Remark 18]{OFB19}). Such a condition is satisfied for many Bregman distances, such as those distances based on Boltzmann–Shannon entropy \cite[Example 40]{OFB19} and Burg entropy \cite[Example 41]{OFB19}.

\begin{EX}\label{ex:bounded-second-order}
  We continue Example~\ref{ex:running-example} to illustrate Assumption~\ref{ass:relative-error}(i). A quick calculation reveals that $\nabla^2g(\bx)$ is bounded over bounded sets. Consider any bounded set $B\subset \R^N$. Define $c := \sup_{\bbx \in B}\vnorm[]{\nabla^2g(\bbx)}$ and choose any $\bbx \in B$, then consider the model function given by :
  \[
    f(\bx;\bbx) := \abs{g(\bbx) + \scal{\nabla g(\bbx)}{\bx - \bbx}}\,.
  \]
  The subdifferential of the model function  is given by 
  \[
  \partial_\bbx f(\bx; \bbx) =   {\bf u}\nabla^2g(\bbx)(\bx -\bbx)\,,
  \]
  where ${\bf u} \in \partial_{g(\bbx) + \scal{\nabla g(\bbx)}{\bx - \bbx}} \abs{g(\bbx) + \scal{\nabla g(\bbx)}{\bx - \bbx}}$. Considering the fact that $\vnorm[]{{\bf u}} \leq 1$ and by the definition of $c$ we have the following: 
  \[
  \inf_{\bv \in \partial_{\bbx} f(\bx; \bbx)}   \vnorm[]{\bv} \leq c \vnorm[]{\bx - \bbx}\,,
  \] 
  which verifies Assumption~\ref{ass:relative-error}(i).
\end{EX}

\medskip
 Now, we look at the relative error condition, which bounds the infimal norm of the subdifferential of the Lyapunov function, i.e., $\inf_{v \in \partial  F^h_{\ub} (\bx\iter\kp, \bx\iter\k)}\vnorm{v}$, with the term $\vnorm[]{\bx\iter\kp - \bx\iter\k}$ upto a scaling factor. Such a bound is useful to achieve stationarity asymptotically, and plays a crucial role in proving global convergence. Note that with the descent property (Proposition~\ref{prop:descent-property}) and Assumption~\ref{ass:relative-error}(iii), we have $\vnorm[]{\bx\iter\kp - \bx\iter\k} \to 0$.

\begin{LEM}[Relative error] \label{lem:sub-diff-vanish} 
  Let Assumptions~\ref{ass:problem}, \ref{ass:model}, \ref{ass:relative-error}  hold. Let the sequence $\seq[\k\in\N]{\bx\iter\k}$ be generated by Model BPG, then there exists a constant $C>0$ such that for certain $k \geq 0$, we have 
  \begin{align} \label{eq:lem:sub-diff-vanish}
      \vnorm[-]{\partial  F^h_{\ub} (\bx\iter\kp, \bx\iter\k)} 
      \leq C \vnorm[]{\bx\iter\kp - \bx\iter\k}\,,
  \end{align}
  where $\vnorm[-]{\partial  F^h_{\ub} (\bx\iter\kp, \bx\iter\k)} := \inf_{v \in \partial  F^h_{\ub} (\bx\iter\kp, \bx\iter\k)}\vnorm{v}$.
\end{LEM}
\begin{proof}
As per \cite[Theorem 2.19]{M2018}, the subdifferential $\partial  F^h_{\ub} (\bx\iter\kp, \bx\iter\k)$ is given by 
\begin{equation}\label{eq:subdiff-part-1}
  \partial  F^h_{\ub} (\bx\iter\kp, \bx\iter\k) = \partial f(\bx\iter\kp;\bx\iter\k) + {\bar L}\nabla D_h(\bx\iter\kp,\bx\iter\k)\,,
\end{equation}
because the Bregman distance is continuously differentiable around $\bx\iter\k \in \dom\fun \cap \sint\dom\leg$. Using \cite[Corollary 10.11]{Rock98}, Assumption~\ref{ass:model}(iv), and using the fact that $h$ is $\mathcal{C}^2$ over $\sint\dom\leg$ (cf. Assumption~\ref{ass:model}) we obtain 
\begin{align}
    \partial F_{\ub}^\leg (\bx\iter\kp, \bx\iter\k) 
    = \Big(&\partial_{\bx\iter\kp} f(\bx\iter\kp;\bx\iter\k) + \ub\big(\nabla \leg(\bx\iter\kp) - \nabla \leg(\bx\iter\k)\big),\nonumber\\
    &\,\label{eq:lem:sub-diff-vanish:proof:1}
    \partial_{\bx\iter\k} f(\bx\iter\kp;\bx\iter\k) - \ub\nabla^2\leg(\bx\iter\k)(\bx\iter\kp - \bx\iter\k) \Big) \,.
\end{align}
Consider the following:
\begin{align}
  \inf_{\zeta \in \partial F(\bx\iter\kp,\bx\iter\k)}\vnorm[]{v} &= \inf_{\xi \in \partial f(\bx\iter\kp;\bx\iter\k)}\vnorm[]{\xi + {\bar L}\nabla D_h(x\iter\kp;x\iter\k)}\,,\nonumber\\
  &= \left(\inf_{(\xi_x ,\xi_y ) \in \partial f(\bx\iter\kp;\bx\iter\k)}\vnorm[]{(\xi_x ,\xi_y)  + {\bar L}\nabla D_h(x\iter\kp,x\iter\k)}\right)\,,\nonumber\\
  &\leq \left(\inf_{\xi_x \in \partial_{\bx\iter\kp} f(\bx\iter\kp;\bx\iter\k)}\vnorm[]{(\xi_x +\ub\big(\nabla \leg(\bx\iter\kp) - \nabla \leg(\bx\iter\k)))} \nonumber  \right)\\
  &+ \left(\inf_{\xi_y \in \partial_{\bx\iter\k} f(\bx\iter\kp;\bx\iter\k)}\vnorm[]{(\xi_y +\ub\nabla^2\leg(\bx\iter\k)(\bx\iter\kp - \bx\iter\k))}   \right)\,,\label{eq:subdiff-main-1}
\end{align}
where in the  first equality we use \eqref{eq:subdiff-part-1}, in the second equality we use the result in \eqref{eq:lem:sub-diff-vanish:proof:1}  with $\xi := (\xi_x , \xi_y)$ such that $\xi_x \in \partial_{\bx\iter\kp} f(\bx\iter\kp,\bx\iter\k)$ and $\xi_y \in \partial_{\bx\iter\k} f(\bx\iter\kp,\bx\iter\k)$, and in the last step we used 
\begin{align}
  \nabla D_h(\bx\iter\kp,\bx\iter\k) = (\nabla h(\bx\iter\kp) - \nabla h(\bx\iter\k), \nabla^2 h(\bx\iter\k)(\bx\iter\kp - \bx\iter\k))\,.
\end{align}
The optimality of $\bx\iter\kp$ in \eqref{eq:alg-BregMin-bt:update} implies the existence of $\xi_{\bx\iter\kp}^{\k+1}\in \partial_{\bx\iter\kp} \fun({\bx\iter\kp;\bx\iter\k})$ such that the following condition holds:
  \begin{equation} \label{eq:update-opt-cond}
    \xi_{\bx\iter\kp}^{\k+1} + \frac{1}{\tau_{\iter\k}}  (\nabla\leg(\bx\iter\kp) - \nabla\leg(\bx\iter\k)) = {\bf 0}\,.
  \end{equation}
Therefore, the first block coordinate in \eqref{eq:lem:sub-diff-vanish:proof:1} satisfies
\begin{align}
\xi^{k+1}_{\bx\iter\kp} + \ub\big(\nabla \leg(\bx\iter\kp) - \nabla \leg(\bx\iter\k)\big) = \eps\iter\k\big(\nabla \leg(\bx\iter\kp) - \nabla \leg(\bx\iter\k)\big)\,.\label{eq:part-1-subdiff}
\end{align}
Now consider the first term of the right hand side  in \eqref{eq:subdiff-main-1}. We have 
\begin{align*}
\inf_{\xi_x \in \partial_{\bx\iter\kp} f(\bx\iter\kp;\bx\iter\k)}\vnorm[]{(\xi_x +\ub\big(\nabla \leg(\bx\iter\kp) - \nabla \leg(\bx\iter\k)))}  &\leq \vnorm[]{\xi^{k+1}_{\bx\iter\kp} + \ub\big(\nabla \leg(\bx\iter\kp) - \nabla \leg(\bx\iter\k)\big)}\,,\\
&\leq \eps\iter\k\vnorm[]{\big(\nabla \leg(\bx\iter\kp) - \nabla \leg(\bx\iter\k)\big)}\,,\\
&\leq \eps\iter\k \tilde{L}_h \vnorm[]{\bx\iter\kp -\bx\iter\k}\,,
\end{align*}
where in the second step we used \eqref{eq:part-1-subdiff} and in the last step we applied mean value theorem along with the fact that the entity $\vnorm[]{\nabla^2h(\bx\iter\kp + s(\bx\iter\kp - \bx\iter\k))}$ is bounded by a constant $\tilde{L}_h > 0$ for certain $s \in [0,1]$, due to Assumption~\ref{ass:relative-error}(ii). Considering the second term of the right hand side in \eqref{eq:subdiff-main-1}, we have 
\begin{align*}
  \inf_{\xi_y \in \partial_{\bx\iter\k} f(\bx\iter\kp;\bx\iter\k)}\vnorm[]{(\xi_y +\ub\nabla^2\leg(\bx\iter\k)(\bx\iter\kp - \bx\iter\k))}  &\leq \inf_{\xi_y \in \partial_{\bx\iter\k} f(\bx\iter\kp;\bx\iter\k)}\vnorm[]{\xi_y} + \vnorm[]{\ub\nabla^2\leg(\bx\iter\k)(\bx\iter\kp - \bx\iter\k)}\,,\\
  &\leq c\vnorm[]{\bx\iter\kp - \bx\iter\k} + \ub L_h\vnorm[]{(\bx\iter\kp - \bx\iter\k)}\,, 
\end{align*}
where in the last step we used Assumption~\ref{ass:relative-error}(i) and the fact that $\vnorm[]{\nabla^2h(\bx\iter\k)}$ is bounded by $L_h$. The result follows from combining the results obtained for \eqref{eq:part-1-subdiff}.
\end{proof}

\subsection{Subsequential convergence}\label{ssec:subsequential-convergence}
We now consider results on generic limit points and show that stationarity can indeed be attained for iterates produced by Model BPG. The set of limit points of some sequence $\seq[\k\in\N]{\bx\iter\k}$  is denoted as follows 
\[
  \omega(\bx\iter0) := \set{\bx\in \X \setsep \exists K\subset\N\colon \bx\iter\k \rto{k\in K} \bx} \,,
\]
and its subset of $f$-attentive limit points 
\[
  \omega_f(\bx\iter0) := \set{\bx\in \X \setsep \exists K\subset\N\colon (\bx\iter\k,f(\bx\iter\k)) \rto{\k\in K} (\bx,f(\bx))} \,.
\]

We explore below certain properties that are generic to any bounded sequence, and are later helpful to quantify properties of the sequence generated by Model BPG.

\begin{PROP}\label{prop:subsequence}
  For a bounded sequence $\seq[\k\in\N]{\bx\iter\k}$ such that $\vnorm[]{\bx\iter\kp - \bx\iter\k} \to 0$ as $k \to \infty$, the following holds:
  \begin{enumerate}
    \item[$\ii1$]  $\omega(\bx\iter0)$ is connected and compact, 
    \item[$\ii2$] $\lim_{\k\to\infty} \dist(\bx\iter\k,\omega(\bx\iter0)) = 0$.
  \end{enumerate}
\end{PROP}
The proof relies on the same technique as the proof of \cite[Lemma 3.5]{BST14} (also see \cite[Remark 3.3]{BST14}).
\medskip
% \begin{proof}
%     $\ii1$ The connectedness of $\omega(\bx\iter0)$  is a simple application of \cite[Lemma 3.5]{BST14} along with assumption $\vnorm[]{\bx\iter\kp - \bx\iter\k} \to 0$ as $k \to \infty$ and compactness is a direct consequence of its definition as an outer set-limit and boundedness. $\ii2$ holds by definition.  
% \end{proof}

We now show that the sequence generated by Model BPG $\seq[\k\in\N]{\bx\iter\k}$ indeed attains $\vnorm[]{\bx\iter\kp - \bx\iter\k} \to 0$ as $k \to \infty$, which in turn enables the application of Proposition~\ref{prop:subsequence} to deduce the properties of the sequence generated by Model BPG, which later proves to be crucial for the proof of global convergence.

\begin{PROP}\label{prop:limit-point-sets-1}
  Let Assumption~\ref{ass:problem}, \ref{ass:model}, \ref{ass:relative-error} hold. Let $\seq[\k\in\N]{\bx\iter\k}$ be a sequence generated by Model BPG. Then, we have 
    \begin{equation} \label{eq:breg-diff-vanish}
      {\underline \eps} \breg[\leg]{\bx\iter\kp}{\bx\iter\k}  \to 0\,,\quad \text{as}\ \k\to\infty\,.
    \end{equation}
    The condition ${\underline \eps}>0$ implies that $\bx\iter\kp-\bx\iter\k \to 0$ as $\k\to\infty$. 
  \end{PROP}
\begin{proof}
Note that the sequence $\seq[\k\in\N]{\bx\iter\k}$ is a bounded sequence (see Remark~\ref{rem:boundedness}).  By the Descent Property (Proposition ~\ref{prop:descent-property}) and using $\eps\iter\k \geq {\underline \eps}$ we have after rearranging
\[
{\underline \eps} \breg[\leg]{\bx\iter\kp}{\bx\iter\k} \leq F_{\ub}^\leg(\bx\iter\k,\bx\iter\km)  -  F_{\ub}^\leg(\bx\iter\kp,\bx\iter\k)\,.
\]
Summing on both sides and due to the convergence of Lyapunov function, using Proposition~\ref{prop:descent-property}, we obtain
\[
\sum_{k=1}^\infty \Big({\underline \eps} \breg[\leg]{\bx\iter\kp}{\bx\iter\k}\Big) \leq  F_{\ub}^\leg(\bx_0,\bx_{-1}) - \lim_{\k\to\infty} F_{\ub}^\leg(\bx\iter\kp,\bx\iter\k) < \infty\,,
\]
which implies \eqref{eq:breg-diff-vanish}.  For ${\underline \eps}>0$, Assumption~\ref{ass:relative-error}(iii)  together with \eqref{eq:breg-diff-vanish} imply $\bx\iter\kp-\bx\iter\k \to 0$ as $\k\to\infty$. 
\end{proof}
Analyzing the full set of limit points of the sequence generated by Model BPG is difficult, as illustrated in \cite{OFB19}. Obtaining the global convergence is still an open problem. Moreover, the work in \cite{OFB19} relies on convex model functions.
\medskip

In order to simplify slightly the setting, we restrict the set of limit points to the set $\sint\dom\leg$. Such a choice may appear to be restrictive, however, Model BPG when applied to many practical problems results in sequences that have  this property as illustrated in Section~\ref{sec:experiments}. 
\medskip

To this regard, denote the following
\[
  \omega^{\sint\dom\leg}(\bx\iter0):=\omega(\bx\iter0)\cap\sint\dom\leg \quad\text{and}\quad\omega_f^{\sint\dom\leg}(\bx\iter0):=\omega_f(\bx\iter0)\cap\sint\dom\leg\,.
\]
The subset of $F_{\ub}^\leg$-attentive (similar to $f$-attentive) limit points is
\[
  \omega_{F_{\ub}^\leg}(\bx\iter0) := \set{(\by,\bx)\in \X\times \X \setsep \exists K\subset\N\colon (\bx\iter\k,F_{\ub}^\leg(\bx\iter\k,\bx\iter\km)) \rto{\k\in K} (\bx,F_{\ub}^\leg(\by,\bx))}\,.
\]
Also, we define $\omega_{F_{\ub}^\leg}^{(\sint\dom\leg)^2} := \omega_{F_{\ub}^\leg} \cap (\sint\dom\leg \times \sint\dom\leg)$.  
% \MMC[inline]{Check again.}
\begin{PROP}\label{prop:limit-points} Let Assumptions~\ref{ass:problem}, \ref{ass:model}, \ref{ass:relative-error}  hold.  Let $\seq[\k\in\N]{\bx\iter\k}$ be a sequence generated by Model BPG. Then, the following holds:
  \begin{enumerate}
    \item[$\ii1$] $\omega^{\sint\dom\leg}(\bx\iter0)=\omega_\fun^{\sint\dom\leg}(\bx\iter0)$, 
    \item[$\ii2$] $\bx\in \omega_f^{\sint\dom\leg}(\bx\iter0)$ if and only if $(\bx,\bx) \in \omega_{F_{\ub}^\leg}^{(\sint\dom\leg)^2}(\bx\iter0)$. 
    \item[$\ii3$] $F_{\ub}^\leg$ is constant and finite on $\omega_{F_{\ub}^\leg}^{(\sint\dom\leg)^2}(\bx\iter0)$ and $\fun$ is constant and finite on $\omega_\fun^{\sint\dom\leg}(\bx\iter0)$ with  same value.
  \end{enumerate}
\end{PROP}
\begin{proof}
  $\ii1$ We show the inclusion $\omega^{\sint\dom\leg}(\bx\iter0)\subset\omega_\fun^{\sint\dom\leg}(\bx\iter0)$ and $\omega_\fun^{\sint\dom\leg}(\bx\iter0)\subset\omega^{\sint\dom\leg}(\bx\iter0)$ is clear by definition. Let $\opt{\bx}\in \omega^{\sint\dom\leg}(\bx\iter0)$, then we obtain the following
  \begin{multline*}
    \fun(\opt{\bx}) +  \left({\lb} +\frac{1}{\tau\iter\k}\right)\breg[\leg]{\opt{\bx}}{\bx\iter\k}  \overset{\eqref{eq:model-ineq}}\geq \fun(\opt{\bx};\bx\iter\k) + \frac{1}{\tau\iter\k}\breg[\leg]{\opt{\bx}}{\bx\iter\k} 
    \overset{\eqref{eq:alg-BregMin-bt:update}}\geq \fun({\bx\iter\kp;\bx\iter\k}) + \frac{1}{\tau\iter\k}\breg[\leg]{\bx\iter\kp}{\bx\iter\k}  \\
    \overset{\eqref{eq:model-ineq}}\geq \fun(\bx\iter\kp) - \Big(\ub-\frac 1{\tau\iter\k}\Big) \breg[\leg]{\bx\iter\kp}{\bx\iter\k} 
    \overset{\eps_k >0}\geq \fun(\bx\iter\kp) \,.
  \end{multline*}
  Obviously, by Assumption~\ref{ass:relative-error}(iii) combined with the fact that $\bx\iter\k \rto{K} \opt{\bx}$, we have $\breg[\leg]{\opt{\bx}}{\bx\iter\k}\to0$ as $\k\rto{K}\infty$,  which, together with the lower semicontinuity of $\fun$, implies 
  \[
    \fun(\opt{\bx}) \geq \liminf_{\k\rto{K}\infty} \fun(\bx\iter\kp) \geq \fun(\opt{\bx}) \,,
  \]
  thus $\opt{\bx}\in\omega_\fun^{\sint\dom\leg}(\bx\iter0)$.
\medskip

$\ii2$ If  $\bx\in\omega_\fun^{\sint\dom\leg}(\bx\iter0)$,  then we have  $\bx\iter\k\rto{\k\in K}\bx$ for $K\subset\N$, and $\fun(\bx\iter\k) \rto{\k\in K} \fun(\bx)$. As a consequence of Proposition~\ref{p:SuffDesc} and Assumption~\ref{ass:relative-error}(iii), $\breg[\leg]{\bx\iter\kp}{\bx\iter\k}\to 0$ as $\k\to\infty$, which implies that $\bx\iter\kp\rto{\k\in K}\bx$. The first part  of the proof implies  $\fun(\bx\iter\kp) \rto{\k\in K} \fun(\bx)$.  We also have $F_{\ub}^\leg(\bx\iter\kp,\bx\iter\k)\rto{\k\in K} \fun(\bx)$ which we prove below, which implies that $(\bx,\bx)\in \omega_{F_{\ub}^\leg}^{\sint\dom\leg}(\bx\iter0)$. Note that by definition of $F_{\ub}^\leg$  we have the following
\begin{align*}
  F_{\ub}^{\leg}(\bx\iter\kp,\bx\iter\k) &= \fun(\bx\iter\kp;\bx\iter\k) +  \ub D_h(\bx\iter\kp,\bx\iter\k)\,,\\
  &= \fun(\bx\iter\kp) + (\fun(\bx\iter\kp;\bx\iter\k) - \fun(\bx\iter\kp)) +  \ub D_h(\bx\iter\kp,\bx\iter\k)\,, 
\end{align*}
and with the MAP property we have
\begin{equation}\label{eq:lyapunov-fun-1}
\fun(\bx\iter\kp)  \leq F_{\ub}^{\leg}(\bx\iter\kp,\bx\iter\k)  \leq \fun(\bx\iter\kp) + (\ub + \lb)D_h(\bx\iter\kp,\bx\iter\k) \,.
\end{equation}
Thus, we have that $F_{\ub}^\leg(\bx\iter\kp,\bx\iter\k)\rto{\k\in K} \fun(\bx)$  as $D_h(\bx\iter\kp,\bx\iter\k) \rto{\k\in K} 0$.  Conversely, suppose $(\bx,\bx)\in \omega_{F_{\ub}^\leg}^{\sint\dom\leg}(\bx\iter0)$ and $\bx\iter\k\rto{\k\in K}\bx$ for $K\subset\N$. This, together with $\breg[\leg]{\bx\iter\kp}{\bx\iter\k}\to 0$ as $k \rto{\k\in K} \infty$, induces $F_{\ub}^\leg(\bx\iter\kp,\bx\iter\k)\rto{\k\in K} \fun(\bx)$, which further implies $\fun(\bx\iter\kp)\rto{\k\in K} \fun(\bx)$ due to the following. Note that we have
\begin{align*}
\fun(\bx\iter\kp) &= F_{\ub}^{\leg}(\bx\iter\kp,\bx\iter\k) +   (\fun(\bx\iter\kp) - \fun(\bx\iter\kp;\bx\iter\k)) +  \ub D_h(\bx\iter\kp,\bx\iter\k)\\
&\geq F_{\ub}^{\leg}(\bx\iter\kp,\bx\iter\k) + (\ub -\lb)D_h(\bx\iter\kp,\bx\iter\k)  \,.
\end{align*}
Finally we have
\begin{align*}
F_{\ub}^{\leg}(\bx\iter\kp,\bx\iter\k) + (\ub -\lb)D_h(\bx\iter\kp,\bx\iter\k)   \leq \fun(\bx\iter\kp) \leq F_{\ub}^{\leg}(\bx\iter\kp,\bx\iter\k) \,.
\end{align*}
Thus, with $D_h(\bx\iter\kp,\bx\iter\k) \to 0$ as $k \rto{\k\in K} \infty$ and $F_{\ub}^{\leg}(\bx\iter\kp,\bx\iter\k) \rto{\k\in K} \fun(\bx)$, we deduce that $\fun(\bx\iter\kp) \rto{\k\in K} \fun(\bx)$. And therefore $\bx\in \omega_\fun^{\sint\dom\leg}(\bx\iter0)$.
\medskip

$\ii3$ By Proposition~\ref{prop:descent-property}, the sequence $\seq[\k\in\N]{F_{\ub}^\leg(\bx\iter\kp,\bx\iter\k)}$ converges to a finite value  $\underline{F}$. Note that $\breg[\leg]{\bx\iter\kp}{\bx\iter\k}\to 0$ as $\k\rto{\k\in K}\infty$ due to Proposition~\ref{p:SuffDesc} (ii), when combined with Assumption~\ref{ass:relative-error}(iii) implies that $\vnorm[]{\bx\iter\kp - \bx\iter\k} \to 0$. For $(\opt{\bx},\opt{\bx})\in \omega_{F_{\ub}^\leg}^{(\sint\dom\leg)^2}(\bx\iter0,\bx\iter0)$ there exists $K\subset\N$ such that $\bx\iter\k\rto{\k\in K} \opt{\bx}$ and $F_{\ub}^\leg(\bx\iter\kp,\bx\iter\k) \rto{\k\in K} F_{\ub}^\leg(\opt{\bx},\opt{\bx})= \fun(\opt{\bx})$, i.e., the value of the limit point is independent of the choice of the subsequence. The result follows directly and by using $\ii1$.
\end{proof}

The following result summarizes that $F_{\ub}^\leg$-attentive sequences converge to a stationary point. 
\begin{THM}[Sub-sequential convergence to stationary points] \label{thm:sub-seq-conv}
  Let Assumptions~\ref{ass:problem}, \ref{ass:model}, \ref{ass:relative-error}   hold. If the sequence $\seq[\k\in\N]{\bx\iter\k}$  is generated by Model BPG, then
  \begin{equation} \label{eq:lem:sub-diff-vanish:crit}
    \omega_{F_{\ub}^\leg}^{(\sint\dom\leg)^2}(\bx\iter0) \subset \crit(F_{\ub}^\leg) \,.
  \end{equation}
\end{THM}
\begin{proof}
  From \eqref{eq:lem:sub-diff-vanish}, we have $\vnorm[-]{\partial  F_{\ub}^\leg (\bx\iter\kp, \bx\iter\k)}   \leq  C \vnorm[]{\bx\iter\kp - \bx\iter\k}$ for some constant $C >0$. Using  $\vnorm[]{\bx\iter\kp-\bx\iter\k}\to0$, convergence of $\seq[\k\in\N]{\tau\iter\k}$, and Proposition~\ref{prop:limit-points}$\ii1$ yields \eqref{eq:lem:sub-diff-vanish:crit}, by the closedness property of the limiting subdifferential \eqref{eq:closedness-lim-subdiff}.
  \end{proof}
  \paragraph{Discussion.}{Subsequential convergence to a stationary point was already considered in few works. In particular, the work in \cite{drusvyatskiy2019nonsmooth} already provides such a result, however, it relies on certain abstract assumptions. Even though such assumptions are valid for some practical algorithms, the authors do not consider a concrete algorithm. Moreover, their abstract update step depends on the minimization of the model function, which can require additional regularity conditions on the problem. For example, if the model function  is linear, then the domain must be compact to guarantee the existence of a solution. A related line-search variant of Model BPG  was considered in \cite{OFB19}, for which subsequential convergence to a stationarity point was proven. The subsequential convergence results in \cite{OFB19} are more general than our work, as they analyse the behavior of limit points in $\dom\leg$, $\scl\dom\leg$, $\sint\dom\leg$ (cf. \cite[Theorem 22]{OFB19}). Our analysis is restricted to limit points in $\sint\dom\leg$, as typically such an assumption holds in practice (see Section~\ref{sec:experiments}). Though subsequential convergence is satisfactory, proving global convergence is nontrivial, in general. It is not yet clear from our work, whether global convergence can be proven if the limit points  lie on the boundary of $\dom\leg$. Both the above-mentioned works rely on function values to obtain a subsequential convergence result. We change this trend. In this paper, we rely on Lyapunov function and obtain an even stronger result, that is  global convergence of the sequence generated by Model BPG to a stationarity point.  } 
\subsection{Global convergence to a stationary point of the Lyapunov function}\label{ssec:global-convergence-lyapunov-func}
The global convergence statement of Model BPG relies on the so-called \KL (KL) property. It has became a standard tool in recent years, and it is essentially satisfied by any function that appears in practice, we just state the definition here and refer to \cite{BDL06,BDLS07,ABS13,BST14,Kurd98} for more details. The following definition is from \cite{ABS13}.
\begin{DEF}[\KL property]\label{def:kl-property}
Let $\map f {\X}{\eR}$ be an extended real valued function and let $\bar {\bf x}\in\dom\partial f$. If there exists $\eta\in(0,\infty]$, a neighborhood $U$ of $\bar {\bf x}$ and a continuous concave function $\map{\phi}{[0,\eta)}{\R_+}$ such that 
\[
  \phi(0)=0,\quad  \phi\in C^1 (0,\eta) ,\quad\text{and}\quad \phi^\prime(s)>0\text{ for all }s\in (0,\eta),
\]
and for all $x\in U\cap [f(\bar {\bf x}) < f({\bf x}) < f(\bar {\bf x}) + \eta]$ the \KL inequality
\begin{equation}\label{eq:KL-ineq}
  \phi^\prime(f({\bf x})-f(\bar {\bf x})) \vnorm[-]{\partial f({\bf x})} \geq 1
\end{equation}
holds, then the function has the \KL property at $\bar {\bf x}$. If, additionally, the function is lower semi-continuous and the property holds for each point in $\dom \partial f$, then $f$ is called a \KL function.
\end{DEF} 

We abbreviate \KL property as KL property. The function $\phi$ in the KL property is known as a desingularizing function. Many functions arising in practical problems satisfy the KL property, such as,  for example,  semi-algebraic functions with a desingularizing function of the following form:
\[
  \phi(s) = cs^{1-\theta}\,, 
\]
for certain $c>0$ and $\theta \in [0,1)$. The KL property is crucial in order to prove the global convergence of sequences generated by many algorithms, for example PALM \cite{BST14}, iPALM \cite{PS16}, BPG \cite{BSTV18}, CoCaIn BPG \cite{MOPS2020} and many others.  For the purpose of simplification of analysis, we use the following uniformization lemma for the KL property from \cite{BST14}.

\begin{LEM}[Uniformized KL property {\cite[Lemma 3.6]{BST14}}]\label{lem:uniform-kl-property}
  Let $\Omega$ be a compact set and let $f : \R^N \to \eR$ be proper and lower semicontinuous function. Assume that $f$ is constant on $\Omega$ and satisfies KL property at each point on $\Omega$. Then, there exist $\eps >0$, $\eta >0$,  a continuous concave function $\map{\phi}{[0,\eta)}{\R_+}$ such that 
  \[
    \phi(0)=0,\quad  \phi\in \mathcal{C}^1 (0,\eta) ,\quad\text{and}\quad \phi^\prime(s)>0\text{ for all }s\in (0,\eta),
  \]
  and for all $\bbx \in \Omega$ and $\bx$ in the following intersection
  \[
  \{\bx \in \R^N : \dist(\bx, \Omega) < \eps \} \cap [f(\bbx) < f(\bx) < f(\bbx) + \eta]
  \]
  one has, 
  \begin{equation}\label{eq:uniform-KL-ineq}
    \phi^\prime(f(\bx)-f(\bbx)) \vnorm[-]{\partial f(\bx)} \geq 1
  \end{equation}
\end{LEM}
\medskip

It is well known that the class of functions definable in an o-minimal structure satisfies KL property \cite[Theorem 14]{BDLS07}. The exact definition of o-minimal structure is given in \cite[Definition 6]{BDLS07}, which we record in Section~\ref{sec:definable-functions} in the appendix.  Numerous functions and sets can be defined in an o-minimal structure, for example, sets and functions that are semi-algebraic and globally subanalytic. For a comprehensive discussion, we refer the reader to \cite[Section 4]{BDLS07} and  \cite[Section 4.5]{Ochs15}. 

\begin{ASS}\label{ass:final-kl}
Let $\mathcal{O}$ be an o-minimal structure. The functions $\tilde{f}:\R^N \times \R^N \to \eR\,,\, (\bx,\bbx) \mapsto  \fun(\bx;\bbx)$ with $\dom \tilde{f} :=  \dom\fun \times \dom\fun$, and  $\tilde{h}:\R^N \times \R^N \to \eR\,,\, (\bx,\bbx) \mapsto  h(\bbx) + \scal{\nabla h(\bbx)}{\bx - \bbx}$ with  $\dom\tilde{h} := \dom\leg \times \sint\dom\leg$ are definable $\mathcal{O}$.
  %  Let $\mathcal{O}$ be an o-minimal structure, and also let the following conditions hold:
  % \begin{enumerate}
  %   \item The function  $\tilde{f}:\R^N \times \R^N \to \eR\,,\, (\bx,\bbx) \mapsto  \fun(\bx;\bbx)$ with $\dom \tilde{f} = \dom\fun \times \dom\fun$ is definable in $\mathcal{O}$. 
  %   \item  The function  $\tilde{h}:\R^N \times \R^N \to \eR\,,\, (\bx,\bbx) \mapsto  h(\bbx) + \scal{\nabla h(\bbx)}{\bx - \bbx}$ with $\dom \tilde{h} = \dom\leg \times \sint\dom\leg$ is definable in $\mathcal{O}$.
  % \end{enumerate}
\end{ASS}

The following result shows that functions definable in an o-minimal structure are closed under pointwise addition and multiplication. This is a standard result which can, for example, be found in {\cite[Corollary 4.32]{Ochs15}}.
\begin{LEM}\label{lem:definable-closed}
  Let $S,T \subset \R^M$, $S \cap T = \emptyset$, and let $f:S \to \R^N$, $g:T \to \R^N$ be maps that belong to $\mathcal{O}$. Then, pointwise addition and multiplication, $f+g$ and $f \cdot g$, restricted to $S \cap T$ belongs to  $\mathcal{O}$.
\end{LEM}
The following result connects KL property to functions that are definable in an o-minimal structure. 
\begin{THM}[{\cite[Theorem 14]{BDLS07}}]\label{thm:definable-kl}
  Any proper lower semi-continuous function $f:\R^N \to \eR$ that is definable in an o-minimal structure $\mathcal{O}$ has the \KL property at each point of $\dom\partial \fun$. Moreover the function $\phi$ in Lemma~\ref{lem:uniform-kl-property} is definable in $\mathcal{O}$.
\end{THM}

\begin{LEM}\label{lem:main-lyapunov-kl} 
  Let Assumptions~\ref{ass:problem}, \ref{ass:model}, \ref{ass:relative-error}, \ref{ass:final-kl}  hold. Then, the Lyapunov function $F^{h}_{\ub}$ is definable in $\mathcal{O}$,  and satisfies KL property at any point of $\dom \partial F^{h}_{\ub}$.
\end{LEM}
\begin{proof}
  As per the conditions of Lemma~\ref{lem:definable-closed}, we deduce that functions that are definable in an o-minimal structure are closed under addition and multiplication. With Assumption~\ref{ass:final-kl}, it is easy to deduce that the $F^{h}_{\ub}$ is also definable in $\mathcal{O}$ using Lemma~\ref{lem:definable-closed}. Invoking Theorem~\ref{thm:definable-kl}, we deduce that $F^{h}_{\ub}$ satisfies KL property at any point of $\dom \partial F^{h}_{\ub}$.
\end{proof}

In the context of additive composite problems, the global convergence analysis of BPG based methods \cite{BSTV18,MOPS2020} relies on strong convexity of $h$. However, in our setting we relax such a requirement on $h$, via the following assumption. Note that imposing such an assumption (Assumption~\ref{ass:local-strong-convexity}) is weaker than imposing the strong convexity of $h$, as we only need the strong convexity property to hold over a compact convex set. Such a property can be satisfied even if $h$ is not strongly convex, for example, Burg's entropy (see Section~\ref{ssec:poisson-linear-inverse-problems}). 

\begin{ASS}\label{ass:local-strong-convexity}
  For any compact convex set $B\subset \sint\dom\leg$, there exists $\sigma_B >0$ such that $h$ is $\sigma_B$-strongly convex over $B$, i.e., for any $\bx, \by \in B$ the condition $D_h(\bx, \by)\geq \frac{\sigma_B}{2}\norm{\bx - \by}^2$ holds. 
\end{ASS}
Now, we present the global convergence result of the sequence generated by Model BPG.
\begin{THM}[Global convergence to a stationary point under KL property] \label{thm:full-conv}
  Let Assumptions~\ref{ass:problem}, \ref{ass:model}, \ref{ass:relative-error}, \ref{ass:final-kl}, \ref{ass:local-strong-convexity}  hold.  Let the sequence $\seq[\k\in\N]{\bx\iter\k}$ be generated by Model BPG (Algorithm~\ref{alg:acc-BregMin-bt}) with $\tau\iter\k \to \tau$ for certain $\tau > 0$  and the condition $\omega^{\sint\dom\leg}(\bx\iter0) =\omega(\bx\iter0)$ holds true. Then, convergent subsequences are  $F_{\ub}^\leg$-attentive convergent, and  
  \[
    \sum_{\k=0}^\infty \vnorm[]{\bx\iter\kp - \bx\iter\k} < +\infty \qquad \text{(finite length property)}\,.
  \]
  Moreover, the sequence $\seq[\k\in\N]{\bx\iter\k}$ converges to $\bx$ such that $(\bx,\bx)$ is the critical point of $F_{\ub}^\leg$. 
\end{THM}
\begin{proof}  
  Note that the sequence $\seq[\k\in\N]{\bx\iter\k}$ generated by Model BPG is a bounded sequence (see Remark~\ref{rem:boundedness}). The proof relies on Theorem~\ref{thm:KL-theorem-descent} provided in Section~\ref{sec:gradient-like-descent} in the appendix, for which we need to verify the conditions \ref{ass:Hs:descent}--\ref{ass:Hs:params}. Due to Lemma~\ref{lem:main-lyapunov-kl}, $F_{\ub}^\leg$ satisfies \KL property at each point of $\dom \partial F^{h}_{\ub}$. 
  \medskip

  Note that as  $\omega^{\sint\dom\leg}(\bx\iter0) =\omega(\bx\iter0)$ holds true, there exists a sufficiently small $\eps >0$ such that $\tilde{B} := \{\bx : \dist(\bx, \omega(\bx\iter0))\leq\eps\} \subset \sint\dom\leg$. As $\omega(\bx\iter0)$ is compact due to Proposition~\ref{prop:subsequence}(i), the set $\tilde{B}$ is also compact. Moreover, the convex hull of the set $\tilde{B}$ denoted by $B:=\text{conv} \,\tilde{B}$ is also compact, as the convex hull of a compact set is also compact in finite dimensional setting. A simple calculation reveals that the set $B$ lies in the set $\sint \dom \leg$. Thus, due to Proposition~\ref{prop:limit-point-sets-1} along with Proposition~\ref{prop:subsequence}(ii), without loss of generality, we assume that the sequence $\seq[\k\in\N]{\bx\iter\k}$ generated by Model BPG lies in the set $B$.  By definition of ${\sigma}_B$ as per Assumption~\ref{ass:local-strong-convexity} we have
  \begin{equation}\label{eq:sc-main-kl}
    \breg[\leg]{\bx\iter\kp}{\bx\iter\k} \geq \frac{\sigma_B}{2} \vnorm[]{\bx\iter\kp -\bx\iter\k}^2 \,,
  \end{equation}
  through which we obtain 
  \[
    F_{\ub}^\leg(\bx\iter\kp,\bx\iter\k) \leq F_{\ub}^\leg(\bx\iter\k,\bx\iter\km) - \frac{\eps_k \sigma_B}{2} \vnorm[]{\bx\iter\kp-\bx\iter\k}^2 \,,
  \]
  which is \ref{ass:Hs:descent} with $d_k =  \frac{\eps_k \sigma_B}{2} \vnorm[]{\bx\iter\kp-\bx\iter\k}^2$ and $a_k= 1$. We also have existence of ${\bf w}\iter\kp\in \partial F_{\ub} ^\leg(\bx\iter\kp,\bx\iter\k)$ due to Lemma~\ref{lem:sub-diff-vanish}  such that for some $C>0$ we have
  \[
    \vnorm[-]{\partial  F^h_{\ub} (\bx\iter\kp, \bx\iter\k)} 
    \leq C \vnorm[]{\bx\iter\kp - \bx\iter\k}\,,  
  \]
 which is \ref{ass:Hs:error} with $b =C$, since the coefficients for both Euclidean distances are bounded from above.
  The continuity condition \ref{ass:Hs:cont} is deduced from a converging subsequence, whose existence is guaranteed by boundedness of $\seq[\k\in\N]{\bx\iter\k}$, and Proposition~\ref{prop:limit-points} guarantees that such convergent subsequences are $F_{\ub}^\leg$-attentive convergent.  The distance condition \ref{ass:Hs:distance} holds trivially as $\eps_k >0$ and $\sigma_B >0$. The parameter condition \ref{ass:Hs:params}, holds because $b_n = 1$ in this setting, hence $\seq[\n\in\N]{b\pit\n}\not\in\ell_1$ and also we have
  \[
 \sup_{n\in\N} \frac 1{b\pit\n a\pit\n} =1 < \infty\,, \quad 
 \inf_\n a\pit\n = 1 > 0\,.
  \]
  Theorem~\ref{thm:KL-theorem-descent} implies the finite length property  from which we deduce that the sequence $\seq[\k\in\N]{\bx\iter\k}$ generated by Model BPG converges to a single point, which we denote by $\bx$. As $\seq[\k\in\N]{\bx\iter\kp}$ also converges to $\bx$, the sequence $\seq[\k\in\N]{(\bx\iter\kp, \bx\iter\k)}$  converges to $(\bx,\bx)$, which is a critical point of $F^h_{{\bar L}}$ due to Theorem~\ref{thm:sub-seq-conv}. 
\end{proof}

\subsection{Global convergence to a stationary point of the objective function}\label{ssec:global-convergence-function}
The global convergence result in Theorem~\ref{thm:full-conv} shows that Model BPG converges to a point, which in turn can be used to represent the critical point of the Lyapunov function. However, our goal is to find a critical point of the objective function $f$. We now establish the connection between a critical point of the Lyapunov function and a critical point of the objective function. Such a connection can later be exploited to conclude that the sequence generated by Model BPG converges to a critical point of $f$.
\medskip

Firstly, we need the following result, which establishes the connection between fixed points of the update mapping and critical points of $f$. 
\begin{LEM}\label{lem:critical-fixed-1}
  Let Assumptions~\ref{ass:problem}, \ref{ass:model} hold. For any $0<\tau<({1}/{\ub})$ and $\bbx \in \dom\fun \cap \sint\dom\leg$, the fixed points of the update mapping $T_{\tau}(\bbx)$ are critical points of $f$.
\end{LEM}
\begin{proof}
    Let $\bbx \in \dom\fun \cap \sint\dom\leg$ be a fixed point of $T_{\tau}$, in the sense the condition $ \bbx \in T_{\tau}(\bbx)$  holds true. By definition of $T_{\tau}(\bbx)$, the following condition holds true:
    \begin{equation*}
      {\bf 0} \in \partial \fun(\bx;\bbx) + \frac{1}{\tau}\left(\nabla h(\bx)-\nabla h(\bbx)\right)
    \end{equation*}
    at $\bx = \bbx$, which implies that ${\bf 0} \in \partial \fun(\bbx;\bbx)$. As a consequence of Lemma~\ref{lem:first-order-info}, we have $\partial \fun(\bbx;\bbx) \subset \partial \fun(\bbx)$, thus $\bbx$ is the critical point of the function $f$. 
  \end{proof}
  
  We also require the following technical result.
  % Secondly, the convergence of $\seq[\k\in\N]{\tau\iter\k}$, convergence of the  $\seq[\k\in\N]{\bx\iter\k}$ (valid due to global convergence) and the continuity of the update mapping $T_{\tau\iter\k}(\bx\iter\k)$ in both $\tau\iter\k, \bx\iter\k$ ensures that $T_{\tau\iter\k}(\bx\iter\k)$ converges to a fixed point $\bx_{\infty} := \lim_{k\to\infty}\bx\iter\k$. Note that in the following result, we use arbitrary convergent sequences for $\tau\iter\k$ and $\bx\iter\k$. The result can later be specialized to the sequence generated by Model BPG.

\begin{LEM}[Continuity property]\label{lem:critical-fixed-2}
  Let Assumptions~\ref{ass:problem}, \ref{ass:model}, \ref{ass:relative-error} hold. Let the sequence $\seq[\k\in\N]{\bx\iter\k}$ be bounded such that  $\bx\iter\k \to \bbx$, where $\bx\iter\k \in \dom\fun \cap \sint\dom\leg$ for all $k \in \N$, and $\bbx \in \dom\fun \cap \sint\dom\leg$. Let $\tau\iter\k \to \tau$, such that  $0 < \underline{\tau}  \leq \tau\iter\k \leq {\bar \tau} < {1}/{{\bar L}}$.  Let there exist a bounded set $B \subset \sint\dom\leg$, such that $T_{\tau\iter\k}(\bx\iter\k) \subset B$, $\bx\iter\k \in B$ for all $k \in \N$. If $\limsup_{k \to \infty}T_{\tau\iter\k}(\bx\iter\k) \subset \dom\fun \cap \sint\dom\leg$, then $\limsup_{k \to \infty}T_{\tau\iter\k}(\bx\iter\k) \subset T_{\tau}(\bbx)$.
  % Assume that the limit points of the sequence $\seq[\k\in\N]{\by\iter\k}$ lie in $\dom\fun \cap \sint\dom\leg$. Then, all limit points lie in the set $T_{\tau}(\bbx)$.
  % Then, we have 
  % \[
  % \lim_{k\to \infty} \by\iter\k  = \by \in T_{\tau}(\bbx) \,.
  % \]
  % If limit points of the sequence $\seq[\k\in\N]{T_{\tau\iter\k}(\bx\iter\k)}$ exist, then assume that the limit points of the sequence $\seq[\k\in\N]{T_{\tau\iter\k}(\bx\iter\k)}$ lie in $\dom\fun \cap \sint\dom\leg$. Then, the update mapping $T_{\tau\iter\k}(\bx\iter\k)$ is continuous in $(\bx\iter\k, \tau\iter\k)$, that is  we have 
  % \[
  % \lim_{k\to \infty} T_{\tau\iter\k}(\bx\iter\k)  = T_{\tau}(\bbx) \,.
  % \]
\end{LEM}
\begin{proof}
  Consider any sequence $\seq[\k\in\N]{\by\iter\k}$ such that for any $k\in \N$, the condition $\by\iter\k \in T_{\tau\iter\k}(\bx\iter\k)$ holds true.  Recall that $\fun(\bx ;\by )$ is continuous on its domain due to Assumption~\ref{ass:model}(iv). By optimality of $\by_k \in T_{\tau\iter\k}(\bx\iter\k)$, for any $\bz \in \R^N$ we have the following:
  \begin{equation}\label{eq:optimality-objective}
    \fun(\by\iter\k;\bx\iter\k) + \frac{1}{\tau\iter\k} \breg[\leg]{\by\iter\k}{\bx\iter\k} \leq \fun(\bz;\bx\iter\k) + \frac{1}{\tau\iter\k} \breg[\leg]{\bz}{\bx\iter\k}\,.
  \end{equation}
  % Substituting $\bz = \bx\iter\k$, and on rewriting \eqref{eq:optimality-objective} we have
  % \begin{align*}
  %    \frac{1}{\tau\iter\k} \breg[\leg]{\by\iter\k}{\bx\iter\k} &\leq \fun(\bx\iter\k) - \fun(\by\iter\k;\bx\iter\k) \,,\\
  %    &\leq \fun(\bx\iter\k) - \left(\fun(\by\iter\k;\bx\iter\k) +{\bar L}\breg[\leg]{\by\iter\k}{\bx\iter\k}\right) + {\bar L} \breg[\leg]{\by\iter\k}{\bx\iter\k}\,.
  % \end{align*}
  % Re-arranging the  terms, defining that $\eps\iter\k = \left(\frac{1}{\tau\iter\k}  - {\bar L}\right)$, using $\eps\iter\k \geq {\underline \eps}:= \left(\frac{1}{{\bar \tau}}  - {\bar L}\right)$  and by definition of $\sigma_B$ we obtain
  % \begin{equation}\label{eq:boundedness-cond}
  %   \frac{\sigma_B{\underline \eps}}{2}\vnorm[]{\by\iter\k - \bx\iter\k}^2 \leq \eps\iter\k\breg[\leg]{\by\iter\k}{\bx\iter\k}  \leq \fun(\bx\iter\k) - \fun(\by\iter\k) \leq \fun(\bx\iter\k) - v(\PPP) < \infty\,.
  % \end{equation}
  % As the sequence $\seq[\k\in\N]{\bx\iter\k}$ is bounded, the sequence $\seq[\k\in\N]{\by\iter\k}$ is also bounded, due to \eqref{eq:boundedness-cond}. Moroever, as the sequence $\seq[\k\in\N]{\bx\iter\k}$ lies in $\dom\fun \cap \sint\dom\leg$, we deduce that the sequence $\seq[\k\in\N]{\by\iter\k}$ also lies in $\dom\fun \cap \sint\dom\leg$.  
  As a consequence of boundedness of the sequence  $\seq[\k\in\N]{\by\iter\k}$, by Bolzano–Weierstrass Theorem there exists a convergent subsequence. Let $\by\iter\k \rto{K} \pi$ such that $\pi \in \dom\fun \cap \sint\dom\leg$.  Note that $\tau\iter\k \rto{K} \tau$ for some $ K \subset \N$. Applying limit on both sides of \eqref{eq:optimality-objective} using the continuity of the model function and the Bregman distance gives 
  \begin{equation}\label{eq:subsequence-pi}
    \fun({\bf \pi};\bbx) + \frac{1}{\tau}\breg[\leg]{{\bf \pi}}{\bbx} \leq \fun(\bz;\bbx) + \frac{1}{\tau} \breg[\leg]{\bz}{\bbx}\,,\quad \forall\, \bz \in \dom\fun \cap \dom\leg\,,
  \end{equation}
  which implies that $\pi$ minimizes the function $\fun(\cdot;\bbx) + \frac{1}{\tau}D_h(\cdot,\bbx)$. This implies that $\pi \in T_{\tau}(\bbx)$ and the result follows.
  % This implies that $T_{\tau}(\bbx)$ is continuous in both $\tau$ and $\bbx$. 
\end{proof}

The following result establishes the fact the sequence generated by Model BPG indeed converges to the critical point of the objective function.
\begin{THM}[Global convergence to a stationary point of the objective function]\label{thm:global-conv-obj-func}
  Under the conditions of Theorem~\ref{thm:full-conv},  the sequence generated by Model BPG converges to a critical point of $f$.
\end{THM}
\begin{proof}
  The sequence $\seq[\k\in\N]{\bx\iter\k}$ generated by Model BPG under the assumptions as in Theorem~\ref{thm:full-conv} is globally convergent, thus let $\bx\iter\k \to \bx$ and also  $\bx\iter\kp \to \bx$. As $\bx\iter\kp \in T_{\tau\iter\k}(\bx\iter\k)$ and $\tau\iter\k$ converges to $\tau$, with Lemma~\ref{lem:critical-fixed-2} we deduce that $\bx \in T_{\tau}(\bx)\,.$ Additionally, with the result in Lemma~\ref{lem:critical-fixed-2},  we deduce that $\bx$ is the fixed point of the mapping $ T_{\tau}(\bx)$, i.e., $\bx \in T_{\tau}(\bx)$. Then, using Lemma~\ref{lem:critical-fixed-1} we conclude that $\bx$ is the critical point of the function $f$.
  % Consider a subsequence  $\seq[\k\in K]{\bx\iter\k}$ where $K \subset \N$, we know that $\bx\iter\k \rto{K} \bx$  and similarly we have $\tau\iter\k \rto{K} \tau$, as $\tau\iter\k$ converges to $\tau$. From the continuity property of the update mapping (Lemma~\ref{lem:critical-fixed-2}), as $(\bx\iter\k, \tau\iter\k) \rto{K} (\bx, \tau)$ we have $\bx\iter\kp \rto{K} \btx \in T_{\tau}(\bx)\,.$
  % % \[
  % % T_{\tau\iter\k}(\bx\iter\k) \rto{K} T_{\tau}(\bx)\,.
  % % \]
  % We already know that $\vnorm[]{\bx\iter\kp - \bx\iter\k} \to  0$ as $k \to \infty$, and  $\bx\iter\k \to \bx$  which implies that $\lim_{k \to \infty}\bx\iter\kp = \bx$, thus $\btx = \bx$. With the result in Lemma~\ref{lem:critical-fixed-2},  we deduce that $\bx$ is the fixed point, as $\bx \in T_{\tau}(\bx)$. Then, using Lemma~\ref{lem:critical-fixed-1} we conclude that $\bx$ is the critical point of the function $f$.
\end{proof}

It is possible to deduce convergence rates for a certain class  of desingularizing functions. Based on \cite{AB09, BST14, FGP14}, we provide the following result, which provides the convergence rates for the sequence generated by Model BPG.
\begin{THM}[Convergence rates]
  Under the conditions of Theorem~\ref{thm:full-conv}, let the sequence  $\seq[\k\in\N]{\bx\iter\k}$ generated by Model BPG converge to $\bx \in \dom\fun \cap \sint\dom\leg$, and let the Lyapunov function $F^{h}_{\ub}$ satisfy \KL property with the following desingularizing function:
  \[
    \phi(s) = cs^{1-\theta}\,, 
  \]
  for certain $c >0$ and $\theta \in [0,1)$. Then, we have the following:
  \begin{itemize}
    \item If $\theta = 0$, then $\seq[\k\in\N]{\bx\iter\k}$ converges in finite number of steps.
    \item If $\theta \in (0, \frac12]$, then there exists $\rho \in [0,1)$ and $G > 0$ such that for all $k\geq 0$ we have 
    \[
    \vnorm[]{\bx\iter\k - \bx} \leq G\rho^k\,.  
    \]
    \item If $\theta \in (\frac12,1)$, then there exists $G>0$ such that for all $k\geq 0$ we have
    \[
      \vnorm[]{\bx\iter\k - \bx} \leq G k^{-\frac{1-\theta}{2\theta  -1}}\,.
    \]
  \end{itemize}
\end{THM}
\begin{proof}
  Here, we consider the same notions as in the proof of Theorem~\ref{thm:full-conv}.   First, using the convexity of the function $-s^{1-\theta}$ we obtain
  \begin{align*}
    & (F_{\ub}^\leg(\bx\iter\k,\bx\iter\km) - v(\PPP))^{1-\theta}  - (F_{\ub}^\leg(\bx\iter\kp,\bx\iter\k)- v(\PPP))^{1-\theta} \\
    &\geq (1-\theta)(F_{\ub}^\leg(\bx\iter\k,\bx\iter\km) -v(\PPP))^{-\theta}(F_{\ub}^\leg(\bx\iter\k,\bx\iter\km)  - F_{\ub}^\leg(\bx\iter\kp,\bx\iter\k))\,,\\
    &\geq (1-\theta)(F_{\ub}^\leg(\bx\iter\k,\bx\iter\km) -v(\PPP))^{-\theta}\frac{\eps_k\sigma_B}{2}\vnorm[]{\bx\iter\kp-\bx\iter\k}^2\,,\\
    &\geq (1-\theta)(F_{\ub}^\leg(\bx\iter\k,\bx\iter\km) -v(\PPP))^{-\theta}\frac{{\underline \eps}\sigma_B}{2}\vnorm[]{\bx\iter\kp-\bx\iter\k}^2\,,
  \end{align*}
  where in the second inequality we used the Proposition~\ref{prop:descent-property} along with the definition of $\sigma_B$, and  in the last step we used $\eps_k \geq {\underline \eps}$. Denote $U := \omega_{F_{\ub}^\leg}^{(\sint\dom\leg)^2}(\bx\iter0)$, and thanks to Theorem~\ref{thm:sub-seq-conv} we have $U \subset \crit(F_{\ub}^\leg)$. Due to Proposition~\ref{prop:subsequence}, we already know that $U$ is a connected compact set and
  \[
    \lim_{k \to \infty}\dist\left((\bx\iter\kp,\bx\iter\k), U\right) = 0\,.
  \]
  Continuing the calculation, following the proof technique of \cite[Theorem 3.1]{BST14}, using Lemma~\ref{lem:uniform-kl-property} with $\Omega = U$, we deduce that there exists $l \in \N$, $C_1>0$ such that for any $k > l$, the following holds:
  \[
  \sum_{i = l + 1}^k \vnorm[]{\bx\iter{i+1} - \bx\iter{i}} \leq \vnorm[]{\bx\iter{l+1} - \bx\iter{l}} + C_1(F_{\ub}^\leg(\bx\iter{l+1},\bx\iter{l}) -v(\PPP))^{1-\theta}\,.   
  \]
  
  Denote $\Delta_{l} := \sum_{i=l}^{\infty}\vnorm[]{\bx\iter{i+1} - \bx\iter{i}}$. On application of Lemma~\ref{lem:uniform-kl-property} with $\Omega= U$, and Lemma~\ref{lem:sub-diff-vanish}, we deduce that there exists $C_2 >0$ such that  
  \[
  \Delta_{l+1} \leq \Delta_l - \Delta_{l+1} + C_2 (\Delta_l - \Delta_{l+1})^{\frac{1-\theta}{\theta}}\,.
  \]  
  The rest of the proof is only a slight modification to the proof of \cite[Theorem 5]{AB09}. 
\end{proof}

\section{Examples}\label{sec:examples}
In this section we consider special instances of $(\PPP)$, namely, additive composite problems and a broad class of composite problems. The goal is to quantify assumptions for these problems such that the global convergence result (Theorem~\ref{thm:global-conv-obj-func}) of Model BPG  is applicable. To this regard, we only consider the functions that satisfy Assumption~\ref{ass:problem}. Typically, function is made up of function components and these components govern the function behavior. Thus, it is beneficial to introduce properties on the components of $f$, for which certain plausible conditions will enable the applicability of Model BPG.   In this section, henceforth we enforce the following blanket assumptions.
\begin{enumerate}[label=(B\arabic*),ref=(B\arabic*), leftmargin=*]
  \item\label{ass:Bs:f-0}  The function $h$ is a Legendre function that is $\mathcal{C}^2$ over $\sint\dom\leg$. For any compact convex set $B\subset \sint\dom\leg$, there exists $\sigma_B >0$ such that $h$ is $\sigma_B$-strongly convex over $B$. Also, $h$ has bounded second derivative on any bounded subset $B_1 \subset \sint\dom\leg$. Moreover, for bounded $\seq[\k\in\N]{{\bf u}\iter\k}$, $\seq[\k\in\N]{\bv\iter\k}$ in $\sint\dom\leg$, the following holds as $\k\to\infty$: 
  \[
    \breg[\leg]{{\bf u}\iter\k}{\bv\iter\k} \to 0 \iff
  \vnorm[]{{\bf u}\iter\k - \bv\iter\k} \to 0 \,.
  \]
  \item\label{ass:Bs:f-1} The function $f$ is coercive and additionally the conditions $\dom\fun \cap \sint\dom\leg \neq \emptyset$, $\crit\fun \cap \sint\dom\leg \neq \emptyset$, $\dom f \subset \scl\dom h$ hold true. 
  % , and the following qualification condition holds true:
  % \begin{equation*}
  %   N_{\dom f}(\bx) \cap (-N_{\dom h}(\bx)) = \{{\bf 0}\}\,,\quad \forall\, \bx \in \dom\fun \cap \dom\leg\,. \textcolor{red}{??}
  % \end{equation*}
  \item\label{ass:Bs:f-2} The functions $\tilde{f}:\R^N \times \R^N \to \eR\,,\, (\bx,\bbx) \mapsto  \fun(\bx;\bbx)$ with $\dom \tilde{f} :=  \dom\fun \times \dom\fun$, and  $\tilde{h}:\R^N \times \R^N \to \eR\,,\, (\bx,\bbx) \mapsto  h(\bbx) + \scal{\nabla h(\bbx)}{\bx - \bbx}$ with  $\dom\tilde{h} := \dom\leg \times \sint\dom\leg$ are definable in an o-minimal structure $\mathcal{O}$.
\end{enumerate}
Note that $h$ satisfies Assumption~\ref{ass:Bs:f-0} which considers the same conditions on $h$ as in Assumptions~\ref{ass:model},~\ref{ass:relative-error},~\ref{ass:local-strong-convexity}. The function satisfies Assumption~\ref{ass:Bs:f-1}, which is a consolidation of function specific assumptions in  Assumptions~\ref{ass:problem},\,\ref{ass:model}. Clearly, Assumption~\ref{ass:Bs:f-2} implies Assumption~\ref{ass:final-kl}. 

\subsection{Additive composite problems}\label{ssec:forward-backward}
We consider the following nonconvex additive composite problem:
\begin{align}\label{eq:additive-setting}
  \inf_{\bx \in \R^N}f(\bx)\,, \quad\, f(\bx) := f_0(\bx) + f_1(\bx)\,,
\end{align}
which is a special case of  $(\PPP)$.  Additive composite problems arise in several applications, such as standard phase retrieval \cite{BSTV18}, low rank matrix factorization \cite{MO2019a}, deep linear neural networks \cite{MWLCO2019}, and many more. 
We impose the following conditions that are common in the analysis of forward--backward algorithms \cite{OCBP14}, which are used to optimize additive composite problems. 

\begin{enumerate}[label=(C\arabic*),ref=(C\arabic*), leftmargin=*]
  \item\label{ass:Cs:f-0} $f_0 : \R^N \rightarrow \eR$ is a proper, lsc function and is regular at any $\bx \in \dom f_0$. Also, the following qualification condition holds true:
  \begin{equation}\label{eq:qualification-condition-additive}
    \partial^{\infty}f_0(\bx) \cap (-N_{\dom h}(\bx)) = \{{\bf 0}\}\,,\quad \forall\, \bx \in \dom\fun_0 \cap \dom\leg\,.
  \end{equation}
  % Restricted to its domain $\dom f_0$, the function $f_0$ is continuous and semi-convex with modulus $\alpha \in \R$.
  \item\label{ass:Cs:f-1} $f_1 : \R^N \rightarrow \eR$ is a proper, lsc function and is $\mathcal{C}^2$ on an open set that contains $\dom f_0$. Also, there exist ${\ub}, {\lb}>0$ such that for any $\bbx \in \dom\fun_0\, \cap\, \sint\dom\leg$, the following condition holds true:
  \begin{equation} \label{eq:model-ineq-add}
    -\lb\breg[\leg]{\bx}{\bar{\bx}} \leq  \fun_1(\bx)- \fun_1(\bbx) - \scal{\nabla f_1(\bbx)}{\bx- \bbx} \leq \ub\breg[\leg]{\bx}{\bar{\bx}} \,, \quad \forall\, \bx\in\dom\fun_0 \cap \dom\leg\,.
  \end{equation}
  % \item\label{ass:Cs:f-3} For all $\bx, \by \in \dom\fun$, the condition $({\bf 0}, \bv) \in \partial^{\infty}f(\bx;\by)$ implies $\bv = {\bf 0}$ and $(\bv, {\bf 0}) \in \partial^{\infty}f(\bx;\by)$ implies $\bv = {\bf 0}$  holds true. Moreover, $f(\bx; \by)$ is regular at any $(\bx, \by) \in \dom\fun \times \dom\fun$. 
\end{enumerate}
Note  that with Assumption~\ref{ass:Cs:f-0}, \ref{ass:Cs:f-1} it is easy to deduce that $\dom f_0 = \dom f$. For $\bbx \in \dom\fun$, the model function $f(\cdot; \bbx): \R^N \to \eR$ which, when evaluated at $\bx \in \dom\fun$  gives
\begin{equation}\label{eq:model-function-additive}
  \fun(\bx;\bbx) :=  f_0(\bx) + f_1(\bbx) +\scal{\nabla f_1(\bbx)}{\bx - \bbx}\,.
\end{equation}
Using the model function in \eqref{eq:model-function-additive} and the condition \eqref{eq:model-ineq-add}, we deduce that there exist $\lb,\ub >0$ such that for any $\bbx \in \dom\fun \cap \sint\dom\leg$, MAP property is satisfied at $\bbx$ with $\lb,\ub$  as the following  holds true:
\begin{equation} \label{eq:model-ineq-add-1}
  -\lb\breg[\leg]{\bx}{\bar{\bx}} \leq  \fun(\bx)- \fun(\bx; \bbx) \leq \ub\breg[\leg]{\bx}{\bar{\bx}} \,, \quad \forall\, \bx\in\dom\fun \cap \dom\leg\,,
\end{equation}
as  $\fun(\bx)- \fun(\bx; \bbx)  := f_1(\bx) - f_1(\bbx) -\scal{\nabla f_1(\bbx)}{\bx - \bbx}$, thus satisfying Assumption~\ref{ass:model}(i).  The condition in \eqref{eq:model-ineq-add-1} is similar to the popular $L$-smad property in \cite{BSTV18}. The main addition is that $\bx \in \dom\fun \cap \dom\leg$ and  $\bbx \in \dom\fun \cap \sint\dom\leg$, whereas the $L$-smad property requires  $\bx, \bbx \in \dom\fun \cap \sint\dom\leg$. We illustrate this below.
\medskip

\textbf{Remark.} Consider  $f_1(x):= 0.5x^2$, $f_0(x):= I_{[0,\infty)}(x)$ and $h(x)=x\log(x)$ with $\dom\leg = [0,\infty)$ under $0\log(0) = 0$. Clearly, $\dom\leg \subset \dom\fun_1$ and $\dom\fun \subset \dom\leg$  hold true. The function $f_1$ is differentiable at $x = 0$, and MAP condition in \eqref{eq:model-ineq-add} holds true for $x=0$. This scenario is not considered in the $L$-smad property (see \cite[Lemma 2.1]{BSTV18}).  
\medskip

We present below Model BPG algorithm that is applicable for additive composite problems. Using the model function in \eqref{eq:model-function-additive} in Model BPG we recover the BPG algorithm from \cite{BSTV18}.
\medskip

{\centering\,
		\fcolorbox{black}{Orange!10}{\parbox{0.98\textwidth}{
      \text{\textbf{BPG} is Model BPG (Algorithm~\ref{alg:acc-BregMin-bt}) with  }  
      \begin{equation}\label{eq:update-step-forward-backward-temp}
        \fun(\bx; \bx\iter\k) := f_0(\bx) + f_1(\bx\iter\k) + \scal{\nabla f_1(\bx\iter\k)}{\bx - \bx\iter\k}\,.
      \end{equation}\vspace{-1.3em}
		}}
}
\medskip

For $h(\bx) = \frac12\vnorm[]{{\bf x}}^2$, Model BPG is equivalent to proximal gradient method. Assumptions \ref{ass:Cs:f-0}, \ref{ass:Cs:f-1} along with  \ref{ass:Bs:f-1}  imply proper, lsc property of $f$ and lower-boundedness of $f$, thus satisfying  Assumption~\ref{ass:problem}. 
Considering \ref{ass:Cs:f-0} we deduce that $f_0(\bx)$ is regular at $\bx \in \dom f_0$. Using  \cite[Proposition 10.5]{Rock98} we note that $f_0(\bx)$  is regular at all $(\bx,\bbx) \in \dom\fun \times \dom\fun$. Let $(\bx, \bbx) \in \dom\fun \times\dom\fun$, using \cite[Proposition 10.5]{Rock98} on $f_0$, we obtain the following result:
\begin{align}
  & \partial_{(\bx,\by)} f_0(\bx) = (\partial_{\bx} f_0(\bx) ,{\bf 0})\,,\quad 
\partial_{(\bx,\by)}^{\infty} f_0(\bx) = (\partial_{\bx}^{\infty} f_0(\bx) ,{\bf 0})\,.\label{eq:subdiffs}
\end{align}
Let $(\bx, \bbx) \in \dom\fun \times\dom\fun$, we consider the following entity:
\[
	\partial \fun(\bx;\bbx) = \partial (f_0(\bx) + f_1(\bbx) + \scal{\nabla f_1(\bbx)}{ \bx - \bbx})\,,
\]
and in order for the summation rule of subdifferential (\cite[Corollary 10.9]{Rock98}) to be applicable at $(\bx, \bbx)$, we need finiteness of $f_0(\bx)$ and continuously differentiability of $\tilde{f}_1(\bx, \bbx) := f_1(\bbx) + \scal{\nabla f_1(\bbx)}{\bx - \bbx}$ (also see \cite[Exercise 8.8]{Rock98}). Clearly, $f_0$ is finite at $(\bx, \bbx)$, and  $\tilde{f}_1$ is  finite and also continuously differentiable around $(\bx, \bbx)$ due to Assumption~\ref{ass:Cs:f-1}. Thus, using \eqref{eq:subdiffs} and \cite[Corollary 10.9]{Rock98} we obtain the following conditions:
\begin{equation}\label{eq:main-func-identities}
  \partial \fun(\bx;\bbx) =  (\partial_{\bx} f_0(\bx) + \nabla f_1(\bbx),  \nabla^2 f_1(\bbx)(\bx-\bbx))\,,\quad   \partial^{\infty} \fun(\bx;\bbx) = (\partial_{\bx}^{\infty} f_0(\bx),  {\bf 0})\,,
\end{equation}
% \[
%   \partial \fun(\bx;\bbx) =  (\partial_{\bx} f_0(\bx) + \nabla f_1(\bbx),  \nabla^2 f_1(\bbx)(\bx-\bbx))\,,\quad   \partial^{\infty} \fun(\bx;\bbx) = (\partial_{\bx}^{\infty} f_0(\bx),  {\bf 0})\,,
% \]
and as a result \eqref{eq:main-seperable-remark} is satisfied. Using the condition \eqref{eq:qualification-condition-additive} and \eqref{eq:main-func-identities}, we deduce that  Assumption~\ref{ass:model}(ii) is satisfied.  Now, we verify Assumption~\ref{ass:relative-error}(i). Consider a bounded subset $S$ in $\dom\fun$.  For fixed $\bx \in \dom\fun$, and for all $\bbx \in S$ we have
\begin{equation}\label{eq:second-argument}
  \partial_\bbx \fun(\bx;\bbx) = \{\nabla_{\bbx}(\fun(\bx;\bbx))\} = \{\nabla^2 f_1(\bbx)(\bx-\bbx)\}.
\end{equation}
Note that $\nabla f_1$ is Lipschitz continuous on any bounded subset of $\dom\fun$, as $f_1$ is $\mathcal{C}^2$ on $\dom\fun$. This implies that the Hessian is bounded on bounded sets of $\dom\fun$. Thus, based on the same notions in \eqref{eq:second-argument}, we deduce that there exists a constant $M>0$ such that
\begin{align*}
\vnorm[]{\nabla_{\bbx}(\fun(\bx;\bbx))} \leq M \vnorm[]{\bx-\bbx}\,,
\end{align*}
holds true, thus verifying Assumption~\ref{ass:relative-error}(i). 
% Assumption~\ref{ass:Cs:f-3} implies Assumption~\ref{ass:model}(iv). \textcolor{red}{TODO}.  Moreover, semi-convexity implies that the model function is regular \cite[Example 7.27, Exercise 8.20]{Rock98}. Additionally, 
As a simple consequence of Assumption~\ref{ass:Cs:f-0},~\ref{ass:Cs:f-1} the condition Assumption~\ref{ass:model}(iv) is satisfied.
\medskip

As  discussed above,  Assumptions~\ref{ass:Cs:f-0}, \ref{ass:Cs:f-1}, \ref{ass:Bs:f-0}, \ref{ass:Bs:f-1}, \ref{ass:Bs:f-2} imply Assumptions~\ref{ass:problem}, \ref{ass:model}, \ref{ass:relative-error}, \ref{ass:final-kl}, \ref{ass:local-strong-convexity}.  Thus, as a consequence of Theorem~\ref{thm:full-conv}, \ref{thm:global-conv-obj-func} we obtain the following result which provides the global convergence of the sequence generated by BPG to a stationary point. 

\begin{THM}[Global convergence of BPG sequence]\label{thm:full-conv-additive}
  Let Assumptions~\ref{ass:Cs:f-0}, \ref{ass:Cs:f-1}, \ref{ass:Bs:f-0}, \ref{ass:Bs:f-1}, \ref{ass:Bs:f-2} hold. Let the sequence $\seq[\k\in\N]{\bx\iter\k}$ be generated by BPG and the condition $\omega^{\sint\dom\leg}(\bx\iter0) =\omega(\bx\iter0)$ holds true. Let $\tau\iter\k \to \tau$ for certain $\tau > 0$.   Then, the sequence $\seq[\k\in\N]{\bx\iter\k}$ has finite length, that is 
  \[
    \sum_{\k=0}^\infty \vnorm[]{\bx\iter\kp - \bx\iter\k} < +\infty \,,
  \]
  and the sequence $\seq[\k\in\N]{\bx\iter\k}$ converges to $\bx$, which is a critical point of $f$.
\end{THM}
\subsection{Composite problems}\label{ssec:general-composite-problems}
We consider the following nonconvex composite problem:
\begin{align}\label{eq:composite-setting}
\inf_{\bx \in \R^N}f(\bx)\,,\quad\, f(\bx) := f_0(\bx) + g(F(\bx))\,,
\end{align}
which is a special case of the problem $(\PPP)$. Composite problems arise in robust phase retrieval, robust PCA, censored $\Z_2$ synchronization \cite{D2017, LW2016, N2007, DL18, DP2019}. We require the following conditions.
\begin{enumerate}[label=(D\arabic*),ref=(D\arabic*), leftmargin=*]
  \item\label{ass:Ds:f-0} $f_0 : \R^N \rightarrow \eR$ is a proper, lsc function and is regular at any $\bx \in \dom f_0$. Also, the following qualification condition holds true:
  \begin{equation}\label{eq:qualification-condition-comp}
    \partial^{\infty}f_0(\bx) \cap (-N_{\dom h}(\bx)) = \{{\bf 0}\}\,,\quad \forall\, \bx \in \dom\fun_0 \cap \dom\leg\,.
  \end{equation}
  % Restricted to its domain $\dom f_0$, the function $f_0$ is continuous and semi-convex with modulus $\alpha_1 \in \R$.
  \item\label{ass:Ds:f-1} $g : \R^M \rightarrow \R$ is a Q-Lipschitz continuous function and a regular function. Also, there exists $P >0$ such that at any $\bx \in \R^M$, the following condition holds true:
  \begin{equation}\label{eq:bounded-subgradients}
    \sup_{\bv \in \partial g(\bx)}\norm{\bv} \leq P\,.  
  \end{equation}
  \item\label{ass:Ds:f-2} $F : \R^N \rightarrow \R^M$ is  $\mathcal{C}^2$ over $\R^N$. Also, there exist $L>0$ such that for any $\bbx \in \dom\fun_0\, \cap\, \sint\dom\leg$, the following condition holds true:
  \[
  \vnorm[]{F(\bx) - F(\bbx) - \nabla F(\bbx)(\bx - \bbx)} \leq   L D_h(\bx,\bbx)\,, \quad  \forall\, \bx \in \dom f_0 \cap \dom \leg\,,
  \]
  where $\nabla F(\bbx)$ is the Jacobian of $F$ at $\bbx$. 
  % \item\label{ass:Ds:f-3} There exist $\alpha_2 \in \R$ such that for any $\bbx \in \dom f \cap \sint\dom\leg$, the function $g(F(\bbx) + \nabla F(\bbx)(\cdot - \bbx))$ is semi-convex with modulus $\alpha_2$. 
  % \item\label{ass:Ds:f-3-a} For any fixed $\bbx \in \dom\fun_0 \cap \dom\leg$, the following qualification condition holds true:
  % \begin{equation}\label{eq:qual-0}
  %   \left(\partial^{\infty}f_0(\bx) + \nabla F(\bbx)^*\partial^{\infty}_{F(\bbx) + \nabla F(\bbx)(\bx - \bbx)}g(F(\bbx) + \nabla F(\bbx)(\bx - \bbx))\right) \, \cap\, (-N_{\dom h}(\bx)) = \{{\bf 0}\}\,,
  % \end{equation}
  % for all $\bx \in \dom\fun_0 \,\cap\, \dom\leg\,.$
  % \item\label{ass:Ds:f-3} For $(\bx,\bbx) \in \dom\fun_0 \times \dom\fun_0$ with $F(\bx;\bbx) := F(\bbx) + \nabla F(\bbx)(\bx - \bbx)$, the following qualification conditions hold true:
  % \begin{equation}\label{eq:qual-1}
  %   \partial^{\infty}_{(\bx,\bbx)} f_0(\bx) \cap (-\partial^{\infty} g(F(\bx;\bbx)))) = \{({\bf 0},{\bf 0})\}\,,
  % \end{equation}
  % and the only $y$ such that 
  % \begin{equation}\label{eq:qual-2}
  %   y \in \partial^{\infty }_{F(\bx;\bbx)} g(F(\bx;\bbx))\text{ with }(\nabla F(\bbx)^{\ast}y, (\nabla F(\bbx)+\nabla_{\bbx}(\nabla F(\bbx)(\bx-\bbx)))^{\ast}y) = ({\bf 0},{\bf 0})\text{ is }y={0}\,,
  % \end{equation}
  % where $\nabla_{\bbx}(\nabla F(\bbx)(\bx-\bbx))$ denotes the Jacobian of the mapping $\nabla F(\bbx)(\bx-\bbx)$ at $\bbx \in \dom f_0$ with fixed $\bx \in \dom f_0$. 
\end{enumerate}

Note that when $M = 1$, $g(x) = x$, the problem in \eqref{eq:composite-setting} is a special case of \eqref{eq:additive-setting}. However, for a generic $g$ satisfying \ref{ass:Ds:f-1}, the problem in \eqref{eq:composite-setting} cannot be captured under the additive composite problem setting given in Section~\ref{ssec:forward-backward}. Thus, in this section we consider a separate analysis for generic composite problems in \eqref{eq:composite-setting}.
\medskip

The properties \ref{ass:Ds:f-0}, \ref{ass:Ds:f-1}, \ref{ass:Ds:f-2} along with \ref{ass:Bs:f-1}  imply proper, lsc property and lower-boundedness of $f$, thus satisfying Assumption~\ref{ass:problem}. Note  that with Assumption~\ref{ass:Ds:f-0}, \ref{ass:Ds:f-1}, \ref{ass:Ds:f-2} it is easy to deduce that $\dom f_0 = \dom f$. Let $\bbx \in \dom\fun$ and we consider the following model function which, when evaluated at $\bx \in \dom\fun$ gives:
\begin{equation}\label{eq:prox-linear-model}
  \fun(\bx;\bbx) = f_0(\bx) + g(F(\bbx) + \nabla F(\bbx)(\bx - \bbx))\,.
\end{equation} 
% As per \ref{ass:Ds:f-0}, the function  $f_0(\cdot) - \alpha_1\frac{\vnorm[]{\cdot}^2}{2}$ is convex. As per \ref{ass:Ds:f-3}, for any $\bbx \in \dom f \cap \dom \leg$, the function $g(F(\bbx) + \nabla F(\bbx)(\cdot - \bbx)) - \alpha_2\frac{\vnorm[]{\cdot}^2}{2}$ is convex. Thus, the semi-convexity condition in Assumption~\ref{ass:model}(iii) is satisfied, as $f_0$ is semi-convex and the function $g(F(\bbx) + \nabla F(\bbx)(\cdot - \bbx))$ is semi-convex, and the semi-convexity modulus of $f(\cdot;\bbx)$ is $\alpha := \alpha_1 + \alpha_2$. 
% \begin{REM}
%   Note that if we exclude the semi-convexity assumption in \ref{ass:Ds:f-3} and consider just $g$ to be semi-convex with certain modulus $\tilde{\alpha} >0$, then for a given $\bbx \in \dom f \cap \sint\dom\leg$ the semi-convexity modulus of $g(F(\bbx) + \nabla F(\bbx)(\cdot - \bbx))$ depends on $\bbx$. Thus, in such a scenario the Assumption~\ref{ass:model}(iii) can potentially be violated.
% \end{REM}
% \medskip

Using \ref{ass:Ds:f-1}, \ref{ass:Ds:f-2} we deduce that there exists ${\bar L} := LQ >0$ such that for any $\bbx \in \dom f \cap \sint\dom h$, the following MAP property holds at $\bbx$ with ${\bar L}$:
\begin{align*}
|\fun(\bx) - \fun(\bx;\bbx)| &= |g(F(\bx)) - g(F(\bbx) + \nabla F(\bbx)(\bx - \bbx))|\,\leq  {\bar L} D_h(\bx, \bbx) \,,
\end{align*}
for all $\bx \in  \dom f\, \cap\, \dom h$,  as $g$ is $Q$-Lipschitz continuous and \ref{ass:Ds:f-2} holds true. Thus, Assumption~\ref{ass:model}(i) is satisfied with ${\bar L} = {\underline L} = LQ$. Before we verify other assumptions, we present Prox-Linear BPG, a specialization of Model BPG that is applicable to composite problems.
\medskip

{\centering\,
		\fcolorbox{black}{Orange!10}{\parbox{0.98\textwidth}{

      \textbf{Prox-Linear BPG} is Model BPG (Algorithm~\ref{alg:acc-BregMin-bt})  with 
        \begin{equation}\label{eq:update-step-prox-linear-temp}
          \fun(\bx; \bx\iter\k):= f_0(\bx) + g(F(\bx\iter\k) + \nabla F(\bx\iter\k)(\bx - \bx\iter\k))\,. 
				\end{equation}\vspace{-1.5em}
		}}
  }
\medskip

  For $h(\bx) = \frac12 \vnorm[]{\bx}^2$, Prox-Linear BPG is related to Prox-Linear method \cite{LW2016, DL18}. Considering \ref{ass:Ds:f-0}, we deduce that $f_0(\bx)$ is regular at $\bx \in \dom f_0$. Using  \cite[Proposition 10.5]{Rock98} we note that $f_0(\bx)$  is regular at all $(\bx,\bbx) \in \dom\fun \times \dom\fun$. Using \cite[Theorem 10.6]{Rock98} and \ref{ass:Ds:f-1} we deduce that $g(F(\bbx) + \nabla F(\bbx)(\bx - \bbx))$ is regular for all $(\bx,\bbx) \in \R^N \times \R^N$. Furthemore, as a consequence of  \cite[Corollary 10.9]{Rock98}, the function $f(\bx;\bbx)$ is regular at $(\bx, \bbx) \in \dom\fun \times \dom\fun$.  
\medskip

Using \cite[Proposition 10.5]{Rock98} on $f_0$ we deduce that for all $(\bx, \bbx) \in \dom\fun \times \dom\fun$, the following conditions hold true:
\begin{equation}\label{eq:additive-composite-2}
  \partial_{(\bx,\bbx)} f_0(\bx) = (\partial_{\bx} f_0(\bx) ,{\bf 0})\,,\quad 
\partial_{(\bx,\bbx)}^{\infty} f_0(\bx) = (\partial_{\bx}^{\infty} f_0(\bx) ,{\bf 0})\,.
\end{equation} 
For this section, henceforth,  we set $(\bx, \bbx) \in \dom\fun \times \dom\fun$ and denote $F(\bx;\bbx) := F(\bbx) + \nabla F(\bbx)(\bx - \bbx)$. Note that as $\partial^{\infty }_{F(\bx;\bbx)} g(F(\bx;\bbx)) = \{{\bf 0}\}$ due to \ref{ass:Ds:f-0} and \cite[Theorem 9.13]{Rock98}, we deduce that the only $y$ such that 
  \begin{equation}\label{eq:qual-2}
    y \in \partial^{\infty }_{F(\bx;\bbx)} g(F(\bx;\bbx))\text{ with }(\nabla F(\bbx)^{\ast}y, (\nabla F(\bbx)+\nabla_{\bbx}(\nabla F(\bbx)(\bx-\bbx)))^{\ast}y) = ({\bf 0},{\bf 0})\text{ is }y={0}\,,
  \end{equation}
  where $\nabla F(\bbx)^{\ast}$ denotes the adjoint of $\nabla F(\bbx)$, and $\nabla_{\bbx}(\nabla F(\bbx)(\bx-\bbx))$ denotes the Jacobian of the mapping $\nabla F(\bbx)(\bx-\bbx)$ at $\bbx$ with fixed $\bx$. Due to \ref{ass:Ds:f-2},  regularity of $g$ and \eqref{eq:qual-2} we have 
\begin{align*}
  &\partial g(F(\bx;\bbx))  = \left(\nabla F(\bbx)^{*}\partial_{F(\bx;\bbx)} g(F(\bx;\bbx)), (\nabla F(\bbx)+\nabla_{\bbx}(\nabla F(\bbx)(\bx-\bbx)))^{\ast}\partial_{F(\bx;\bbx)} g(F(\bx;\bbx))\right)\,.
\end{align*}
% \begin{align*}
%   \partial^{\infty} g(F(\bbx) + \nabla F(\bbx)(\bx - \bbx)) = \{{\bf 0},{\bf 0}\}\,.
% \end{align*}
A similar statement also holds for $\partial^{\infty} g(F(\bx;\bbx))$ which on using $\partial^{\infty }_{F(\bx;\bbx)} g(F(\bx;\bbx)) = \{{\bf 0}\}$ due to \cite[Theorem 9.13]{Rock98} results in $\partial^{\infty} g(F(\bx;\bbx)) = \{{\bf 0},{\bf 0}\}$. This further implies that the following qualification condition holds true:
\begin{equation}\label{eq:qual-1}
  \partial^{\infty}_{(\bx,\bbx)} f_0(\bx) \cap (-\partial^{\infty} g(F(\bx;\bbx)))) = \{({\bf 0},{\bf 0})\}\,.
\end{equation}
Using the  qualification condition \eqref{eq:qual-1} along with \cite[Corollary 10.9]{Rock98}, we obtain the following:
\begin{equation}\label{eq:additive-composite-1}
  \partial f(\bx,\bbx) =  \partial_{(\bx,\bbx)} f_0(\bx)  + \partial g(F(\bx;\bbx))))\,,\quad  \partial^{\infty} \fun(\bx;\bbx) = (\partial_{\bx}^{\infty} f_0(\bx),  {\bf 0})\,.
\end{equation}
  Thus, \eqref{eq:main-seperable-remark} is satisfied. Additionally, using the condition \eqref{eq:qualification-condition-comp} in \ref{ass:Ds:f-0}, we deduce that  Assumption~\ref{ass:model}(ii) is satisfied.   Now, we verify Assumption~\ref{ass:relative-error}(i). Let's consider a bounded subset $S$ in $\dom\fun$. For $\bbx \in \dom\fun$, there exists a constant $M_S>0$ (dependent on $S$) such that for all ${\bf w} \in \partial_{\bbx} \fun(\bx;\bbx) :=  (\nabla F(\bbx)+\nabla_{\bbx}(\nabla F(\bbx)(\bx-\bbx)))^{*}\partial_{F(\bx;\bbx)} g(F(\bx;\bbx))$ the following condition holds true:
\begin{align}
\vnorm[]{{\bf w}} \leq M_S \vnorm[]{\bx - \bbx}\,,\quad  \forall\, \bx \in S\,,
\end{align}
where we have used the boundedness of second order derivatives of components of $F$ over $S$, as $F$ is a twice continuously differentiable mapping, and boundedness of subgradients of $g$ as per \eqref{eq:bounded-subgradients}. As a simple consequence of Assumption~\ref{ass:Ds:f-0},~\ref{ass:Ds:f-1},~\ref{ass:Ds:f-2} the condition Assumption~\ref{ass:model}(iv) is satisfied. 
\medskip

As  discussed above,  Assumptions~\ref{ass:Ds:f-0}, \ref{ass:Ds:f-1}, \ref{ass:Ds:f-2},  \ref{ass:Bs:f-0}, \ref{ass:Bs:f-1}, \ref{ass:Bs:f-2} imply Assumptions~\ref{ass:problem}, \ref{ass:model}, \ref{ass:relative-error}, \ref{ass:final-kl}, \ref{ass:local-strong-convexity}.  Thus, as a consequence of Theorem~\ref{thm:full-conv}, \ref{thm:global-conv-obj-func} we obtain the following result which provides the global convergence of the sequence generated by Prox-Linear BPG to a stationary point. 
\begin{THM}[Global convergence of Prox-Linear BPG sequence]\label{thm:full-conv-composite}
  Let Assumptions~\ref{ass:Ds:f-0}, \ref{ass:Ds:f-1}, \ref{ass:Ds:f-2}, \ref{ass:Bs:f-0}, \ref{ass:Bs:f-1}, \ref{ass:Bs:f-2} hold. Let the sequence $\seq[\k\in\N]{\bx\iter\k}$ be generated by Prox-Linear BPG  and the condition $\omega^{\sint\dom\leg}(\bx\iter0) =\omega(\bx\iter0)$ holds true. Let $\tau\iter\k \to \tau$ for certain $\tau > 0$. Then, the sequence $\seq[\k\in\N]{\bx\iter\k}$ has finite length, that is 
  \[
    \sum_{\k=0}^\infty \vnorm[]{\bx\iter\kp - \bx\iter\k} < +\infty\,,
  \]
  and the sequence $\seq[\k\in\N]{\bx\iter\k}$ converges to $\bx$, which is a critical point of $f$.
\end{THM}

\section{Experiments}\label{sec:experiments}
For the purpose of empirical evaluation we consider many practical problems, namely, standard phase retrieval problems, robust phase retrieval problems and Poisson linear inverse problems. We compare our algorithms with Inexact Bregman Proximal Minimization Line Search (IBPM-LS) \cite{ODBP13}, which is a popular algorithm to solve generic nonsmooth nonconvex problems. Before we provide the empirical results, we comment below on a variant of Model BPG based on the backtracking technique, which we used in the experiments.

\paragraph{Model BPG with backtracking.} It is possible that the value of ${\bar L}$ in the MAP property is unknown. This issue can be solved by using a backtracking technique, where in each iteration a local constant ${\bar L}\iter\k$ is found such that the following condition holds:
\begin{equation}\label{eq:upper-bound}
  f(\bx\iter\kp) \leq f(\bx\iter\kp;\bx\iter\k) + {\bar L}\iter\k D_h(\bx\iter\kp, \bx\iter\k)\,.
\end{equation}
The value of ${\bar L}\iter\k$ is found by taking an initial guess ${\bar L}\iter\k^0$. If the condition \eqref{eq:upper-bound} fails to hold, then with a scaling parameter $\nu > 1$, we set ${\bar L}\iter\k$ to the smallest value in the set $\{\nu {\bar L}\iter\k^0,\nu^2 {\bar L}\iter\k^0, \nu^3 {\bar L}\iter\k^0,\ldots \}$ such that \eqref{eq:upper-bound} holds true.
Enforcing ${\bar L}\iter\k \geq {\bar L}\iter\km$ for $k\geq 1$ ensures that after finite number of iterations there is no change in the value of ${\bar L}\iter\k$, which takes us to the situation that we analyzed in the paper. The condition ${\bar L}\iter\k \geq {\bar L}\iter\km$ can be  enforced by choosing ${\bar L}\iter\k^0 = {\bar L}\iter\km$.

\subsection{Standard phase retrieval}
The phase retrieval problem involves approximately solving a system of quadratic equations. Let $b_i \in \R$ and ${\bf A}_i \in \R^{N\times N}$  be a symmetric positive semi-definite matrix, for all $i = 1,\ldots,M$.  The goal of standard phase retrieval problem is to find $\bx \in \R^N$ such that the following system of quadratic equations is satisfied:
\begin{equation}\label{eq:robust-phase-retrieval}
  \bx^T{\bf A}_i\bx \approx b_i, \quad \text{ for }i = 1,\ldots,M.
\end{equation}
In standard terminology, $b_i$'s are measurements and ${\bf A}_i$'s are so-called sampling matrices. In the context of Bregman proximal algorithms,  regarding the phase retrieval problem, we refer the reader to \cite{BSTV18, MOPS2020}. Further references regarding the phase retrieval problem include \cite{CLS2015, WGE2018,L2017}. The standard technique to solve such system of quadratic equations is to  solve the following optimization problem:
\begin{equation}\label{eq:phase-retrieval}
  \min_{\bx \in \R^N} \mathcal{P}_0(\bx)\,, \quad \mathcal{P}_0(\bx) := \frac{1}{M}\sum_{i=1}^M{(\bx^T{\bf A}_i\bx - b_i)^2} + \mathcal{R}(\bx)\,,
\end{equation}
where $\mathcal{R}(\bx)$ is the regularization term. We consider here L1 regularization with $\mathcal{R}(\bx) =\lambda \vnorm[1]{\bx}$ and squared L2 regularization with $\mathcal{R}(\bx) = \frac{\lambda}{2}\vnorm[]{\bx}^2$, with some $\lambda >0$. We consider two model functions in order to solve the problem in \eqref{eq:phase-retrieval}.

\paragraph{Model 1.}  Here, the analysis falls under the category of additive composite problems given in Section~\ref{ssec:forward-backward}, where we set the following:
\[
f_0(\bx):= \mathcal{R}(\bx)\,, \text{ and } \quad f_1(\bx):= \frac{1}{M}\sum_{i=1}^M{(\bx^T{\bf A}_i\bx - b_i)^2} \,.
\]
We consider the standard model for additive composite problems from \cite{BSTV18}, where around $\by \in \R^N$, the model function ${\mathcal{P}_0}(\cdot;\by): \R^N  \to \R$  at $\bx \in \R^N$ is given by 
\begin{equation}
  {\mathcal{P}_0}(\bx;\by) := \frac{1}{M}\sum_{i=1}^M\left((\by^T{\bf A}_i\by  - b_i)^2 + (\by^T{\bf A}_i\by  - b_i)\scal{2{\bf A}_i\by}{\bx-\by}\right)  + \mathcal{R}(\bx) \,.\label{eq:p_0_1}
\end{equation}
Consider the following Legendre function:
\[
  h(\bx) = \frac{1}{4}\vnorm[]{\bx}^4 + \frac12 \vnorm[]{\bx}^2\,.
\]
Then, due to \cite[Lemma 5.1]{BSTV18} the following $L$-smad property or the MAP property is satisfied :
\begin{equation}
  \abs{{\mathcal{P}_0(\bx) - {\mathcal{P}_0}(\bx;\by)} } 
  \leq L_0D_h(\bx, \by)\,,\text{ for all }\bx,\by \in \R^N,
\end{equation}
where $L_0\geq \sum_{i=1}^M(3\vnorm[F]{{\bf A}_i}^2 + \vnorm[F]{{\bf A}_i}\abs{b_i})$. In this setting, Model BPG subproblems have closed form solutions (see \cite{BSTV18,MOPS2020}).
\medskip

\paragraph{Model 2.} The importance of finding better models suited to a particular problem was emphasized in \cite{AD2019}. The above provided model function in \eqref{eq:p_0_1} is satisfactory, however, we would like take advantage of the structure of the function \eqref{eq:phase-retrieval}. Taking inspiration from \cite{AD2019}, a simple observation that the objective is nonnegative can be exploited to create a new model function. We incorporate such a behavior in our second model function provided below. We use the Prox-Linear setting described in Section~\ref{ssec:general-composite-problems}, where for any $\bx \in \R^N$ we set the following: 
\[
  f_0(\bx) := \mathcal{R}(\bx)\,,
\]
\[
(F(\bx))_i =   (\bx^T{\bf A}_i\bx - b_i)^2\,, \text{ for all } i= 1,\ldots, M\,,
\]
and for any $\tilde{\by} \in \R^M$ we set
\[
  g(\tilde{\by}) := \frac{1}{M}\vnorm[1]{\tilde{\by}}  \,, \text{ for } \tilde{\by} \in \R^M\,.
\]
Based on the model function \eqref{eq:prox-linear-model}, for fixed $\by \in \R^N$, we consider the model function $\mathcal{P}_1(\cdot;\by) : \R^N \to \R$ which, when evaluated at $\bx \in \R^N$ gives
\begin{equation}\label{eq:p_0_1-2}
  \mathcal{P}_1(\bx;\by) := \frac{1}{M}\sum_{i=1}^M\abs{(\by^T{\bf A}_i\by  - b_i)^2 + (\by^T{\bf A}_i\by  - b_i)\scal{2{\bf A}_i\by}{\bx-\by}}  + \mathcal{R}(\bx) \,.
\end{equation}
Considering the Legendre function $h(\bx) = \frac{1}{4}\vnorm[]{\bx}^4 + \frac12 \vnorm[]{\bx}^2$ and \cite[Lemma 5.1]{BSTV18}, a simple calculation reveals that the following MAP property holds true:
\begin{equation}
  \abs{{\mathcal{P}_0(\bx) - \mathcal{P}_1(\bx;\by)} } 
  \leq L_0D_h(\bx, \by)\,, \text{ for all } \bx, \by \in \R^N\,,
\end{equation}
with $L_0\geq \sum_{i=1}^M(3\vnorm[F]{{\bf A}_i}^2 + \vnorm[F]{{\bf A}_i}\abs{b_i})$. In this setting,  Model BPG subproblems are solved using Primal-Dual Hybrid Gradient Algorithm (PDHG) \cite{PC2011}.
\medskip

\begin{figure}[t]
  \centering
  \begin{subfigure}[b]{0.24\textwidth}
    \includegraphics[width=\textwidth]{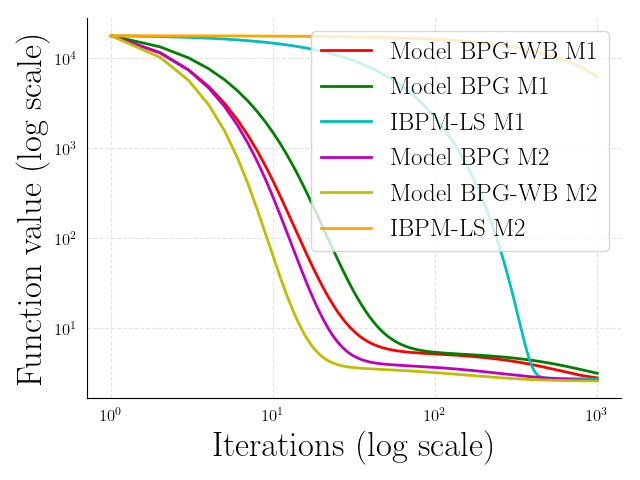}
    \caption{L1 reg}
    % \label{fig:func_plot_1}
  \end{subfigure}
  \begin{subfigure}[b]{0.24\textwidth}
    \includegraphics[width=\textwidth]{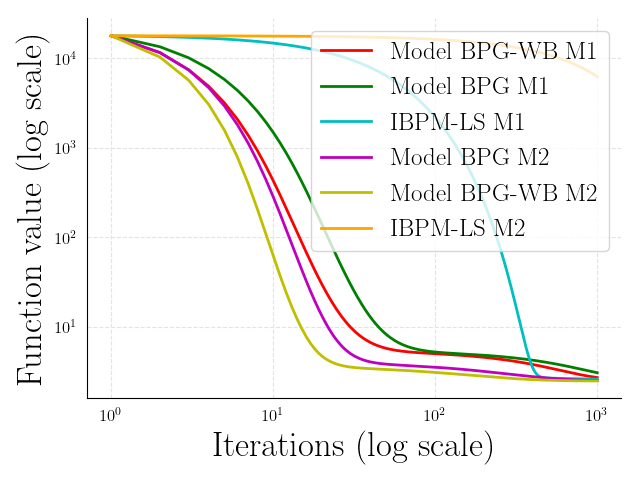}
    \caption{Squared L2 reg}
    % \label{fig:func_plot_2}
  \end{subfigure}
  \begin{subfigure}[b]{0.24\textwidth}
    \includegraphics[width=\textwidth]{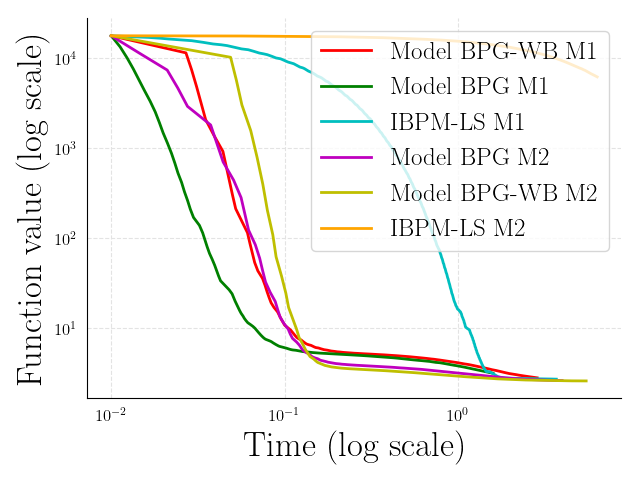}
    \caption{L1 reg}
    % \label{fig:func_plot_1_1}
  \end{subfigure}
  \begin{subfigure}[b]{0.24\textwidth}
    \includegraphics[width=\textwidth]{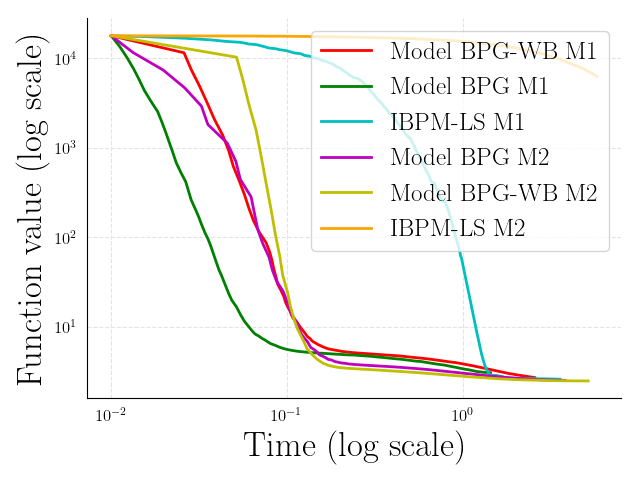}
    \caption{Squared L2 reg}
    % \label{fig:func_plot_2_1}
  \end{subfigure}
  \vspace{1em}
   \caption{In this experiment we compare the  performance of Model BPG, Model BPG with Backtracking (denoted as Model  BPG-WB), and IBPM-LS \cite{ODBP13} on standard phase retrieval problems, with both L1 and squared L2 regularization. For this purpose, we consider M1 model function as in \eqref{eq:p_0_1} without absolute sign (which is the same setting as \cite{BSTV18}), and with M2 model function as in \eqref{eq:p_0_1-2}. Model BPG with M2 \eqref{eq:p_0_1-2} is  faster in both the settings and Model BPG variants perform significantly better than IBPM-LS. By \textit{reg}, we mean \textit{regularization}. } 
   \label{fig:standard_phase_retrieval}
\end{figure}

\begin{figure}[t]
  \centering
  \begin{subfigure}[b]{0.24\textwidth}
    \includegraphics[width=\textwidth]{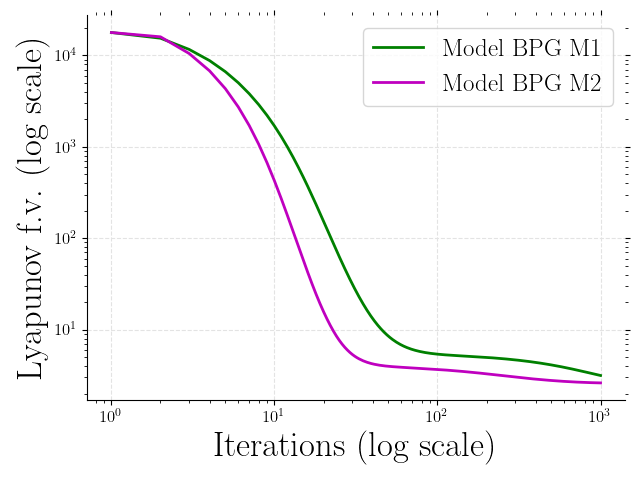}
    \caption{L1 reg}
  \end{subfigure}
  \begin{subfigure}[b]{0.24\textwidth}
    \includegraphics[width=\textwidth]{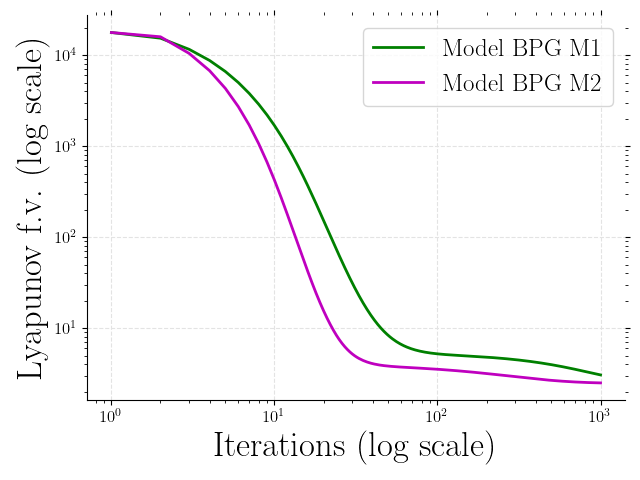}
    \caption{Squared L2 reg}
  \end{subfigure}
  \begin{subfigure}[b]{0.24\textwidth}
    \includegraphics[width=\textwidth]{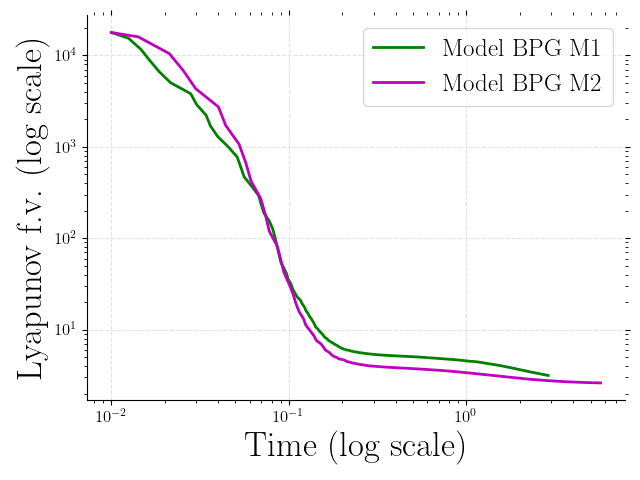}
    \caption{L1 reg}
  \end{subfigure}
  \begin{subfigure}[b]{0.24\textwidth}
    \includegraphics[width=\textwidth]{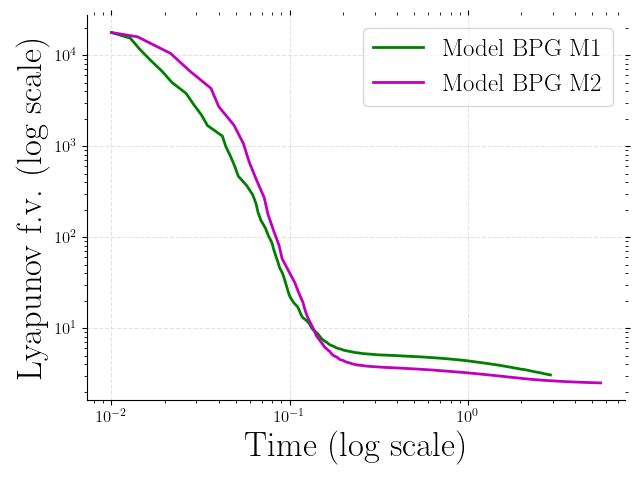}
    \caption{Squared L2 reg}
  \end{subfigure}
  \vspace{1em}
   \caption{ We illustrate that when Model BPG applied to standard phase retrieval problem in \eqref{eq:phase-retrieval}, with model function chosen to be  either Model 1 in \eqref{eq:p_0_1} or  Model 2 in \eqref{eq:p_0_1-2}, result in sequences where the Lyapunov function value evaluations are monotonically nonincreasing. In terms of iterations,  Model BPG with Model 2 (Model BPG M2) is better than Model BPG with Model 1 (Model BPG M1). In terms of time, Model BPG M1 and Model BPG M2 perform almost the same, however, towards the end Model BPG M2 is faster in both the cases. By \textit{reg} we mean \textit{regularization}, and by \textit{Lyapunov f.v.} we mean Lyapunov function values.} 
   \label{fig:standard_phase_retrieval_lyapunov}
\end{figure}

We provide empirical results in Figure~\ref{fig:standard_phase_retrieval}, where we show superior performance of Model BPG variants compared to IBPM-LS, in particular, with the model function provided in \eqref{eq:p_0_1-2}. For simplicity, we choose a constant step-size $\tau$ in all the iterations, such that $\tau \in (0,{1}/{L_0})$. We empirically validate Proposition~\ref{prop:descent-property} in Figure~\ref{fig:standard_phase_retrieval_lyapunov}. All the assumptions required to deduce the global convergence of Model BPG are straightforward to verify, and we leave it as an exercise to the reader. Note that here $\sint\dom\leg = \R^N$, thus the condition $\omega^{\sint\dom\leg}(\bx\iter0) =\omega(\bx\iter0)$ holds trivially.

\subsection{Robust phase retrieval} 
\begin{figure}[t]
  \centering
  \begin{subfigure}[b]{0.24\textwidth}
    \includegraphics[width=\textwidth]{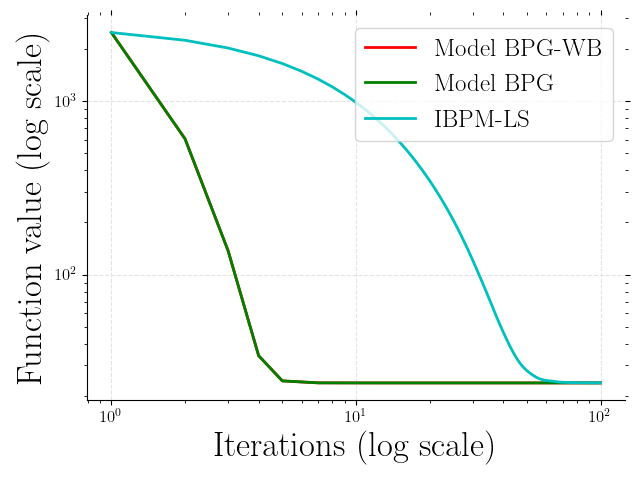}
    \caption{L1 reg}
  \end{subfigure}
  \begin{subfigure}[b]{0.24\textwidth}
    \includegraphics[width=\textwidth]{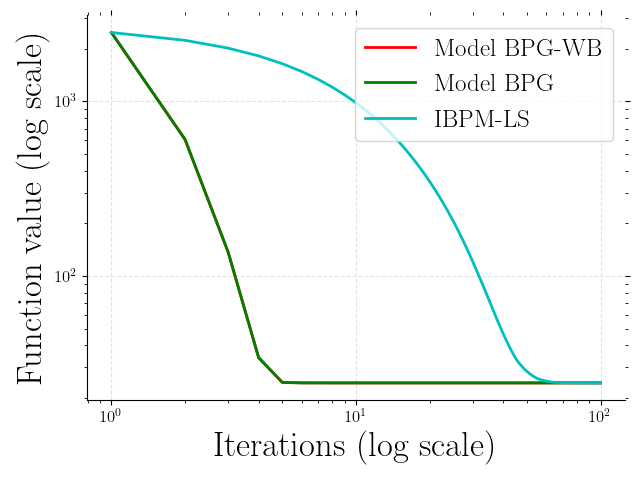}
    \caption{Squared L2 reg}
  \end{subfigure}
  \begin{subfigure}[b]{0.24\textwidth}
    \includegraphics[width=\textwidth]{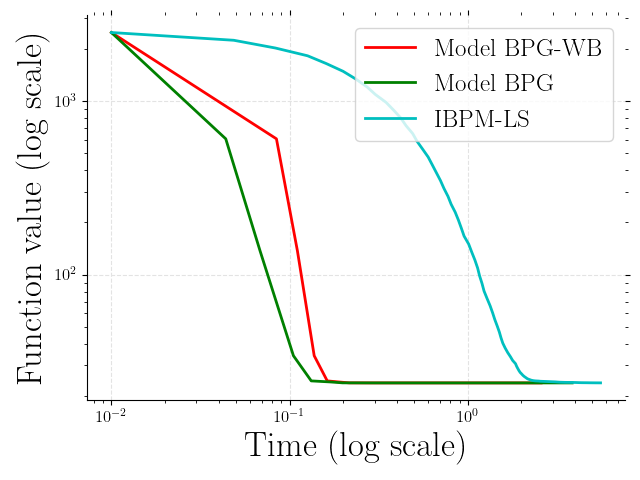}
    \caption{L1 reg}
  \end{subfigure}
  \begin{subfigure}[b]{0.24\textwidth}
    \includegraphics[width=\textwidth]{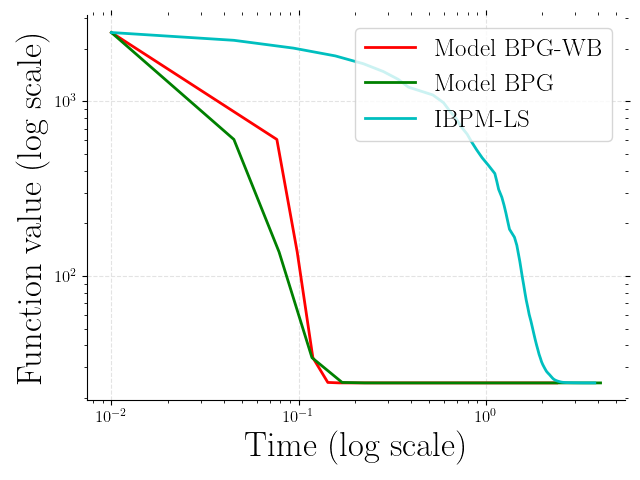}
    \caption{Squared L2 reg}
  \end{subfigure}
  \vspace{1em}
  \caption{In this experiment we consider the  performance of Model BPG vs Model BPG with Backtracking (denoted as Model  BPG-WB) vs IBPM-LS \cite{ODBP13} on robust phase retrieval problems, with both L1 and squared L2 regularization. Model BPG variants perform similarly and are better than IBPM-LS. By \textit{reg}, we mean \textit{regularization}.}
  \label{fig:robust_phase_retrieval}
\end{figure}

\begin{figure}[t]
  \centering
  \begin{subfigure}[b]{0.24\textwidth}
    \includegraphics[width=\textwidth]{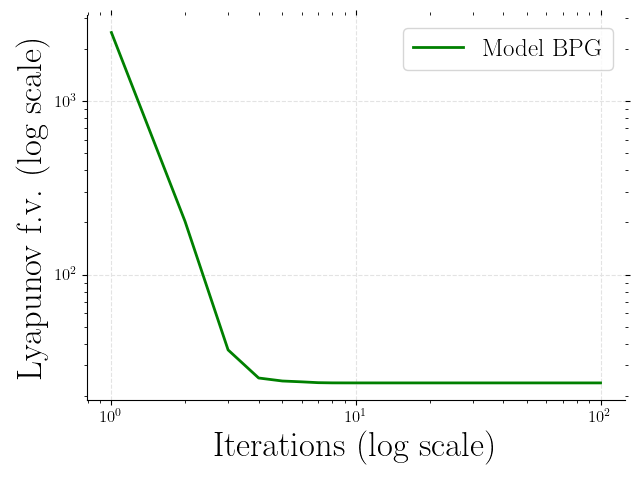}
    \caption{L1 reg}
  \end{subfigure}
  \begin{subfigure}[b]{0.24\textwidth}
    \includegraphics[width=\textwidth]{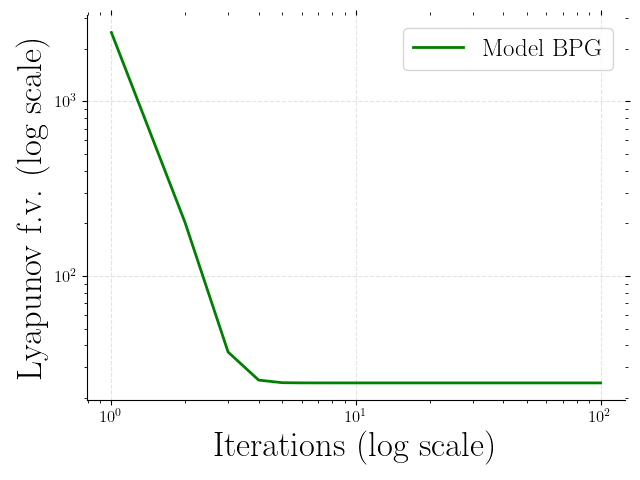}
    \caption{Squared L2 reg}
  \end{subfigure}
  \begin{subfigure}[b]{0.24\textwidth}
    \includegraphics[width=\textwidth]{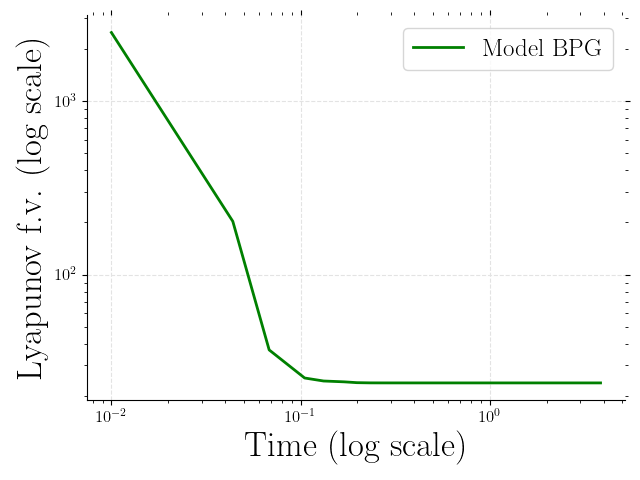}
    \caption{L1 reg}
  \end{subfigure}
  \begin{subfigure}[b]{0.24\textwidth}
    \includegraphics[width=\textwidth]{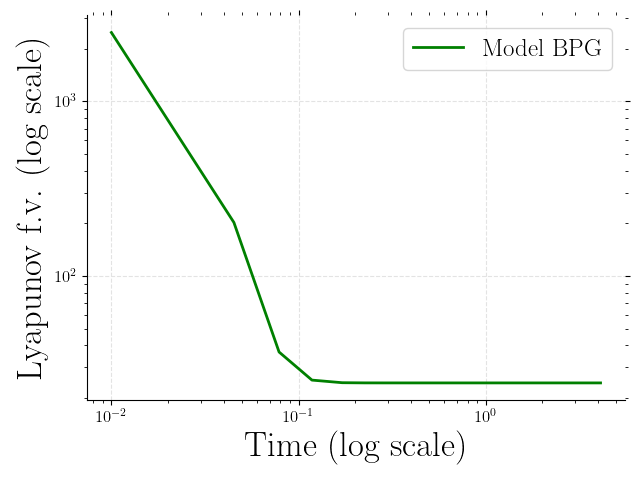}
    \caption{Squared L2 reg}
  \end{subfigure}
  \vspace{1em}
   \caption{ Under the same setting as in Figure~\ref{fig:robust_phase_retrieval}, we  illustrate that Model BPG when applied on robust phase retrieval problems, with both L1 and squared L2 regularization, results in sequences with monotonically decreasing Lyapunov function evaluations, thus validating Proposition~\ref{prop:descent-property}. By \textit{reg} we mean \textit{regularization}, and by \textit{Lyapunov f.v.} we mean Lyapunov function values.}
   \label{fig:robust_phase_retrieval_lyapunov}
\end{figure}

Now, we consider the robust phase retrieval problem, where the goal is the same as standard phase retrieval problem, that is to solve the system of quadratic equations in \eqref{eq:robust-phase-retrieval}.
It is well known that L1 loss is more robust to noise compared to squared L2 loss \cite{FHT2001}. The problem in  \eqref{eq:phase-retrieval} uses squared L2 loss. Here, we consider L1 loss based robust phase retrieval problem, which involves solving the following optimization problem :
\begin{equation*}
  \min_{\bx \in \R^N} f(\bx),\quad f(\bx) := \frac{1}{M}\sum_{i=1}^M\abs{\bx^T{\bf A}_i\bx - b_i} + \mathcal{R}(\bx)\,,
\end{equation*}
where we set $\mathcal{R}(\bx) =\lambda \vnorm[1]{\bx}$ (L1 regularization) or $\mathcal{R}(\bx) = \frac{\lambda}{2}\vnorm[]{\bx}^2$ (squared L2 regularization), for some $\lambda >0$. Such an objective is preferred if the data obtained is noisy, and we require the solution that is robust to noise. We use the Prox-Linear setting described in Section~\ref{ssec:general-composite-problems}, where for any $\bx \in \R^N$ we set the following: 
\[
  f_0(\bx) := \mathcal{R}(\bx)\,,
\]
\[
(F(\bx))_i =   \bx^T{\bf A}_i\bx - b_i\,, \text{ for all } i= 1,\ldots, M\,,
\]
and for any $\tilde{\by} \in \R^M$ we set
\[
  g(\tilde{\by}) := \frac{1}{M}\vnorm[1]{\tilde{\by}}  \,, \text{ for } \tilde{\by} \in \R^M\,.
\]
We consider the following model function. For fixed $\by \in \R^N$, the model function $f(\bx;\by)$ at $\bx \in \R^N$ is given by
\begin{equation}
    f(\bx;\by) := \frac{1}{M}\sum_{i=1}^M\abs{\by^T{\bf A}_i\by  - b_i + \scal{2{\bf A}_i\by}{\bx-\by}}  + \mathcal{R}(\bx)\,.\label{eq:p_0_2-a}
\end{equation}
With the Legendre function $h(\bx) = \frac12 \vnorm[]{\bx}^2$ and as a consequence of triangle property,  a simple calculation reveals that for all $\bx, \by \in \R^N$ we have 
\begin{align*}
  \abs{{f(\bx) - f(\bx;\by)} } \leq 0.5L_1\vnorm[]{\bx-\by}^2\,,
\end{align*}
with  $L_1 \geq \frac{2\sum_{i=1}^M\lambda_{\max}({\bf A}_i)}{M}$. We use a constant step-size $\tau\iter\k = \tau$ such that  $\tau \in (0,{1}/{L_1})$. All the other assumptions of Model BPG are straightforward to verify and we leave it as an exercise to the reader. In each iteration of Model BPG, subproblems take the following form:
\begin{equation*}
  \Argmin_{\bx \in \R^N}\left\{ \frac{1}{M}\sum_{i=1}^M\abs{\by^T{\bf A}_i\by  - b_i + \scal{2{\bf A}_i\by}{\bx-\by}}  + \mathcal{R}(\bx) + \frac{1}{2\tau}\vnorm[]{\bx - \by}^2 \right\}\,,
\end{equation*}
which we solve using Primal-Dual Hybrid Gradient 
Algorithm (PDHG) \cite{PC2011}.  The empirical results are reported in Figure~\ref{fig:robust_phase_retrieval}, where we illustrate the better performance of Model BPG based methods compared to IBPM-LS \cite{ODBP13} on robust phase retrieval problems. We empirically validate Proposition~\ref{prop:descent-property} in Figure~\ref{fig:robust_phase_retrieval_lyapunov}. Note that here $\sint\dom\leg = \R^N$, thus the condition $\omega^{\sint\dom\leg}(\bx\iter0) =\omega(\bx\iter0)$ holds trivially.

\subsection{Poisson linear inverse problems}\label{ssec:poisson-linear-inverse-problems}
We now consider a broad class of problems with varied practical applications, known as Poisson inverse problems \cite{BBDV09, BBT16, OFB19, N05}. The problem setting is as follows. For all $i=1,\ldots,M$, let $b_i >0$, $\ba_i \neq 0$ and $\ba_i \in \R_{+}^N$ be known. Moreover, we  have for any $\bx \in \R_{+}^N$, $\scal{\ba_i}{\bx} >0$ and $\sum_{i=1}^M(\ba_{i})_j > 0$, for all $j=1,\ldots,N$, $i=1,\ldots,M$.  Equipped with these notions, one can write the optimization problem of Poisson linear inverse problems as following:
\begin{equation}\label{eq:main_objective}
\min_{\bx \in \R_{+}}  \left\{ f(\bx) := \sum_{i=1}^M \left( \scal{\ba_i}{\bx} - b_i \log(\scal{\ba_i}{\bx})   \right) + \reg(\bx) \right\}\,,
\end{equation}
where $\reg$ is the regularizing function, which is potentially nonconvex. For simplicity, we set $\reg = 0$.   The function $f_1 : \R^N \to \eR$ at any $\bx \in \R^N$ is defined as following:
\[
  f_1(\bx) := \sum_{i=1}^M \left( \scal{\ba_i}{\bx} - b_i \log(\scal{\ba_i}{\bx})   \right)\,.
\]
Note that the function $f_1$ is coercive. Since $f_1$ is a continuous function, its level set restricted to $\R_+$, i.e., $C := \{\bx \geq 0 : f_1(\bx) \leq f_1(\bx_0)\}$ is compact, for any $\bx_0 \in \R_+$. In order to apply Model BPG, we need $h$ such that the MAP property is satisfied. We consider the Legendre function $h: \R_{++}^N \to \R$ that is given by 
\begin{equation}\label{eq:kernel-dist}
h(\bx) = - \sum_{i=1}^N\log(\bx_i)\,, \quad \text{ for all } \bx \in \R^N_{++},
\end{equation}
where $\bx_i$ is the $i^{\text{th}}$ coordinate of $\bx$. The above given function $h$ is also known as Burg's entropy. Consider the following lemma.
\begin{LEM}\label{lem:map-property}
  Let $h$ be defined as in \eqref{eq:kernel-dist}. For  $L \geq \sum_{i=1}^M b_i$,  the function $Lh - f_1$ and $Lh + f_1$ is convex on $\R^N_{++}$, or equivalently the following $L$-smad property or the MAP property holds true:
  \begin{equation}\label{eq:map-like-property}
  -LD_h(\bx, \bbx) \leq f_1(\bx) - f_1(\bbx) - \scal{\nabla f_1(\bbx)}{\bx - \bbx} \leq LD_h(\bx, \bbx)\,, \text{ for all } \bx, \bbx \in \R^N_{++}\,.
  \end{equation}
\end{LEM}
\begin{proof}
  The proof of convexity of $Lh - f_1$ follows from \cite[Lemma 7]{BBT16}. The function $Lh + f_1$ is convex as $f_1$ is convex.
\end{proof}

When Model BPG is applied to solve \eqref{eq:main_objective} with $h$ given in \eqref{eq:kernel-dist}, if the limit points of the sequence generated by Model BPG lie in $\sint \dom \leg$, our global convergence result is valid. However, it is difficult to guarantee such a condition.  This is because, there can exist subsequences for which certain components of the iterates can tend to zero. In such a scenario, some components of $\nabla^2h(\bx\iter\k)$ will tend to $\infty$, which will lead to the failure of the relative error condition in Lemma~\ref{lem:sub-diff-vanish}. In that case, our analysis cannot guarantee the global convergence of the sequence generated by Model BPG.   
\medskip

Thus, in such a scenario it is important to guarantee that the iterates of Model BPG lie in $\R_{++}^N$. To this regard, we modify the problem \eqref{eq:main_objective}, by adding certain constraint set, such that all the limit points lie in $\sint\dom\leg$.  Then, the global convergence of the sequence generated by Model BPG sequence can be guaranteed. The full objective after the modification is provided below
\begin{equation}\label{eq:main_objective-1}
  \min_{\bx \in \R^N}  \left\{ f(\bx) := \delta_{C_{\eps}}(\bx)  + \sum_{i=1}^M \left( \scal{\ba_i}{\bx} - b_i \log(\scal{\ba_i}{\bx})   \right) + \reg(\bx) \right\}\,,
  \end{equation}
where for certain $\eps >0$ we denote
\[
C_{\eps}  = \{\bx : \bx_i\geq \eps,\, \forall i=1,\ldots,N\}\,,
\]
and $\delta_{C_{\eps}}(\cdot)$ is the indicator function of the set $C_{\eps}$. We consider $\reg = 0$ or $\reg(\bx) = \lambda \vnorm[1]{\bx}$ or $\reg(\bx) = \lambda \frac{\norm{\bx}^2}{2}$, with certain $\lambda >0$. Note that $C_{\eps} \subset \R_{+}$. For practical purposes, $C_{\eps}$ is almost the same as $\R_{+}$, when $\eps$ is chosen sufficiently small. Note that the choice of $\eps$ is only heuristic.  To this end, with $\bbx \in C_{\eps}$, we consider the following model function which, when evaluated at $\bx$ gives:
\begin{equation}\label{eq:model-function}
  f(\bx ; \bbx) := \delta_{C_{\eps}}(\bx)  + f_1(\bbx) + \scal{\nabla f_1(\bbx)}{\bx - \bbx} + \reg(\bx)\,.
\end{equation}
% Well definedness of the Model BPG update step is guaranteed by  Lemma~\ref{lem:well-posedness-2}. 
% \medskip

%  Using the MAP property, with certain $\ub > 0$ for all $\bx, \bbx \in C_{\eps}$ we obtain the following:
% \begin{align*}
%   f(\bx) \leq  f(\bx;\bbx) + {\bar L} D_h(\bx, \bbx) \leq  f(\bx;\bbx) + \frac{1}{\lambda} D_h(\bx, \bbx)\,,
% \end{align*}
% where in the last step we used the inequality $\frac{1}{\lambda} > {\bar L}$. It is straightforward to deduce that the following condition holds true:
%   \begin{equation*}
%     \lim_{\vnorm{\bx} \to \infty} \lambda\fun(\bx;\bbx) + D_h(\bx, \bbx) = \infty\,,
%   \end{equation*}
%   as $\lim_{\vnorm{\bx} \to \infty} f(\bx) = \infty$ for all the considered choices of $\reg$. Well definedness of the Model BPG update step can then be gauranteed by the coercivity of $f(\cdot; \bbx) + D_h(\cdot,\bbx)$ and strict convexity of $h$.
% \medskip

The Legendre function in \eqref{eq:kernel-dist} is still valid as $C_{\eps} \subset \R_{+}$, and the MAP property holds true as a consequence of Lemma~\ref{lem:map-property}. The coercivity of the function $f$ along with Proposition~\ref{prop:function-descent-property} implies that the iterates of Model BPG will lie in the compact convex set $\{\bx : f(\bx) \leq f(\bx_0)\}$. Thus, the sequence generated by Model BPG is bounded. The analysis falls under the category of additive composite problems given in Section~\ref{ssec:forward-backward}, where we set $f_1 := f_1$  and $f_0(\cdot):= \delta_{C_{\eps}}(\cdot) + \reg(\cdot)$. In the earlier discussion, we have proved the crucial assumptions for applying Model BPG to Poisson linear inverse problems. The rest of the assumptions in Theorem~\ref{thm:full-conv-additive} are straightforward to verify and we leave it as an exercise to the reader. We now provide closed form expressions for the update step \eqref{eq:alg-BregMin-bt:update} in three settings of $\reg$.

  \paragraph{Closed form update step - No regularization.} Set $\reg = 0$. The update step of Model BPG involves solving the following subproblem: 
  \[
  \bx\iter\kp \in \argmin_{\bx} \delta_{C_{\eps}}(\bx) + f(\bx\iter\k) + \scal{\nabla f(\bx\iter\k)}{\bx-\bx\iter\k} + \frac{1}{\tau_k}D_h(\bx, \bx\iter\k)\,.
  \]
  The optimality condition for the $i^{\text{th}}$ component of $\bx\iter\kp$ due to Fermat's rule is given by
  \[
  0 =  (\bv\iter\kp)_i +  \nabla f(\bx\iter\k)_{i} + \frac{1}{\tau_k} \Big( \frac{1}{(\bx\iter\k)_i} - \frac{1}{(\bx\iter\kp)_i }\Big)\,,
  \]
  for some  $\bv\iter\kp  \in N_{C_{\eps}}(\bx\iter\kp)$. Thus, we deduce that with $\tau\iter\k$ chosen such that $1 + \tau\iter\k \nabla f(\bx\iter\k)_i(\bx\iter\k)_i >0$, for $i=1,\ldots, N$, the solution is given by
  \begin{equation}\label{eq:approx-burg}
    \bx\iter\kp = \max\left\{{\eps}, \frac{\bx\iter\k}{1 + \tau\iter\k \nabla f(\bx\iter\k)\bx\iter\k}\right\}\,,
  \end{equation}
  where all the operations are performed element-wise.  
  \paragraph{Closed form update step - L1 regularization.} We consider here the standard L1 regularization setting, where with certain $\lambda >0$ we set $\reg(\bx) = \lambda \vnorm[1]{\bx}$.  The update step of Model BPG involves solving the following subproblem: 
  \[
  \bx\iter\kp \in \argmin_{\bx} \delta_{C_{\eps}}(\bx) +\lambda \vnorm[1]{\bx} + f(\bx\iter\k) + \scal{\nabla f(\bx\iter\k)}{\bx-\bx\iter\k} + \frac{1}{\tau_k}D_h(\bx, \bx\iter\k)\,.
  \]
  Based on  \cite[Section 5.2]{BBT16} and Fermat's rule we deduce that with $\tau\iter\k$ chosen such that $1 + \tau\iter\k\lambda(\bx\iter\k)_i  +\tau\iter\k \nabla f(\bx\iter\k)_i(\bx\iter\k)_i >0$, for $i=1,\ldots, N$, the closed form solution is given by 
  \begin{equation}
    \bx\iter\kp = \max\left\{{\eps}, \frac{\bx\iter\k}{1 + \tau\iter\k\lambda\bx\iter\k +\tau\iter\k \nabla f(\bx\iter\k)\bx\iter\k}\right\}\,,
  \end{equation}
  where all the operations are performed element-wise.  
\begin{figure}[t]
  \centering
  \begin{subfigure}[b]{0.32\textwidth}
    \includegraphics[width=\textwidth]{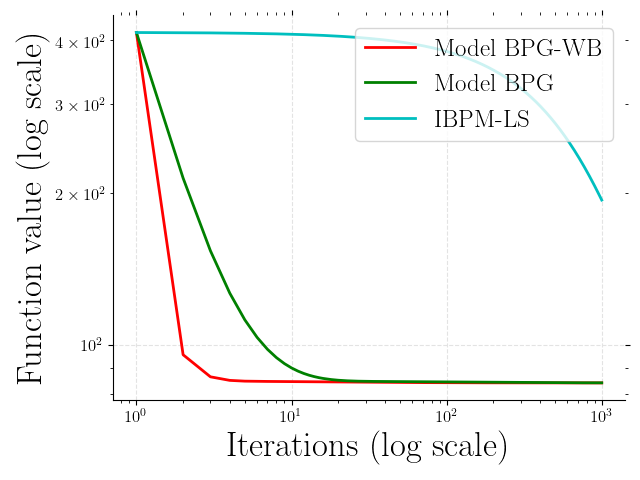}
    \caption{L1 regularization}
  \end{subfigure}
  \begin{subfigure}[b]{0.32\textwidth}
    \includegraphics[width=\textwidth]{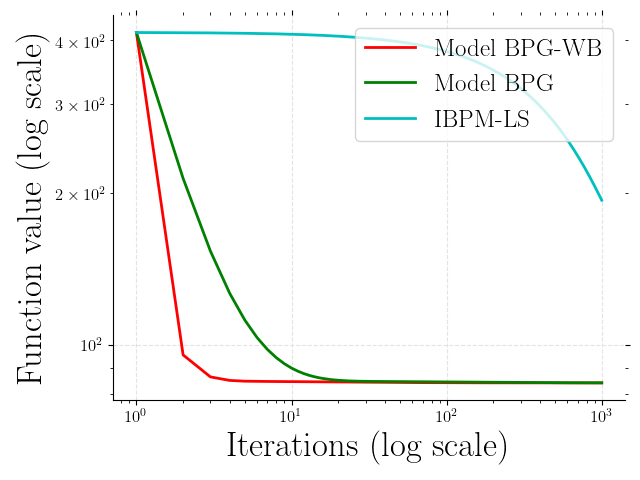}
    \caption{Squared L2 regularization}
  \end{subfigure}
  \begin{subfigure}[b]{0.32\textwidth}
    \includegraphics[width=\textwidth]{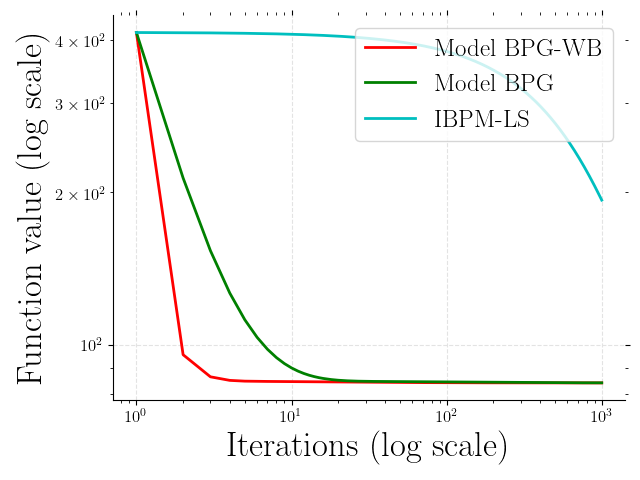}
    \caption{No regularization}
  \end{subfigure}
  \vspace{1em}
  \caption{In this experiment we compare the  performance of Model BPG, Model BPG with Backtracking (denoted as Model  BPG-WB) and IBPM-LS \cite{ODBP13} on Poisson linear inverse problems with L1 regularization, squared L2 regularization and with no regularization. We set the regularization parameter $\lambda$ to $0.1$. The plots illustrate that Model BPG-WB is faster in all the settings, followed by Model BPG. }
  \label{fig:poisson_linear_inverse_problems_function}
\end{figure}
\begin{figure}[t]
  \centering
  \begin{subfigure}[b]{0.32\textwidth}
    \includegraphics[width=\textwidth]{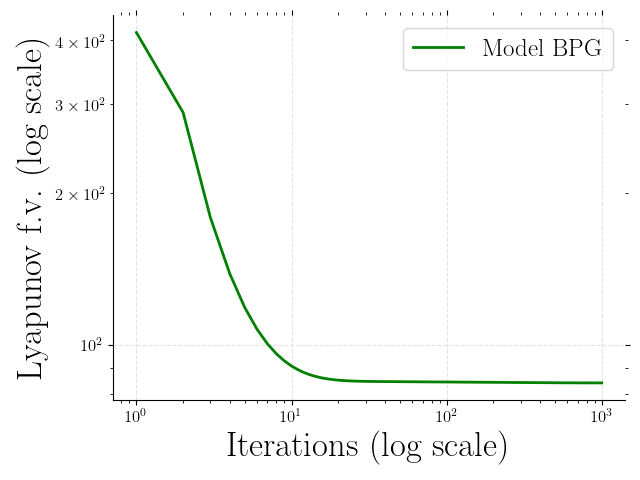}
    \caption{L1 regularization}
  \end{subfigure}
  \begin{subfigure}[b]{0.32\textwidth}
    \includegraphics[width=\textwidth]{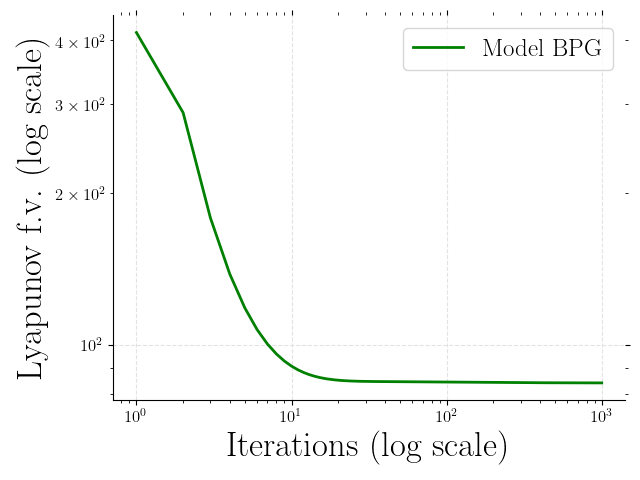}
    \caption{Squared L2 regularization}
  \end{subfigure}
  \begin{subfigure}[b]{0.32\textwidth}
    \includegraphics[width=\textwidth]{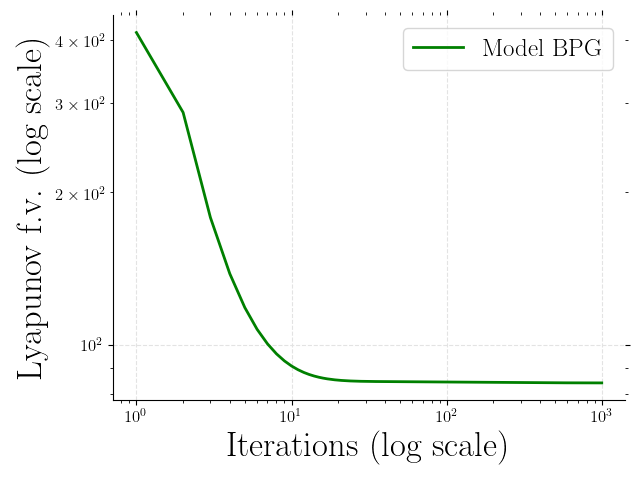}
    \caption{No regularization}
  \end{subfigure}
  \vspace{1em}
\caption{By \textit{Lyapunov f.v.} we mean Lyapunov function values. Under the same setting as in Figure~\ref{fig:poisson_linear_inverse_problems_function}, we illustrate here that Model BPG results in sequences that have monotonically nonincreasing Lyapunov function value evaluations.}
\label{fig:poisson_linear_inverse_problems_lyapunov}
\end{figure}

\paragraph{Closed form update step -  L2 regularization.}  We consider here the standard L2 regularization setting, where with certain $\lambda >0$ we set $\reg(\bx) = \frac{\lambda}{2}\vnorm[2]{\bx}^2$.  The update step of Model BPG involves solving the following subproblem: 
\[
\bx\iter\kp \in \argmin_{\bx} \delta_{C_{\eps}}(\bx) +\frac{\lambda}{2}\vnorm[2]{\bx}^2 + f(\bx\iter\k) + \scal{\nabla f(\bx\iter\k)}{\bx-\bx\iter\k} + \frac{1}{\tau_k}D_h(\bx, \bx\iter\k)\,.
\]
The optimality condition for the $i^{\text{th}}$ component of $\bx\iter\kp$ due to Fermat's rule is given by
  \[
  0 =  (\bv\iter\kp)_i + \lambda (\bx\iter\kp)_i + \nabla f(\bx\iter\k)_{i} + \frac{1}{\tau_k} \Big( \frac{1}{(\bx\iter\k)_i} - \frac{1}{(\bx\iter\kp)_i }\Big)\,,
  \]
  for some  $\bv\iter\kp  \in N_{C_{\eps}}(\bx\iter\kp)$.  Based on  \cite[Section 5.2]{BBT16} we deduce that with $\tau\iter\k$ chosen such that $1  +\tau\iter\k \nabla f(\bx\iter\k)_i(\bx\iter\k)_i  + \tau\iter\k \lambda \eps  >0$, for $i=1,\ldots, N$, the closed form solution is given by 
\begin{equation}
  \bx\iter\kp = \max\left\{{\eps}, \frac{\sqrt{ (1+ \tau\iter\k \bx\iter\k \nabla f(\bx\iter\k))^2 + 4\lambda\tau\iter\k \bx\iter\k^2} - (1+ \tau\iter\k \bx\iter\k \nabla f(\bx\iter\k))}{2\lambda \tau\iter\k \bx\iter\k} \right\}\,,
\end{equation}
where all the operations are performed element-wise.  
\medskip

  The empirical results are reported in Figure~\ref{fig:poisson_linear_inverse_problems_function}, where we illustrate the better performance of Model BPG based methods compared to IBPM-LS \cite{ODBP13}, when applied on Poisson linear inverse problems.  We empirically validate Proposition~\ref{prop:descent-property} in Figure~\ref{fig:poisson_linear_inverse_problems_lyapunov}. Note that here $\sint\dom\leg = \R^N_{++}$. Based on the aforementioned closed form solutions it is clear that the sequence generated by Model BPG lies in $C_{\eps}$. The condition $C_{\eps} \subset \sint\dom\leg$ implies that $\omega^{\sint\dom\leg}(\bx\iter0) =\omega(\bx\iter0)$ holds true.

\section{Conclusion}
Bregman proximal minimization framework is prominent in solving additive composite problems, in particular, using BPG \cite{BSTV18} algorithm or its variants \cite{MOPS2020}. However, extensions to generic composite problems was an open problem. To this regard, based on foundations of \cite{drusvyatskiy2019nonsmooth, OFB19}, we proposed Model BPG algorithm that is applicable to a vast class of nonconvex nonsmooth problems, including generic composite problems. Model BPG relies on certain function approximation, known as model function, which preserves first order information about the function. The model error is bounded via certain Bregman distance, which drives the  global convergence analysis of the sequence generated by Model BPG. The analysis is nontrivial and requires significant changes compared to the standard analysis of \cite{BSTV18,BST14, AB09, ABS13}.  Moreover, we numerically illustrate the superior performance of Model BPG on various real world applications.

\section{Acknowledgments}
Mahesh Chandra Mukkamala and Peter Ochs thank German Research Foundation for providing financial support through DFG Grant OC 150/1-1. 

{
\small

{
\bibliographystyle{plain}
\bibliography{ochs}
}
}
\newpage
\appendix
\section{Additional preliminaries}\label{sec:additional-preliminaries}
We work with extended-valued functions $\map{f}{\X}{\eR}$, $\eR:= \R\cup\set{+\infty}$. The \emphdef{domain} of $f$ is $\dom f:= \set{x\in \X\,\vert\, f({\bf x}) < +\infty}$ and a function $f$ is \emphdef{proper}, if $\dom f\neq\emptyset$. It is \emphdef{lower semi-continuous} (or \emphdef{closed}), if $\liminf_{{\bf x}\to {\bf \bar {\bf x}}} f({\bf x}) \geq f(\bar {\bf x})$ for any $\bar {\bf x}\in \X$. Let $\sint \Omega$ denote the \emphdef{interior} of $\Omega\subset\X$. We use the notation of \emphdef{$f$-attentive convergence} $x \fto f \bar {\bf x} \Leftrightarrow ({\bf x},f({\bf x})) \to (\bar {\bf x}, f(\bar {\bf x}))$ and the notation $\k\fto K\infty$ for some $K\subset\N$ to represent $\k\to\infty$ where $\k\in K$.   The \emphdef{indicator function} $\ind C$ of a set $C\subset\X$ is defined by $\ind C({\bf x})=0$, if ${\bf x}\in C$ and $\ind C({\bf x})=+\infty$, otherwise.  The \emphdef{(orthogonal) projection} of $\bar {\bf x}$ onto $C$, denoted $\proj_C(\bar {\bf x})$, is given by a minimizer of $\min_{{\bf x}\in C}\, \vnorm{{\bf x}-\bar {\bf x}}$. A \emphdef{set-valued mapping} $\smap{T}{\X}{\Y}$ is defined by its \emphdef{graph} $\Graph T:=\set{({\bf x},{\bf v})\in \X\times\Y \setsep {\bf v}\in T({\bf x})}$ with domain given by $\dom T:=\set{{\bf x}\in\X\setsep T({\bf x})\neq\emptyset}$. Following \cite[Def. 8.3]{Rock98}, we introduce subdifferential notions for nonsmooth functions.   The \emphdef{Fr\'echet subdifferential} of $f$ at $\bar {\bf x} \in\dom f$ is the set $\rpartial f(\bar {\bf x})$ of elements ${\bf v} \in \X$ such that
\[
  \liminf_{\substack{{\bf x}\to \bar {\bf x}\\ {\bf x}\neq \bar {\bf x}}} \frac{f({\bf x}) - f(\bar {\bf x}) - \scal{{\bf v}}{{\bf x}-\bar {\bf x}}}{\vnorm{{\bf x}-\bar {\bf x}}} \geq 0  \,.
\]
For $\bar {\bf x}\not\in \dom f$, we set $\rpartial f(\bar {\bf x}) = \emptyset$. 
The \emphdef{(limiting) subdifferential} of $f$ at $\bar {\bf x}\in\dom f$ is defined by
\[
  \partial f(\bar {\bf x}) := \set{{\bf v}\in \X \setsep \exists\, {\bf y}_k \fto{f} \bar {\bf x},\;{\bf v}\iter\k\in \rpartial f({\bf y}_k),\;{\bf v}\iter\k \to {\bf v}} \,,
\]
and $\partial f(\bar {\bf x}) = \emptyset$ for $\bar {\bf x} \not\in \dom f$.  As a direct consequence of the definition of the limiting subdifferential, we have the following closedness property at any $\bar {\bf x}\in \dom f$: 
\begin{equation} \label{eq:closedness-lim-subdiff}
  {\bf y}_k \fto{f} \bar {\bf x},\ {\bf v}\iter\k\to {\bar {\bf v}},\ \text{and for all } \k\in\N\colon {\bf v}\iter\k \in \partial f({\bf y}_k)\quad \Longrightarrow\quad  {\bar {\bf v}}\in \partial f(\bar {\bf x}) \,.
\end{equation}
A point $\bar {\bf x}\in \dom f$ for which ${\bf 0}\in \partial f(\bar {\bf x})$ is a called a \emphdef{critical point}, which is a necessary optimality condition (Fermat's rule \cite[Thm. 10.1]{Rock98}) for $\bar {\bf x}$ being a local minimizer. We define the \emphdef{set of (global) minimizers} of a function $f$ by 
\[
  \Argmin_{x\in \X}\, f({\bf x}) := \set{x\in \X\,\vert\, f({\bf x}) = \inf_{\bar {\bf x}\in\X} f(\bar {\bf x}) }\,,
\]
and the \emphdef{(unique) minimizer} of $f$ by $\argmin_{x\in\X}\, f({\bf x})$ if $\Argmin_{x\in \X}\, f({\bf x})$ consists of a single element. As shorthand, we also use $\Argmin f$ and $\argmin f$.
\medskip

% \begin{LEM}[subgradients versus gradients {\cite[Ex. 8.8]{Rock98}}] \label{lem:subdiff-C1}
%     If $f$ is differentiable at $\bar {\bf x}$, then $\rpartial f(\bar {\bf x}) = \set{\nabla f(\bar {\bf x})}$ and $\nabla f(\bar {\bf x})\subset \partial f(\bar {\bf x})$. If $f$ is smooth on a neighborhood of $\bar {\bf x}$, then $\partial f(\bar {\bf x}) = \set{\nabla f(\bar {\bf x})}$. If $f= f_0+f_1$ with $f_0$ finite at $\bar {\bf x}$ and $f_1$ smooth on a neighborhood of $\bar {\bf x}$, then $\partial f(\bar {\bf x}) = \partial f_0(\bar {\bf x}) + \nabla f_1(\bar {\bf x})$.
%   \end{LEM}
In order to conveniently work with Bregman distances, we collect a few properties.
\begin{PROP} \label{prop:basic-Breg-properties}
  Let $\leg\in\clLeg$ and $\bregmap[\leg]$ be the associate Bregman distance.
  \begin{itemize}
    \item[\ii1] $\bregmap[\leg]$ is strictly convex on every convex subset of $\dom\partial \leg$ with respect the first argument.
    \item[\ii2] For ${\bf y}\in \sint\dom \leg$, it holds that $\breg[\leg]{\bf x}{{\bf y}} = 0$ if and only if ${\bf x}={\bf y}$.
    \item[\ii3] For ${\bf x}\in \X$ and ${\bf u},{\bf v}\in \sint\dom\leg$ the following \emphdef{three point identity} holds:
    \begin{equation} \label{eq:breg-id:1}
      \breg[\leg]{{\bf x}}{{\bf u}} = \breg[\leg]{{\bf x}}{{\bf v}} + \breg[\leg]{{\bf v}}{{\bf u}} + \scal{{\bf x}-{\bf v}}{\nabla\leg({\bf v}) - \nabla\leg({\bf u})} \,.
    \end{equation}
  \end{itemize}
\end{PROP}
\begin{proof}
$\ii1$ and $\ii2$ follow directly from the definition of $h$ being essentially strictly convex. $\ii3$ is stated in \cite[Prop. 2.3]{BBC03}. It follows from the definition of a Bregman distance.
\end{proof}
Associated with such a distance function is the following proximal mapping.
\begin{DEF}[{Bregman proximal mapping \cite[Def. 3.16]{BBC03}}] \label{def:Bregman-prox}
  Let $\map {f}{\X}{\eR}$ and $\bregmap[\leg]$ be a Bregman distance associated with $\leg\in\clLeg$. The \emphdef{$\bregmap[\leg]$-prox} (or Bregman proximal mapping) associated with $f$ is defined by
  \begin{equation} \label{eq:def-Breg-Prox}
    \Prox[\leg]_{f} ({\bf y}) := \argmin_{x}\, f({\bf x}) + \breg[\leg]{{\bf x}}{{\bf y}} \,.
  \end{equation}
\end{DEF}
In general, the proximal mapping is set-valued, however for a convex function, the following lemma simplifies the situation.
\begin{LEM} \label{lem:breg-prox-single-valued}
  Let $\map{f}{\X}{\eR}$ be a proper, closed, convex function that is bounded from below, and $\leg\in\clLeg$ such that $\dom f \cap \sint\dom\leg \neq\emptyset$. Then the associated Bregman proximal mapping $\Prox[\leg]_{f}$ is single-valued on its domain and maps to $\dom f \cap \sint\dom\leg$.
\end{LEM}
\begin{proof}
Single-valuedness follows from \cite[Cor.~3.25(i)]{BBC03}. The second claim is from \cite[Prop. 3.23(v)(b)]{BBC03}.
\end{proof}
\begin{PROP} \label{prop:three-point-ineq}
  Let $\map{f}{\X}{\eR}$ be a proper, closed, convex function that is bounded from below, and $\leg\in\clLeg$ such that $\sint\dom\leg \cap \dom f\neq\emptyset$. For ${\bf y}\in \sint\dom\leg$, $\hat{{\bf x}} = \Prox[\leg]_f({\bf y})$, and any ${\bf x}\in \dom f$ the following inequality holds:
  \[
      f({\bf x}) + \breg[\leg]{{\bf x}}{\bar {\bf x}} \geq f(\hat {\bf x}) + \breg[\leg]{\hat {\bf x}}{ \bar {\bf x}} + \breg[\leg]{{\bf x}}{\hat {\bf x}} \,.
  \]
\end{PROP}
\begin{proof}
See \cite[Lem.~3.2]{CT93}.
\end{proof}
For examples and more useful properties of Bregman functions, we refer the reader to \cite{BB97,BBC03,BBT16,Nguyen17}.

\section{Gradient-like Descent Sequence}\label{sec:gradient-like-descent}
We briefly review the concept of gradient-like descent sequence, which is given below. For ease of global convergence analysis of Model BPG we use following results from \cite{Ochs16b}. Let $\map{\F}{\R^\dimN\times \R^\dimP}{\eR}$ be a proper, lower semi-continuous function that is bounded from below, then assume the following assumption from \cite{Ochs16b} holds.
\begin{ASS}[Gradient-like Descent Sequence \cite{Ochs16b}]\label{ass:Hs}
  Let $\seq[\n\in\N]{{\bf u}_\n}$ be a sequence of parameters in $\R^\dimP$ and let
  $\seq[n\in\N]{\eps_n}$ be an $\ell_1$-summable sequence of non-negative real
  numbers. Moreover, we assume there are sequences $\seq[n\in\N]{a_n}$,
  $\seq[n\in\N]{b_n}$, and $\seq[\n\in\N]{d_\n}$ of non-negative real numbers,
  a non-empty finite index set $I\subset\Z$ and $\theta_i\geq 0$, $i\in I$, with 
  $\sum_{i\in I}\theta_i = 1$ such that the following holds:

\begin{enumerate}[label=(H\arabic*),ref=(H\arabic*)]
\item\label{ass:Hs:descent} (Sufficient decrease condition) 
      For each $n\in \N$, it holds that
\[
     \F({\bf x}\iter\np,{\bf u}\iter\np) + a\pit\n d\pit\n^2 \leq  \F({\bf x}\iter\n,{\bf u}\iter\n)\,.
\]
\item\label{ass:Hs:error} (Relative error condition) 
      For each $n\in\N$, the following holds: (set $d\pit{j}=0$ for $j\leq0$)
\[
    b\pit{\np} \vnorm[-]{\partial \F({\bf x}\iter\np,{\bf u}\iter\np) } 
    \leq b \sum_{i\in I} \theta\pit{i}d\pit{\np-i} + \eps\pit{n+1} \,.
\]
\item\label{ass:Hs:cont} (Continuity condition) 
      There exists a subsequence $\seq[j\in\N]{({\bf x}\iter{n_j},{\bf u}\iter{n_j})}$ and 
      $({\tilde {\bf x}},{\tilde {\bf u}})\in\R^\dimN\times\R^\dimP$ such that 
\[
  ({\bf x}\iter{n_j},{\bf u}\iter{n_j}) \Fto[\F] ({\tilde {\bf x}},{\tilde {\bf u}})\quad \text{as}
  \quad j\to\infty \,.
\]
\item\label{ass:Hs:distance} (Distance condition) 
      It holds that
  \begin{align*}
    &d\pit\n \to 0 \Longrightarrow \vnorm[2]{{\bf x}\iter\np-{\bf x}\iter\n} \to 0 
  \qquad\text{and}\\
  &
%  d\pit\n = 0 \Longrightarrow {\bf x}\iter\np={\bf x}\iter\n \,,
  \exists \n^\prime\in\N\colon \forall \n\geq \n^\prime \colon d\pit\n = 0 
  \Longrightarrow 
  \exists \n^{\prime\prime}\in\N\colon \forall \n\geq \n^{\prime\prime} \colon 
      {\bf x}\iter\np={\bf x}\iter\n 
  \end{align*}

\item\label{ass:Hs:params} (Parameter condition) 
      It holds that
 \[
 \seq[\n\in\N]{b\pit\n}\not\in\ell_1\,, \quad  
 \sup_{n\in\N} \frac 1{b\pit\n a\pit\n} < \infty\,, \quad 
 \inf_\n a\pit\n =: \underline a > 0\,.
\]
\end{enumerate}
\end{ASS}
Such an assumption is crucial in order to obtain global convergence of the sequences generated by Model BPG. Assumption~\ref{ass:Hs} is more general compared to the conditions that arise in standard gradient-like sequence \cite{BSTV18}, which is basically based on the first three conditions. 
\medskip

We now provide the global convergence statement from  \cite{Ochs16b}, based on Assumption~\ref{ass:Hs}. Firstly, denote the following. The set of limit points of a bounded sequence $\seq[\n\in\N]{({\bf x}\iter\n,{\bf u}\iter\n)}$ is given by  $\omega({\bf x}\iter0,{\bf u}\iter0) := \limsup_{\n \to \infty}\, \set{({\bf x}\iter\n,{\bf u}\iter\n)}\,,$ and the subset of $\F$-attentive limit points is denoted by 
\[
  \omega_{\F}({\bf x}\iter0,{\bf u}\iter0) := \set{(\bar {\bf x},\bar {\bf u})\in \omega({\bf x}\iter0,{\bf u}\iter0)
  \setsep ({\bf x}\iter{\n_j},{\bf u}\iter{\n_j}) \Fto[\F] (\bar {\bf x},\bar {\bf u}) \text{ for } j\to \infty} \,.
\]

\begin{THM}[Global convergence {\cite[Theorem 10]{Ochs16b}}]\label{thm:KL-theorem-descent}
  Suppose $\F$ is a proper lower semi-continuous \KL function that is bounded 
  from below. Let $\seq[\n\in\N]{{\bf x}\iter\n}$ be a bounded sequence generated by 
  an abstract algorithm parametrized by a bounded sequence 
  $\seq[\n\in\N]{{\bf u}\iter\n}$ that satisfies Assumption~\ref{ass:Hs}. Assume that
  $\F$-attentive convergence holds along converging subsequences of 
  $\seq[\n\in\N]{({\bf x}\iter\n,{\bf u}\iter\n)}$, i.e.
  $\omega({\bf x}\iter0,{\bf u}\iter0)=\omega_{\F}({\bf x}\iter0,{\bf u}\iter0)$. 
  Then, the following holds:
  \begin{enumerate}
    \item\label{thm:KL-theorem-finite-length-abstr-dist} 
          The sequence $\seq[\n\in\N]{d\pit\n}$ satisfies  $\sum_{\k=0}^{\infty} d\pit\k < +\infty\,,$ i.e., the trajectory of the sequence $\seq[\n\in\N]{{\bf x}\iter\n}$ has 
          finite length with respect to the abstract distance measures 
          $\seq[\n\in\N]{d\pit\n}$.
    \item\label{thm:KL-theorem-finite-length-dist} 
          Suppose $d\pit\k$ satisfies 
          $\vnorm[2]{{\bf x}\iter\kp-{\bf x}\iter\k}\leq \bar c d\pit{\k+\k^\prime}$ for some 
          $\k^\prime\in\Z$ and $\bar c\in\R$, then 
          $\sum_{\k=0}^{\infty} \vnorm[2]{{\bf x}\iter\kp-{\bf x}\iter\k} < +\infty\,,$ and the trajectory of the sequence $\seq[\n\in\N]{{\bf x}\iter\n}$ has 
          a finite Euclidean length, and thus $\seq[\n\in\N]{{\bf x}\iter\n}$ 
          converges to $\tilde x$ from \ref{ass:Hs:cont}.
    \item\label{thm:KL-theorem-critical-point} 
    Moreover, if $\seq[\n\in\N]{{\bf u}_\n}$ is a converging sequence, then 
    each limit point of $\seq[n\in\N]{({\bf x}_\n,{\bf u}_\n)}$ is a critical point, which in the situation of \ref{thm:KL-theorem-finite-length-dist} is the unique point $(\tilde {\bf x}, \tilde {\bf u})$ from \ref{ass:Hs:cont}.
  \end{enumerate}
\end{THM}

\section{Proof of Example~\ref{ex:running-example}}\label{sec:proof-running-example}
The model error  is given by
  \begin{align*}
    \abs{f(\bx) - f(\bx;\bbx)} &\leq \abs{g(\bx) - g(\bbx) - \scal{\nabla g(\bbx)}{\bx - \bbx}}\,,\\
    &\leq \abs{\scal{\nabla g(\bbx + s (\bx - \bbx)) - \nabla g(\bbx)}{\bx- \bbx}}\,,\\
    & \leq \vnorm[]{\nabla g(\bbx + s (\bx - \bbx)) - \nabla g(\bbx)}\vnorm[]{\bx- \bbx}\,.
  \end{align*}
  where in the second inequality we use mean value theorem with  $s \in [0,1]$, the third inequality is a simple application of Cauchy-Schwarz rule. On further application of the fundamental theorem of calculus, we have 
  \begin{align*}
    \vnorm[]{\nabla g(\bbx + s (\bx - \bbx)) - \nabla g(\bbx)} &= \vnorm[]{\int_{0}^1\nabla^2 g(\bbx + s(\bx- \bbx))(\bx- \bbx)ds }\,,\\
    &\leq \int_{0}^1\vnorm[]{\nabla^2  g(\bbx + s(\bx- \bbx))}\vnorm[]{\bx- \bbx} ds\,.
  \end{align*}
  Using the fact that $\nabla^2g(\bx) = 4 \vnorm[]{\bx}^2 I + 8 \bx\bx^T$, and $\vnorm[]{\nabla^2g(\bx)} \leq 12\vnorm[]{\bx}^2$ we obtain
  \begin{align*}
    \vnorm[]{\nabla g(\bbx + s (\bx - \bbx)) - \nabla g(\bbx)} &\leq 12\int_{0}^1\vnorm[]{\bbx + s(\bx- \bbx)}^2 \vnorm[]{\bx- \bbx} ds\,,\\
    &\leq 12\int_{0}^1 \left(2\vnorm[]{\bbx}^2 + 2s^2\vnorm[]{(\bx- \bbx)}^2 \right)\vnorm[]{\bx- \bbx} ds\,,\\
    &\leq 24\vnorm[]{\bbx}^2\vnorm[]{\bx- \bbx} + 8\vnorm[]{\bx- \bbx}^3\,,
  \end{align*}
  where in the second step we used the inequality $\vnorm[]{{\bf a} + {\bf b}}^2 \leq 2\vnorm[]{{\bf a}}^2 + 2\vnorm[]{{\bf b}}^2$ for any ${\bf a},{\bf b} \in \R^N$. For any model center $\bbx \in \R^N$, the growth function is then given by $\omega_{\bbx}(t) = 24\vnorm[]{\bbx}^2 t^2 + 8t^4$.

\section{Model function preserves first order information}\label{sec:first-order-model}

\begin{LEM}\label{lem:first-order-info}
  Let Assumption~\ref{ass:problem}, \ref{ass:model} hold true. For any $\bx \in \dom\fun$, the following condition holds true:
  \[
    \partial_{\by} \fun({\by};\bx)|_{{\by} = \bx} = \widehat{\partial} f(\bx)\,.
  \]
\end{LEM}
\begin{proof}
  We follow the proof strategy of \cite[Lemma 14]{OM2019}. Let $\btx \in \dom\fun$ and let $\bv \in \widehat{\partial} f(\btx)$, then, by definition we have 
  \[
  f(\bx) \geq f(\btx) +  \scal{\bv}{\bx - \btx} + o(\vnorm[]{\bx - \btx}) \quad \forall\, \bx \in \dom\fun.
  \]
  Using the Definition~\ref{def:model-function}, with $f(\btx; \btx) = f(\btx)$ we have the following
  \[
  f(\bx; \btx) + \omega_{\btx}(\vnorm[]{\bx - \btx}) \geq   f(\btx; \btx) +  \scal{\bv}{\bx - \btx} + o(\vnorm[]{\bx - \btx})\,.
  \]
  For any $t>0$, note that $\omega_{\btx}(t) = o(t)$ as $\omega_{\btx}$ is a growth function, using which we obtain 
  \[
  f(\bx; \btx) \geq   f(\btx; \btx) +  \scal{\bv}{\bx - \btx} + o(\vnorm[]{\bx - \btx})\,.
  \]
  This implies that $\bv \in \widehat{\partial} f(\btx; \btx)$ and by regularity of $f(\cdot; \btx)$ we also obtain that $\bv \in \partial f(\btx; \btx)$. For the second part of the proof, let $\bv  \in \partial f(\btx; \btx)$ with $\btx \in \dom\fun$, thus satisfying:
  \begin{align*}
    f(\bbx; \btx) \geq f(\btx; \btx) + \scal{\bv}{\bbx - \btx} + {o(\vnorm[]{\bbx - \btx})}\,,\quad \forall\, \bbx \in \dom\fun\,.
  \end{align*}
  Using the definition of model function (Definition~\ref{def:model-function}), we obtain
  \begin{align*}
      f(\bbx) + \omega_{\btx}(\vnorm[]{\bbx - \btx}) \geq f(\btx; \btx) + \scal{\bv}{\bbx - \btx} + o( \vnorm[]{\bbx - \btx})\,,\quad \forall\, \bbx \in \dom\fun\,,
  \end{align*}
  which on using the fact that $\omega_{\btx}(t) = o(t)$ results in 
  \begin{align*}
    f(\bbx) \geq f(\btx) + \scal{\bv}{\bbx - \btx} + o( \vnorm[]{\bbx - \btx})\,,\quad \forall\, \bbx \in \dom\fun\,.
  \end{align*}
\end{proof}

\section{Proof of  Proposition~\ref{prop:function-descent-property}}\label{sec:function-descent}
By global optimality of $\bx\iter\kp$ as in \eqref{eq:alg-BregMin-bt:update}, we have
\begin{equation}\label{eq:main-descent-func-1}
  \fun(\bx\iter\kp;\bx\iter\k)+  \frac{1}{\tau\iter\k}D_h(\bx\iter\kp,\bx\iter\k) \leq \fun(\bx\iter\k;\bx\iter\k) = \fun(\bx\iter\k)\,.
\end{equation}
We have the following inequality from MAP property
\begin{equation}\label{eq:main-descent-func-2}
  f(\bx\iter\kp) \leq f(\bx\iter\kp;\bx\iter\k) + {\bar L}D_h(\bx\iter\kp,\bx\iter\k)\,.
\end{equation}
Thus, the result follows by combining \eqref{eq:main-descent-func-1} and \eqref{eq:main-descent-func-2}. \qed

\section{Definable functions}\label{sec:definable-functions}
  We require the definition of canonical projection $\Pi : \R^{N+1} \to \R^{N}$  onto $\R^N$, which is defined by 
\[
  \Pi(x_1,\ldots, x_N, t) = (x_1,\ldots,x_N)\,.
\]
\begin{DEF}[o-minimal structure {\cite[Definition 6]{BDLS07}}]
  An o-minimal structure on $(\R, +, .)$ is a sequence of boolean algebras $\mathcal{O}_N$ of ``definable'' subsets of $\R^N$, such that for each $N \in \N$
  \begin{enumerate}
    \item if $A$ in $\mathcal{O}_N$, then $A \times R$ and $R \times A$ belong to $\mathcal{O}_{N+1}$;
    \item if $\Pi: \R^{N+1} \to \R^N$ is the canonical projection onto $\R^N$ then for any $A$ in $\mathcal{O}_{N+1}$, the set $\Pi(A)$ belongs to $\mathcal{O}_{N}$; 
    \item $\mathcal{O}_{N}$ contains a family of algebraic subsets of $\R^N$, that is, every set of the form
    \[
    \{\bx\in \R^N : p(\bx) = 0 \} \,,
    \]
    where $p:\R^N \to \R$ is a polynomial function;
    \item the elements of $\mathcal{O}_1$ are exactly the finite unions of intervals and points.
  \end{enumerate}
\end{DEF}
\begin{DEF}[definable function {\cite[Definition 7]{BDLS07}}]
  Given an o-minimal structure $\mathcal{O}$ (over $(\R,+,.)$), a function $f:\R^N \to \eR$ is said to be definable in $\mathcal{O}$ if its graph belongs to $\mathcal{O}_{N+1}$.
\end{DEF}
\end{document}